\documentclass{svjour3}
\usepackage{geometry}
\usepackage{epsfig}
\usepackage{amsmath}
\usepackage{cuted}
\usepackage{caption}
\usepackage{subfiles}
\usepackage{subcaption}
\usepackage{hyperref}
\usepackage[numbers,sort&compress]{natbib}
\usepackage{pgfplots}
\usepackage{filecontents}
\usepackage{tikz}
\usepackage{tikzscale}
\usepackage{pgf}
\usetikzlibrary{patterns}
\usetikzlibrary{calc}
\usetikzlibrary{decorations.pathmorphing,decorations.markings}
\RequirePackage{fix-cm}
\definecolor{Color_FS}{rgb}{.7,0,.7}
\definecolor{Color_NF}{rgb}{0,0,.7}
\definecolor{Color_NFTO}{rgb}{0,.5,.2}
\definecolor{Color_MD}{rgb}{.75,.35,0}
\definecolor{Color_SMD}{rgb}{1,.7,0}
\newcommand{\q}{R}
\newcommand{\vonkar}{von K\'arm\'an }
\newcommand{\red}[1]{{\color{black} {#1}}}

\begin{document}
\title{Comparison of nonlinear mappings for reduced-order modelling of vibrating structures: normal form theory and quadratic manifold method with modal derivatives}
\titlerunning{Comparison of nonlinear mappings for reduced-order modelling of vibrating structures} 
\author{Alessandra Vizzaccaro\and
	Lo{\"i}c Salles \and
	Cyril Touz\'e
}

\institute{A. Vizzaccaro \at
	Imperial College London, Exhibition Road, SW7 2AZ London, UK\\
	\email{a.vizzaccaro@imperial.ac.uk}
	\and
	L. Salles	\at
	Imperial College London, Exhibition Road, SW7 2AZ London, UK \\
	\and
	C. Touz\'e	\at
	IMSIA, ENSTA Paris, CNRS, EDF, CEA, Institut Polytechnique de Paris, 828 Boulevard des Mar{\'e}chaux, 91762 Palaiseau Cedex, France
}

%


\maketitle

\begin{abstract}
The objective of this contribution is to compare two methods proposed recently in order to build efficient reduced-order models for geometrically nonlinear structures. \red{The first method relies on the normal form theory that allows one to obtain a nonlinear change of coordinates for expressing the reduced-order dynamics in an invariant-based span of the phase space}. The second method is the modal derivative (MD) approach, and more specifically the quadratic manifold defined in order to derive a second-order nonlinear change of coordinates. Both methods share a common point of view, willing to introduce a nonlinear mapping to better define a reduced-order model that could take more properly into account the nonlinear restoring forces. However the calculation methods are different and the quadratic manifold approach has not the invariance property embedded in its definition. Modal derivatives and static modal derivatives are investigated, and their distinctive features in the treatment of the quadratic nonlinearity is underlined. Assuming a slow/fast decomposition allows understanding how the three methods tend to share equivalent properties. While they give proper estimations for flat symmetric structures having a specific shape of nonlinearities and a clear slow/fast decomposition between flexural and in-plane modes, the treatment of the quadratic nonlinearity makes the predictions different in the case of curved structures such as arches and shells. In the more general case, normal form approach appears preferable since it allows correct predictions of a number of important nonlinear features, including for example the hardening/softening behaviour, whatever the relationships between slave and master coordinates are.
\keywords{Reduced Order modelling \and Normal Form \and Quadratic Manifold \and Modal Derivatives}
\end{abstract}

\section{Introduction}

Reduced-order modelling of thin structures experiencing large amplitude vibration is a topic that has attracted a large amount of research in the last years. A number of methods have been proposed, with variants driven either by the structure under study and its peculiarity~\cite{WANG2015}, the dynamical behaviour exhibited by the system~\cite{Weeger2016}, the model~\cite{TOUZE:CMAME:2008} or the discretisation method~\cite{mignolet13}. 

Roughly speaking, one can divide the techniques proposed in the literature into two different categories, the first one using linear change of coordinates, while in the second family nonlinear mappings are defined. When referring to linear methods, one can also distinguish techniques where the best orthogonal basis is computed once and from all. \red{Modal basis~\cite{touze01-JSV,muravyov,Amabmodal,givois2019}}, Ritz vectors~\cite{Kapania1993}, \red{dual modes~\cite{KIM2013}}, and Proper orthogonal decomposition (POD)~\cite{KryslPOD,AmabiliPOD1,Kerschen2005POD} falls into that family. The proper generalized decomposition (PGD)~\cite{Chinesta2011,GroletPGD} under its progressive variant (pPGD) as defined in~\cite{MEYRAND2019} also belongs to that case since additional vectors are added when the dynamics is becoming more complex. On the other hand, the linear change of coordinate can be adaptive, depending on the dynamics, the computation (single point or a whole branch of solution) or the location in phase space. Nonlinear principal component analysis (NLPCA)~\cite{KERSCHEN2002NLPCA} as well as the optimized PGD (oPGD)~\cite{MEYRAND2019} belongs to this family of improved linear methods, sometimes coined as nonlinear since the basis may change depending on some parameter. 

In the third class of methods, a nonlinear change of coordinate is derived once and from all. Nonlinear normal modes \cite{Rosenberg62,ShawPierre93,nayfehnayfeh94,touze03-NNM,KRACK2013}, Spectral submanifolds~\cite{Haller2016,PONSIOEN2018}, and the quadratic manifold derived from modal derivatives~\cite{Jain2017,Rutzmoser}, belongs to this family. As shown in \cite{TOUZE:JFS:2007}, when a linear method ({\em e.g.} POD) tries to find the best orthogonal axis fitting a learning set that have a complex shape, then the number of vectors will be larger than the number of curved subspaces one can use to describe the same datasets. In this particular example, it was shown that invariant manifolds pass exactly through the learning set thus diminishing the number of coordinates needed to describe the dynamics.

Nonlinear normal modes (NNMs) and spectral submanifolds (SSM) offer a rigorously established conceptual framework for reducing geometrically nonlinear structures. In particular, the invariance property of reduction spaces is encapsulated in their definition, ensuring that the dynamical solutions computed from a reduced-order model (ROM) also exist for the full system~\cite{ShawPierre91,Steindl01,TouzeCISM,HallerSF}. This key ingredient allows deriving accurate ROMs, which, for example, are able to predict the correct hardening/softening behaviour of nonlinear structure, which is not the case for their linear counterparts~\cite{touze03-NNM}. \red{More specifically, recent contributions by Haller and collaborators have shown that SSMs are unique continuations of spectral subspaces of the linear system under the nonlinear terms~\cite{Haller2016}, and are thus the best mathematical object to be used in the present context. For nonlinear conservative vibratory systems, SSMs simplify to the classic Lyapunov subcenter manifolds (LSM) that are filled with periodic orbits, thus unifying a number of definitions given for NNMs in the past decades, see {\em e.g.}~\cite{Rosenberg62,Kelley69,VakakisNNM,ShawPierre91}.}

On the other hand, modal derivatives (MDs) have been proposed independently~\cite{IDELSOHN1985,Weeger2016}, and  they share a number of common points with NNMs. In particular, MDs are defined by assuming that the mode shape (eigenvector) together with its eigenfrequency, have a dependence on amplitude, so that one can differentiate the classical Sturm-Liouville eigenvalue problem that defines linear normal modes, in order to make appear a quantity which is defined as the modal derivative. Symmetrically, NNMs also relies on the fact that modal quantities depends on amplitude, a key feature in nonlinear oscillations. The backbone curve and the dependence of the eigenmode shape with amplitude, is then a result from the computation of NNMs, defined as invariant manifold in phase space. However, a complete comparison of both method has not been drawn out yet. The only related paper  uses the modal derivatives as a reduction method, from which the NNM, seen in this case as the family of periodic orbit in phase space --and thus reducing their information to the backbone curve only, without using the geometrical information in phase space-- can be computed~\cite{SOMBROEK2018}.

A recent development in the use of modal derivatives is to form a quadratic manifold for more accurate model order reduction. The properties of this nonlinear mapping are such that it is tangent to a subspace spanned by the most relevant vibration modes, and its curvature is provided by modal derivatives~\cite{Jain2017}. An idea also claimed in~\cite{Rutzmoser} is that such a quadratic manifold should be able to cancel the quadratic forces in the ROM. Incidentally, NNMs defined in the framework of normal form theory, as proposed in~\cite{touze03-NNM,TOUZE:JSV:2006}, already present these features. Indeed, a third-order nonlinear change of coordinate is given, which has the property to be identity-tangent when the initial model is expressed in modal coordinates, thus conserving the linear modes as first approximation. Also, in case of no second-order internal resonance, the mapping exactly cancels all quadratic terms. Finally, the invariance property  is directly inherited from the definition of an NNM as an invariant manifold in phase space, while the invariance of the quadratic manifold computed from MDs is not at hand. 

The aim of this contribution is thus to investigate more properly the common points
and differences of the two methods, and explain their advantages and drawbacks
in the context of building reduced-order models for geometrically nonlinear structures.
The paper is organized as follows. Section 2 is concerned with the theoretical developments. The framework of geometrically nonlinear structures is briefly recalled, then both methods of interest, normal form theory, modal derivatives and their extension to quadratic manifold (QM), are recalled and analysed in depth. The general derivation of the QM framework for both modal derivatives (MDs) and static modal derivatives (SMDs) is highlighted, whereas previous contributions generally use the simplifying assumption of SMDs in the developments. As a consequence of this development, the distinctive treatment of the quadratic nonlinearity between MDs and SMDs is specifically underlined. Of particular interest is the comparison of methods when a slow/fast decomposition of the system can be assumed. In the course of the paper, we will contrast the results given by MDs, SMDs and normal form and underlines that the simplifying assumption of slow/fast approximation allows retrieving partly the correct results.
By doing so, an illustration of the general theorem given in~\cite{HallerSF} is thus provided for a more restrictive framework. Indeed, theorems given in~\cite{HallerSF} encompasses more generality and exact results, allowing to deal with the case of damping and forcing. We give however here more detailed comparisons, and in particular analyse how the SMD can produce incorrect predictions for structures having a strong quadratic coupling such as arches and shells. Section 3 illustrates the findings of the previous section on two simple two degrees-of-freedom (dofs) systems. Finally section 4 applies the previous results to continuous structures discretised with the finite element (FE) procedure.

\section{Models and methods}

\subsection{Framework}

Geometric nonlinearity refers to the case of thin structures vibrating with large amplitudes while the material behaviour remains linear elastic. In this framework, the semi-discretised version of the equations of motion, generally obtained from a finite-element procedure, reads~:
\begin{equation}
\vec{M}\vec{\ddot{u}} + \vec{F}(\vec{u})=\vec{Q},
\label{eq:first}
\end{equation}
where $\vec{M}$ is the mass matrix, $\vec{u}$ the displacement vector at the nodes, $\vec{F}$ the nonlinear restoring force and $\vec{Q}$ the external force. The number of degrees of freedom (dofs) is $N$, being thus the dimension of vectors $\vec{u}$, $\vec{F}$  and $\vec{Q}$. Note that damping is presently not taken into account since most of the presented work deals with efficient treatments of nonlinearities in the restoring force. While the concepts of NNMs and spectral submanifolds (SSM) can be straightforwardly extended to the cases with damping, as already shown for example in~\cite{TOUZE:JSV:2006} for normal form or in~\cite{Haller2016} for SSM, a clear extension of MDs to damped systems does not seem to be present in the literature, to the best of our knowledge. \red{Consequently, we restrict ourselves in this contribution to the treatment of the nonlinear stiffness without considering the effect of damping, but we acknowledge that damping have important effects in nonlinear vibrations that should thus need further investigations.}
 
Geometric nonlinearity for slender structures is assumed so that $\vec{F}$, for the sake of simplicity, only depends on the displacement vector $\vec{u}$, but other cases can also be treated. 
More particularly, a number of models have been derived for thin structures such as plates and shells, relying on simplifying assumptions ({\em e.g.} \vonkar models for beams and plates~\cite{Landau1986,Bazant,ThomasBilbao08}, Donnell's assumption for shallow shells \cite{Amabilibook,AmabiliLagrange}), showing that the partial differential equations of motion only contains quadratic and cubic terms with respect to the displacement. On the other hand, general equations for three-dimensional elasticity with geometric nonlinearity (linear stress/strain relationship but nonlinear strain/displacement relationship) also show that the nonlinear terms in the restoring force should be of this type~\cite{mignolet13,LazarusThomas2012,Touze:compmech:2014,givois2019}. Consequently  we consider in this contribution a nonlinear force that can be expressed as a function of the displacement up to cubic order terms, reading:
\begin{equation}\label{eq:restoforceexp}
\vec{F}(\vec{u})=\vec{K}\vec{u}+ \vec{G}\vec{u}\vec{u} + \vec{H}\vec{u}\vec{u}\vec{u}.
\end{equation}
In this last equation, we use a simplified notation of the tensor product for the quadratic and cubic terms, already introduced in~\cite{Jain2017,Rutzmoser}. The notation is fully explained in Appendix~\ref{app:products}, where the indicial expressions of the products are detailed for the sake of clarity.  $\vec{G}$ is a third-order tensor of quadratic coefficients with current term $G^p_{ij}$, while $\vec{H}$ is the fourth-order tensor grouping the cubic coefficients $H^p_{ijk}$. For example, the vector $\vec{G}\vec{u}\vec{u}$ of the quadratic terms writes:
\begin{equation}\label{eq:product_GuuText}
\vec{G}\vec{u}\vec{u}=\sum_{i=1}^N \sum_{j=1}^N \vec{G}_{ij}{u}_{i}{u}_{j},
\end{equation}
with $\vec{G}_{ij}$ the N-dimensional vector of coefficients $G^p_{ij}$, for $p=1,\, ...,\, N$. 
Note also that in this contribution, 
the representation of quadratic and cubic terms does not use the fact that the usual product is commutative 
($u_{i}u_{j}=u_{j}u_{i}$), 
so that the second summation in~\eqref{eq:product_GuuText} 
could be limited to the indices such as  $j\geq i$, 
assuming also $\vec{G}_{ij}=0$ 
for $i\geq j$. 
In the representation selected throughout the paper, all summations will be full, as in~\eqref{eq:product_GuuText} with a fully populated tensor of coefficient  $\vec{G}$. The same rule applies for the cubic term also. 
This choice has been made since it allows shorter and simpler expressions for a number of equations given in the presentation, but of course it is not a limiting assumption and the other choice could have also be done.

The first (linear) term in Eq.~\eqref{eq:restoforceexp} makes appear the usual tangent stiffness matrix $\vec{K}$ defined by~:
\begin{equation}
\vec{K} = \left. \frac{\partial \vec{F}}{\partial \vec{u}}\right|_{\vec{u}=0},
\end{equation}
from which one can define the eigenmodes, solution of the eigenvalue problem:
\begin{equation}
(\vec{K} - \omega_i^2 \vec{M})\vec{\phi}_i = 0,
\label{linear_eigsys}
\end{equation}
with $\vec{\phi}_i$ the $i^{\mathrm{th}}$ eigenvector and $\omega_i$ its associated eigenfrequency, for $i=1,\, ...,\, N$.  Using $\vec{u} = \vec{\Phi}\vec{X}$, with $\vec{\Phi}$ the matrix of all eigenvectors $\vec{\phi}_i$, and $\vec{X}$ the modal coordinates, the problem can be rewritten in the modal basis by premultiplying 
Eq.~\eqref{eq:first} by $\vec{\Phi}^T$, arriving at:
\begin{equation}
\vec{\ddot{X}} + \vec{\Omega}^2 \vec{X} + \vec{g}\vec{X}\vec{X}+\vec{h}\vec{X}\vec{X}\vec{X}=\vec{q},
\end{equation}
where we have introduced $\vec{\Omega}$ the matrix of eigenfrequencies $\omega_i$,  $\vec{g}$ and 
$\vec{h}$ the tensors of quadratic and cubic coefficients in the modal basis, and $\vec{q} = \vec{\Phi}^T \vec{Q}$ the modal external force. The equation of motion in modal space can be written in explicit form with these coefficients as: 
\begin{equation}
 \forall \, p=1, \, ..., N  \; :\quad  
 \ddot{X}_p + \omega_p^2 X_p 
 + 
 \sum_{i=1}^{N} \sum_{j =1}^{N} 
 g_{ij}^p X_i X_j 
 +
  \sum_{i=1}^{N} \sum_{j=1}^{N} \sum_{k =1}^{N}  
  h_{ijk}^p X_iX_jX_k = q_p.
\label{eq:EDO}
\end{equation}
The relations between the nonlinear tensors in physical coordinates $\vec{G}$ and $\vec{H}$, and those in modal coordinates $\vec{g}$ and $\vec{h}$ are derived from the linear change of coordinates and involves products with the matrix of eigenvectors $\vec{\Phi}$. They are  provided in Appendix~\ref{app:phys_to_mod} for the sake of completeness.

\subsection{NNMs and normal form}\label{subsec:NNMNF}

Nonlinear normal modes or NNMs have been used since the pioneering work by Rosenberg~\cite{Rosenberg62} in numerous vibratory problems. It offers a sound conceptual framework in order to understand the organization of the dynamics in the phase space. \red{Different definitions  have been given in the past, {\em e.g.} family of periodic orbits~\cite{Rosenberg62,KerschenNNM09}, invariant manifold in phase space, tangent at the linear eigenspaces near the origin~\cite{ShawPierre93}. More recently, a mathematically well-justified definition of NNM has been provided~\cite{Haller2016}, allowing to settle down the different treatments in an unified way. For that purpose, Haller and Ponsioen proposed to refer to the smoothest member of an invariant manifold family tangent to a modal subbundle along an NNM as a spectral submanifold (SSM). In that sense, SSMs provides a rigorous framework allowing to define the corresponding concepts in all the situations encountered in mechanical vibrations: conservative or dissipative systems, autonomous or non-autonomous systems. Interestingly, the authors also provide in \cite{PONSIOEN2018} automated formulations in order to derive SSMs up to large order, allowing them to draw out comparisons with numerous other methods proposed in the recent years, see {\em e.g.}~\cite{BreunungHaller18}.  Enforcing the invariant property  is  key in a perspective of reduced-order modelling, since it is the only way to ensure that the trajectories of the ROM will also exist for the full system. Elaborating on this idea, NNMs has been used in the perspective of model-order reduction using either center manifold theorem~\cite{ShawPierre93,PesheckJSV}, normal form theory~\cite{touze03-NNM,TOUZE:CMAME:2008,TouzeCISM}, or spectral submanifolds~\cite{Haller2016,PONSIOEN2018,BreunungHaller18,VERASZTO}.

In this contribution, the normal form theory,  as defined in~\cite{touze03-NNM,TOUZE:JSV:2006}, is used. The main idea is to define a nonlinear change of coordinates, from the modal coordinates to new ones defined as the {\em normal} coordinates. The nonlinear mapping is inherited from Poincar{\'e} and Poincar{\'e}-Dulac theorems, based on the idea of finding out a nonlinear relationship capable of eliminating as much as possible of nonlinear terms.} In this contribution, only the main results are recalled, the interested reader is referred to~\cite{touze03-NNM,TOUZE:JSV:2006,TouzeCISM} for more details. The nonlinear change of coordinates is identity-tangent, and formally reads:
\begin{subequations}
\label{eq:zecvNL}
\begin{align}
X_p & = R_p + {\mathcal P}_p (R_i,S_j), \label{eq:zecvNL-a}\\ 
Y_p & = S_p + {\mathcal Q}_p (R_i,S_j), \label{eq:zecvNL-b} 
\end{align}
\end{subequations}
where ${\mathcal P}_p$ and ${\mathcal Q}_p$ are third-order polynomials, the analytical expressions of which are given in~\cite{touze03-NNM} for the undamped case and in \cite{TOUZE:JSV:2006} for the damped case. $X_p$ is the modal coordinate, $Y_p$ the modal velocity, and $(R_i,S_j)$ are the new coordinates related to the invariant manifolds, and called {\em normal} coordinates.

The method used to derive the nonlinear mapping is based on the recognition of nonlinear resonances involving the eigenfrequencies of the system. In case where no internal resonance is present, one can show for example that all the quadratic terms can be cancelled from the normal form which is thus much simpler than the original system. 

The dynamics, expressed with the newly introduced normal variables $(R_i,S_j)$, is written in an invariant-based span of the phase space, and thus prone to open the doors to efficient reduced-order models, as already shown in \cite{TouzeCISM}. The general equation for the dynamics expressed in the new coordinates reads:
\begin{align}\label{normadynr}
 \forall \; p \, = & \, 1, \, . .  . \, , \, n: \nonumber \\
\ddot{R}_p+&\omega_p^2 R_p +
(A_{ppp}^p + h_{ppp}^p) R_p^3 + B_{ppp}^p R_p S_p^2 + \nonumber \\
& + R_p \underset{j\neq p}{\sum_{j=1}^n}
\left((3\,h_{pjj}^p+2\,A_{jjp}^p + A_{pjj}^p)R_j^2 + B_{pjj}^p S_j^2 \right)
+ S_p \underset{j\neq p}{\sum_{j=1}^{n}}
\left(2\,B_{jjp}^p R_j S_j\right)
 =0.
\end{align}
where $n$ is the number of master modes retained for the ROMS, $\vec{\q} = (\q_1, ..., \q_n)$; in most cases $n\ll N$, but the formula are given for $n$ arbitrary and can be used also for $n=N$.
Note that the expression in slightly different from the one proposed in~\cite{touze03-NNM}, a direct consequence of the choice of the representation of quadratic and cubic terms, with full summations. 
The coefficients $A_{ijk}^p$ and $B_{ijk}^p$ stems from the cancellation of the quadratic terms. Their expressions read:
\begin{subequations} \label{AB}
\begin{gather}
A_{ijk}^p = \sum_{s =1}^{N} 2\;\bar{g}_{is}^p a_{jk}^s, \label{AB1} \\
B_{ijk}^p = \sum_{s =1}^{N} 2\;\bar{g}_{is}^p b_{jk}^s,\label{AB2}
\end{gather}
\end{subequations}
where $\bar{g}_{is}^p = (  g^p_{is}  +  g^p_{si}  )/2$ is the mean value between two adjacent terms implying the same monomial term. The coefficients $a_{jk}^s$ and $b_{jk}^s$ appearing in the expression of $A_{ijk}^p$ and $B_{ijk}^p$ are related to the quadratic terms of the change of coordinate. \red{For the sake of completeness, the interested reader can find their full expressions in Appendix~\ref{app:NF_coeff}.}
As known from the theory, these second-order coefficients have a singular behaviour in the vicinity of internal resonances. In this case, a strong coupling is present between the nonlinear oscillators whose eigenfrequencies are commensurate, and the associated coefficient in the change of coordinate is set to zero, so that the corresponding monomial terms stay in the normal form.

From Eqs.~\eqref{normadynr}, one can observe that invariant-breaking terms are no longer present in the equations of motion. Invariant-breaking terms are defined as quadratic monomials of the form $g^k_{pp}X_p^2$ and cubic monomials $h^k_{ppp}X_p^3$ on $k$-th oscillator equation. As soon as mode $p$ has some energy, then these invariant-breaking terms directly excite oscillator $k$, thus breaking the invariance of the linear mode subspace. As these terms are no longer present in Eqs.~\eqref{normadynr}, it shows that the dynamics is now expressed in an invariant-based span. One can also note that the only monomial terms present in Eqs.~\eqref{normadynr} are those related to trivially resonant terms.

A ROM is simply selected by keeping in the truncation only the normal coordinates $(R_p,S_p)$ of interest, depending on the problem at hand. By doing so, one restricts the motion in the invariant manifold described by the master normal coordinates retained, giving rise to efficient reduced models, that simulate trajectories existing in the complete phase space, and allowing to recover the correct type of nonlinearity~\cite{touze03-NNM,touze-shelltypeNL} as well as nonlinear frequency response curves~\cite{TOUZE:CMAME:2008}. The simplest ROM is built by restricting the motion to a single NNM by keeping only one pair  $(R_p,S_p)$ and cancelling all the other: $\forall \; k\neq p, \quad R_k = S_k = 0$. In this case the nonlinear change of coordinates for the master coordinates reads:
\begin{subequations}\label{eq:transformNF}
\begin{align}
X_p &= R_p+a_{pp}^p R_p^2 + b_{pp}^p S_p^2,  \label{eq:transformXp} \\
Y_p &= S_p + \gamma_{pp}^p R_p S_p, \label{eq:transformYp}
\end{align}
\end{subequations}
whereas for the slave coordinates one has:
\begin{subequations}
\begin{align}
\forall \, k \neq & p:  \nonumber \\
X_k &= a_{pp}^k R_p^2 + b_{pp}^k S_p^2+ r_{ppp}^k R_p^3 + u_{ppp}^k R_pS_p^2,  \label{eq:transformXk} \\
Y_k &= \gamma_{pp}^k R_pS_p+ \mu_{ppp}^k S_p^3 + \nu_{ppp}^k S_pR_p^2. \label{eq:transformYk}
\end{align}\label{eq:invarman}
\end{subequations}
Again, all the introduced coefficients, $\gamma_{pp}^p$, $r_{ppp}^k $, $u_{ppp}^k $, $\mu_{ppp}^k $ and $\nu_{ppp}^k$, originate from the explicit expression of the polynomials ${\mathcal P}_p$ and ${\mathcal Q}_p$ of Eq.~\eqref{eq:zecvNL}. They are all analytic and their expressions are given in~\cite{touze03-NNM}. \red{Interestingly, Eqs.~\eqref{eq:invarman} describes the geometry of the invariant manifold in phase space, up to order three, but of course one can limit the development of this equation to second-order only.}  The dynamics on the invariant manifold ($p^{\mbox{th}}$ NNM) is found by cancelling all $(R_k,S_k)$ for $k\neq p$ in Eqs.~\eqref{normadynr}. 
In the case of a single NNM motion the equation is particularly simple and reads:
\begin{equation} 
\ddot{R}_p + \omega_p^2 R_p + (A_{ppp}^p + h_{ppp}^p) R_p^3 + B_{ppp}^p R_p \dot{R}_p^2 \; = \; 0 \; . 
\label{dynsingledof}
\end{equation}
\red{Of particular interest here is the fact that the correcting coefficients $A_{ppp}^p$ and $B_{ppp}^p$ appearing in this last equation are provided by the second-order terms in the nonlinear change of coordinates. Consequently, the third-order terms have no influence on this reduced dynamics, which is thus exactly the one given by the second-order truncation of the normal form nonlinear mapping.
}

All these formulas can be used to reconstruct the mode shape dependence on amplitude, assuming the motion is enslaved to a single NNM, {\em i.e.} lying in the invariant manifold associated to mode $p$. Assuming this single-NNM motion, the physical displacement is reconstructed from
\begin{equation}
\vec{u} = \sum_{k=1}^N X_k \vec{\phi}_k = X_p \vec{\phi}_p + \underset{k\neq p}{\sum_{k=1}^N}X_k \vec{\phi}_k,
\end{equation}
where $X_p$ is replaced using Eq.~\eqref{eq:transformXp} and $X_k$ using~Eq.~\eqref{eq:transformXk}, so that one finally obtains the amplitude-dependent mode shape as~:
\begin{equation}
\vec{u} = \left( R_p+a_{pp}^p R_p^2 + b_{pp}^p S_p^2 \right) \vec{\phi}_p + \underset{k\neq p}{\sum_{k=1}^N} \left(a_{pp}^k R_p^2 + b_{pp}^k S_p^2+ r_{ppp}^k R_p^3 + u_{ppp}^k R_pS_p^2\right) \vec{\phi}_k.
\end{equation}
This formula has already been used in order to represent the amplitude dependence of mode shapes on amplitude, see {\em e.g.}~\cite{ShawPierre94,touze03-NNM}, and will be further analysed and compared to the prediction given by the method of quadratic manifold from modal derivatives in Sect.~\ref{sec:drift}.

\red{Note that, as a comparison to quadratic manifold is targeted, a detailed description of the effects of order truncation in the normal form approach is in order. In the present approach of the normal form, the change of coordinates is up to order three, but the reduced-order dynamics can be considered as up to the second order, since the effect of cancelling the cubic terms to the higher-orders have not been taken into account due to the third-order truncation of all asymptotic developments. Also, most of the comparisons in the remainder of the paper will be drawn between single-mode reduced-order dynamics. In this simplified context, Eq.~\eqref{dynsingledof} clearly shows that the cancellation of the third-order non-resonant monomials have absolutely no effect on this equation which is left unchanged. Consequently, 
Eq.~\eqref{dynsingledof} is the reduced dynamics obtained with a second-order normal form nonlinear mapping. The only difference between second-order and third-order is in Eq.~\eqref{eq:invarman}, which describes how the exact invariant manifold is approximated in phase space, and one can analyse the effect of either second-order or third-order nonlinear mapping in this respect. In the remainder of the paper, a clear attention will be devoted to these two specific truncations in order to draw out a fair comparison with the quadratic manifold approach.}

We now turn to the definition of modal derivatives and the associated nonlinear mapping: the so-called quadratic manifold, before comparing the two methods in detail.

\subsection{Modal Derivatives}

Modal derivatives have been first introduced by Idehlson and Cardona to solve structural vibrations problems with a nonlinear stiffness matrix~\cite{IDELSOHN1985}. They have been used in recent years in the context of reduced-order modelling~\cite{Weeger2016}, and the last developments propose to use them in order to create a nonlinear mapping with a quadratic manifold~\cite{Jain2017,Rutzmoser}. In this section, we derive again the most important definitions, make the distinction between modal derivatives (MDs) and static modal derivatives (SMDs), and introduce the quadratic manifold approach.

\subsubsection{Definition of Modal Derivatives and Static Modal Derivatives}\label{subsec:defMD}

The modal derivatives have been first introduced with the aim of offering a framework taking into account the dependence of mode shapes and eigenfrequencies on amplitude for nonlinear system. This is a common point with nonlinear normal modes, that also recognizes this fact as a major outcome that needs to be addressed correctly in the modelling. \red{The introduction of the modal derivatives proposed in this section is mostly heuristic and based on previous works.} Let us denote $\tilde{\vec{\phi}}_i(\vec{u})$ this amplitude-dependent eigenvector.  The already introduced eigenvector $\vec{\phi}_i$,  solution of the Sturm-Liouville problem, Eq.~\eqref{linear_eigsys}, represents the value of $\tilde{\vec{\phi}}_i(\vec{u})$ when $\vec{u}=\vec{0}$. The  $ij$-th modal derivative (MD) is defined as the derivative of  $\tilde{\vec{\phi}}_i$ with respect to the $j$-th coordinate used for the reduced basis, denoted here as $R_j$. For the sake of clarity,  $X_i$ is the modal coordinates, and $R_j$ the reduced coordinates, following the notations introduced for the normal form approach. At first order, one has $X_i=R_i$, but as we consider nonlinear change of coordinates, these relationships will be enriched by higher-order terms. For the quadratic manifold approach, this will be explained in the next subsections, so that for the present definitions, one can assume $R_i=X_i$. In that context, the $ij$-th MD $\vec{\Theta}_{ij}$ is the derivative of $\tilde{\vec{\phi}}_i$ with respect to a displacement enforced along the direction of the $j$-th eigenvector $\vec{\phi}_j$ as introduced in~\cite{IDELSOHN1985,Weeger2016,Jain2017,Rutzmoser}, and writes:
\begin{equation}
\vec{\Theta}_{ij} 
\doteq \frac{\partial \tilde{\vec{\phi}}_i(\vec{u})}{\partial \q_j}
\bigg\rvert_{\vec{u}=\vec{0}} \, .
\end{equation}
In order to derive an equation from which the MD can be computed, one has to rewrite the eigenproblem given by Eq.~\eqref{linear_eigsys} assuming the known dependencies on the amplitude, as:
\begin{equation}
\left(\dfrac{\partial\vec{F}(\vec{u})}{\partial \vec{u}} - \tilde{\omega}_i^2(\vec{u}) \vec{M}\right)\tilde{\vec{\phi}}_i(\vec{u}) = \vec{0},
\label{nonlinear_eigprob}
\end{equation}
where the linear stiffness matrix is replaced by the full nonlinear restoring force,  and both eigenvalues and eigenvectors are amplitude-dependent. Note that, in this contribution, the mass matrix is assumed to be independent of the amplitude, since this is the selected framework for this paper focused on geometric nonlinearity. However further development could include a dependence of the mass matrix on the amplitude in order to extend the use of MDs to other cases. The nonlinear eigenproblem of Eq.~\eqref{nonlinear_eigprob} must be complemented with the nonlinear mass normalisation equation:
\begin{equation}
\tilde{\vec{\phi}}_i(\vec{u})^T 
\vec{M}
\tilde{\vec{\phi}}_i(\vec{u})=1.
\label{nonlinear_mass_norm}
\end{equation}
The last two equations, \eqref{nonlinear_eigprob}-\eqref{nonlinear_mass_norm} can then be Taylor-expanded as function of the amplitude, assuming moderate vibrations in the vicinity of the position at rest defined by $\vec{u}=\vec{0}$. Assuming that the displacement $\vec{u}$ depends on the coordinates introduced for the reduced basis, $R_1$ to $R_n$, each term can then be expanded along these new coordinates. The full derivation of this Taylor expansion is given in Appendix~\ref{app:taylor_MD}.

The Taylor expansion of Eq.~\eqref{nonlinear_eigprob} and Eq.~\eqref{nonlinear_mass_norm} in the $\q_j$  coordinates, up to first order, generates constant terms that coincide with the linear eigenproblem and mass normalisation. The next order terms, linear in $\q_j$, allows deriving the following system, where the two unknowns are the MD vector $\vec{\Theta}_{ij}$, and the scalar describing the variation of the squared eigenfrequency with respect to amplitude, $\frac{\partial \omega_i^2}{\partial \q_j}$:
\begin{equation}
\begin{bmatrix}
\vec{K}-\omega_i^2 \vec{M}
&\;
-\vec{M}
\vec{\phi}_i
\\[10pt]
-\vec{\phi}_i^T 
\vec{M}
&
\; 0
\end{bmatrix}
\left\lbrace
\begin{matrix}
\vec{\Theta}_{ij}
\\[10pt]
\frac{\partial \omega_i^2}{\partial \q_j}
\end{matrix}\right\rbrace
=
\left\lbrace
\begin{matrix}
-2\vec{G}
\vec{\phi}_j \vec{\phi}_i
\\[15pt]
0
\end{matrix}\right\rbrace,
\label{eq:MD_sys}
\end{equation}
where the quadratic tensor $\vec{G}$ of the restoring force introduced in Eq.~\eqref{eq:restoforceexp}, has been used. The detailed proof for the derivation of this system is given in Appendix~\ref{app:taylor_MD}.

In most of the studies concerned with application of modal derivatives to model order reduction, the so-called {\em static modal derivatives} (SMDs) are used instead. Let us denote as $\vec{\Theta}^{(S)}_{ij}$
the SMD of $\vec{\Theta}_{ij}$, obtained by neglecting the terms related to the mass matrix in~\eqref{eq:MD_sys}, which then simplifies to:
\begin{equation}
\vec{K}\vec{\Theta}^{(S)}_{ij}=-2\vec{G}
\vec{\phi}_j \vec{\phi}_i.
\label{eq:SMD_sys}
\end{equation}
This last simplification evidently highlights the fact that MDs and SMDs are able to retrieve the quadratic coupling generated by the nonlinear restoring force, since being directly proportional to the tensor of coefficients $\vec{G}$. Eq.~\eqref{eq:SMD_sys} also shows that the computation of SMDs is drastically reduced as compared to MDs, for two main reasons. The first one is that, given the usual symmetry of the quadratic tensor $\vec{G}$ at hand in structural problems, one has $\vec{G}
\vec{\phi}_j \vec{\phi}_i=\vec{G}
\vec{\phi}_i \vec{\phi}_j$, so that the SMDs are symmetric $\vec{\Theta}^{(S)}_{ij}=\vec{\Theta}^{(S)}_{ji}$. This involves that the number of calculations for indexes $i\neq j$ is then halved in the case of SMDs as compared to MDs. 
The second reason lies in the fact that, despite the sizes of the systems to solve are comparable (the size of system \eqref{eq:MD_sys} is $N+1$ and the size of system \eqref{eq:SMD_sys} is $N$), the computation of a SMD
can be done with a standard operation in a commercial FE software whereas the computation of a MD cannot.
Indeed, the non-intrusive computation of a SMD requires to solve a linear system $\vec{K}\vec{u} = \vec{f}$, where the applied force $\vec{f}$ is the right-hand side of Eq.~\eqref{eq:SMD_sys} and the resultant displacement $\vec{u}$ is the SMD. Solving such linear system coincides with operating a simple linear static analysis on the structure with imposed force and unknown displacement. Conversely, the linear system to compute a MD is the one in Eq.~\eqref{eq:MD_sys}. The solution of this system does not correspond to the standard operation one could easily perform in a FE software. Consequently to compute the MD, one needs not only to access to the full stiffness and mass matrices but also to export them in an external code to be able to solve the linear system. When the structure is discretised with a large number of dofs, such operation can be memory and time consuming when not infeasible.

\subsubsection{Expression of MDs as function of the quadratic coefficients from the modal basis}\label{sec:MDmodal}

In this section, the relation between MDs and SMDs and the coefficients of the quadratic tensor in modal basis $\vec{g}$ is derived. This relation will help to draw comparisons between the normal form method and the quadratic manifold method that will be introduced in the next section. For that purpose, the $ij$-th MD in the modal basis, denoted as $\vec{\theta}_{ij}$, is introduced as
\begin{equation}\label{eq:MDmodalchange}
\quad \vec{\Theta}_{ij}=\vec{\Phi}\vec{\theta}_{ij}=\sum_{s=1}^N \vec{\phi}_s \theta^s_{ij},
\end{equation}
following the linear change of basis from physical to modal space, where the summation thus spans over all the modes of the structure, being $\vec{\Phi}$ the full eigenvector matrix. In the modal basis, the eigenvector $\vec{\phi}_i$ coincides with the $i$-th vector basis $\vec{e}_i$, where the entries of $\vec{e}_i$ are all zero except 1 in position $i$, so that: $\vec{\phi}_i=\vec{\Phi}\vec{e}_i$.

The system of equations \eqref{eq:MD_sys}, can be now written in modal coordinates by premultiplying the first $N$ rows by $\vec{\Phi}^T$ and by substituting the values of $\vec{\phi}_i$ and $\vec{\Theta}_{ij}$ with their values in modal coordinates. One finally  obtains:
\begin{equation}
\begin{bmatrix}
\vec{\Omega}^2-\omega_i^2 \vec{I}
&\;
-
\vec{e}_i
\\[10pt]
-\vec{e}_i^T 
&
\; 0
\end{bmatrix}
\left\lbrace
\begin{matrix}
\vec{\theta}_{ij}
\\[10pt]
\frac{\partial \omega_i^2}{\partial \q_j}
\end{matrix}\right\rbrace
=
\left\lbrace
\begin{matrix}
-2\vec{g}_{ij}
\\[15pt]
0
\end{matrix}\right\rbrace,
\label{eq:MD_sys_mod}
\end{equation}

where the right-hand side has been simplified using the relationship $\vec{g}_{ij}=\vec{\Phi}^T \vec{G} \vec{\phi}_i \vec{\phi}_j$, demonstrated in Eq.~\eqref{eq:g_ij} of Appendix~\ref{app:phys_to_mod}. 

The system \eqref{eq:MD_sys_mod} is easier to understand when written term by term:
\begin{subequations}\label{eq:mdmodalexpr}
\begin{align}
& (\omega_s^2-\omega_i^2)\theta^s_{ij}=-2g^s_{ij},\qquad \mbox{for} \quad s\neq i,\label{eq:mdmodalexpr-a}\\
& \frac{\partial \omega_i^2}{\partial \q_j}=2g^i_{ij}, \qquad \mbox{for} \quad s=i, \label{eq:mdmodalexpr-b}\\
&\theta^i_{ij}=0.\label{eq:mdmodalexpr-c}
\end{align}
\end{subequations}

One can notice that the $ij$-th modal derivatives is then directly proportional to the $ij$-th component of the quadratic tensor in modal coordinates. This clearly shows that the $ij$-th MD is able to retrieve a strong quadratic coupling occurring between slave mode $s$ and the  master modes $i$ and $j$.   The value of the modal derivative in physical coordinates can be now easily reconstructed from the preceding development, and reads:
\begin{equation}
\vec{\Theta}_{ij}
=
\underset{s\neq i}{\sum_{s=1}^N} \vec{\phi}_s
\;\dfrac{-2\;g^s_{ij}}{\omega_s^2-\omega_i^2}.
\label{eq:MD}
\end{equation}
If one follows a similar procedure for the case of static modal derivative, Eq.~\eqref{eq:SMD_sys} is written in modal coordinates as $
\vec{\Omega}^2\vec{\theta}^{(S)}_{ij}=-2\vec{g}_{ij}$
and the static modal derivative in physical coordinates is directly given as:
\begin{equation}
\vec{\Theta}^{(S)}_{ij}=\sum_{s=1}^N
\vec{\phi}_s
\;\dfrac{-2\;g^s_{ij}}{\omega_s^2}.
\label{eq:SMD}
\end{equation}
In both cases, MDs and SMDs can be simply defined as a linear combination of modes weighted by a factor proportional to $g^s_{ij}$, the quadratic modal coupling coefficient. In the case of modal derivative, the method shows a divergent behaviour in case of 1:1 internal resonance between two eigenfrequencies, a feature that will be further commented in Sect.~\ref{sec:HS}. One can also note that the weighting factors have larger values for the modes, the eigenfrequencies of which are closer to the eigenfrequency of the $i$-th mode. On the other hand for static modal derivatives, the weighting factors are simply proportional to the inverse of the squared eigenfrequencies, and thus should decrease for higher modes. Note however that this fact can be severely compensated by the values of the quadratic coefficients, which scales according to the linear stiffness. Consequently, as shown for example in~\cite{SOMBROEK2018,Vizza3d} for thin and flat symmetric structures (beams and plates), the SMD is able to recover the most important couplings with in-plane modes. 

As a conclusion, MDs and SMDs can be seen as a displacement field that takes into account the contribution of all \textit{quadratically coupled} modes into one equivalent vector. From this perspective, the use of a reduced basis composed of MDs is equivalent to using a basis composed of all quadratically coupled modes, with the supplementary condition that the quadratic couplings makes appear new directions in phase space, that are independent of the already selected mode. If the quadratic coupling is only dependent on modes already present in the reduced basis, then the new vector will not bring out new eigendirections.

\subsubsection{Quadratic manifold}

The quadratic manifold approach has been introduced in~\cite{Jain2017,Rutzmoser} in order to extend the use of modal derivatives in the context of model order reduction, and propose a nonlinear mapping from initial to reduced coordinates. The nonlinear mapping is quadratic in nature and does not account for nonlinear internal resonance as the normal form theory does. In this section, the derivation of reduced-order models using the quadratic manifold is given, following the previous results obtained in~\cite{Jain2017,Rutzmoser}. A particular attention is paid on writing the differences one can await when using the quadratic manifold with MDs and SMDs, with the comparison to the results provided by normal form theory in mind, thus giving rise to new developments. The coordinates describing the reduced-order models are denoted as $R_p$ for all the methods in order to compare more directly the equations. One has however to keep in mind that the meaning of these coordinates is not the same for each method.

Since the MDs are defined from a second-order Taylor expansion of the nonlinear eigenvalue problem, it is intuitive to use them in a quadratic nonlinear mapping. If one operates a Taylor expansion of the approximate solution $\vec{u}$ in the reduced coordinates $\vec{R}$ up to quadratic order, one finds:
\begin{equation}\label{eq:TayloruR}
\vec{u}(\vec{\q})=\vec{u}(\vec{0})
+
\sum^n_{i=1}\dfrac{\partial \vec{u}(\vec{\q})}{\partial \q_i}\bigg\rvert_{\vec{0}}\q_i
+
\frac{1}{2}
\sum^n_{i=1}\sum^n_{j=1}
\left(
\dfrac{\partial^2 \vec{u}(\vec{\q})}{\partial \q_j \partial \q_i}
\right)\bigg\rvert_{\vec{0}}\q_i\q_j
+
\mathcal{O}(|\vec{\q}|^3),
\end{equation}
where $n$ is the number of master modes retained for the ROMS, $\vec{\q} = (R_1, ..., R_n)$. By extending the definition of linear eigenvectors to the nonlinear ones, the nonlinear eigenvector spans the tangent space of the displacement with respect to the reduced coordinates, so that:
\begin{equation}
\dfrac{\partial\vec{u}}{\partial R_i} =  \tilde{\vec{\phi}}_i (\vec{R}).
\end{equation}

In Eq.~\eqref{eq:TayloruR}, we can then substitute $\vec{u}(\vec{0}) = \vec{0}$ (the position at rest is at the origin of the coordinates), and
\begin{align}
\dfrac{\partial \vec{u}(\vec{\q})}{\partial \q_i}\bigg\rvert_{\vec{0}} &= \vec{\phi}_i, \\
\dfrac{\partial^2 \vec{u}(\vec{\q})}{\partial \q_i \partial \q_j}\bigg\rvert_{\vec{0}} &= \vec{\Theta}_{ij},
\end{align}

However, this series of operations would lead to an inconsistent formulation in the case of MDs due to their asymmetry, as already outlined in \cite{Jain2017}. In fact, since $\vec{\Theta}_{ij} \neq \vec{\Theta}_{ji}$, it implies that the Schwarz's identity $\partial^2 \vec{u}/ \partial\q_i \partial \q_j \neq \partial^2 \vec{u}/ \partial\q_j \partial \q_i$ is not fulfilled anymore. To overcome this issue, and given the independence of the quadratic mapping on the asymmetric part of each MD shown in~\cite{Jain2017}, the correct strategy proposed in~\cite{Jain2017} is to express both the mapping and its tangent space by means of symmetrized MDs $\bar{\vec{\Theta}}_{ij}=(\vec{\Theta}_{ij}+\vec{\Theta}_{ji})/2$, leading to:
\begin{align}\label{eq:QMNLmap}
&\vec{u}(\vec{\q})
\approx
\sum^n_{i=1}\vec{\phi}_i\q_i
+
\frac{1}{2}\sum^n_{i=1}\sum^n_{j=1}
\bar{\vec{\Theta}}_{ij} \q_i\q_j
=
\vec{\phi}\vec{\q}+\frac{1}{2}\bar{\vec{\Theta}}\vec{\q}\vec{\q},
\\
&\tilde{\vec{\phi}}_i(\vec{\q})
\approx
\vec{\phi}_i
+
\sum^n_{j=1}\bar{\vec{\Theta}}_{ij} \q_j=
\vec{\phi}_i+\bar{\vec{\Theta}}\vec{\q}.
\end{align}
Note that these expressions are used in order to define the reduced-order model, so the dimension $n$ of $\vec{\q}$ is much smaller than the dimension $N$ of $\vec{u}$, $n\ll N$, since only the master coordinates of the ROM are present in $\vec{\q}$. Consequently $\vec{\phi}$ is the matrix of eigenvectors relative to the master coordinates, and should be distinguished from the full matrix of eigenvectors $\vec{\Phi}$ used {\em e.g.} in \eqref{eq:MDmodalchange}.  Finally,  $\bar{\vec{\Theta}}$ is the third-order tensor gathering the MDs $\bar{\vec{\Theta}}_{ij}$.

For future comparison with the normal form method, it is useful to also define the quadratic mapping in modal coordinates:
\begin{equation}\label{eq:QMNLmap_modal}
\vec{X}(\vec{\q})
\approx
\vec{\q}+\frac{1}{2}\bar{\vec{\theta}}\vec{\q}\vec{\q},
\end{equation}
and by components:
\begin{equation}\label{eq:QMNLmap_modal_comp}
X_k
\approx
\q_k+\frac{1}{2}\sum^n_{i=1}\sum^n_{j=1}
\bar{\theta}^k_{ij}\q_i\q_j.
\end{equation}

\subsubsection{Reduced-order model obtained with quadratic manifold}

The nonlinear mapping can then be used in order to derive the reduced-order equations by directly applying Eq.~\eqref{eq:QMNLmap} to the original equations of motion, Eq.~\eqref{eq:first}, and using a standard Galerkin projection. For that purpose, one has to compute the derivatives of Eq.~\eqref{eq:QMNLmap} with respect to time, and finally left-multiply Eq.~\eqref{eq:first} by $\tilde{\vec{\phi}}_i^T$. These derivations have already been proposed in~\cite{Jain2017,Rutzmoser}, we refer the interested reader to these articles for details about the procedure. Here we give the reduced-order dynamics obtained once the projection realised, as a function of the modal coupling coefficients $\vec{g}$ and $\vec{h}$, a derivation that is not given in~\cite{Jain2017,Rutzmoser} and will allow drawing out more direct comparisons with the normal form approach. 

The dynamics for each reduced coordinates $R_p$ finally reads, for $p=1...n$:
\begin{equation}\label{eq:dynaromMD}
\begin{split}
&\ddot{\q}_p
\;+\;
\omega^2_{p} \q_p
\;+\;
\sum_{i,j=1}^n
\left(
(g^p_{ij} +\frac{\omega^2_{p}}{2}\,\bar{\theta}^p_{ij})\;\q_i \q_j
\;+\;
\bar{\theta}^p_{ij}\;(\dot{\q}_i \dot{\q}_j +\q_i\ddot{\q}_j)
\;+\;
\bar{\theta}^j_{pi}\;(\omega^2_{j}\q_i \q_j +\q_i \ddot{\q}_j )
\right)
\;+
\\
&\quad+\;
\sum_{i,j,k=1}^n
\left(
\left(
h^p_{ijk}+
\sum_{s=1}^N 
\left(
\bar{g}^p_{is}\;\bar{\theta}^s_{jk}+
\bar{\theta}^s_{pk}\;
(g^s_{ij}+\frac{\omega^2_{s}}{2}\,\bar{\theta}^s_{ij})
\right)
\right)
\;
\q_i \q_j \q_k
\;+\;
\sum_{s=1}^N 
\left(
\bar{\theta}^s_{pk}\;\bar{\theta}^s_{ij}
\right)
(\dot{\q}_i \dot{\q}_j \q_k+\ddot{\q}_i \q_j\q_k
)
\right) \;=\; 0,
\end{split}
\end{equation}
where the following notations have been introduced for simplifying the expressions~: $\bar{g}^p_{is}=\frac{{g}^p_{is}+{g}^p_{si}}{2}$. Note that this formula simplifies in the case of a symmetric quadratic tensor, which is generally the case in structural mechanics.

One can observe that the linear part is uncoupled, resulting from the fact that the first term of the quadratic manifold in Eq.~\eqref{eq:QMNLmap} is the usual expansion onto the eigenmodes, thus implying, at linear order, uncoupled linear oscillators. The nonlinear terms can be compared to  those obtained when using the normal form approach as nonlinear mapping, Eqs.~\eqref{normadynr}. In particular, one can observe that the normal form approach completely cancels all quadratic terms, provided that no second-order internal resonance are present, a key feature embedded in the derivation which makes the distinction between resonant and non-resonant terms. On the other hand, quadratic terms are always present in~\eqref{eq:dynaromMD}. A second comment is on the presence of terms depending on accelerations in~\eqref{eq:dynaromMD}, not present in the reduced-order dynamics given by the normal form approach.

The restriction to a single master dof is provided, so that one could draw out a term-by-term comparison between the reduced-order dynamics provided by the two methods. Assuming that only mode $p$ is present as reduced coordinates, thus $R_i = 0$, for all $i\neq p$, Eq.~\eqref{eq:dynaromMD} simplifies to:
\begin{equation}\label{eq:dynaromMDSDOF}
\begin{split}
&
\ddot{\q}_p
\;+\;
\omega^2_{p} \q_p
\;+\;
(g^p_{pp} +\frac{\omega^2_{p}}{2}\,\theta^p_{pp})\;\q_p^2
\;+\;
\theta^p_{pp}\;(\dot{\q}_p^2 +\q_p\ddot{\q}_p)
\;+\;
\theta^p_{pp}\;(\omega^2_{p}\q_p^2 +\q_p \ddot{\q}_p )
\;+
\\
&\quad+\;
h^p_{ppp}\q_p^3
\;+\;
\sum_{s=1}^N
\left(
\bar{g}^p_{ps}\;\theta^s_{pp}\;\q_p^3
\;+\;
\theta^s_{pp}\;
(g^s_{pp}+\frac{\omega^2_{s}}{2}\,\theta^s_{pp})\;
\q_p^3
\;+\;
\left(\theta^s_{pp}\right)^2\;
(\dot{\q}_p^2 \q_p+\ddot{\q}_p \q_p^2
)
\right)\;=\;0
\end{split}
\end{equation}
This last equation can then be used either for MD or SMD, so that one can contrast the results obtained by using one of these two strategies (modal derivatives, be they static or dynamic) with the nonlinear change of coordinates provided by normal form theory, which is the aim of the next section.

\subsection{Comparison of the methods and slow/fast approximation}\label{sec:compaSF}

This section aims at comparing the different nonlinear mappings used to derive reduced-order models on the different outcomes they provide: reduced-order dynamics, and prediction of typical nonlinear features such as hardening/softening behaviour, and dependence of mode shapes on amplitude. For that purpose, we restrict ourselves to a single master mode.
Moreover, from now on, we introduce the symmetry property of the quadratic tensor $\vec{g}$ that results from the fact that the internal force derives from a potential, thus leading to $g^i_{jk}=g^i_{kj}$ and $g^i_{jk}=g^j_{ki}=g^k_{ij}$. Note however that, due to our initial choice of fully populated sums and tensors without assuming commutativity of the product, the symmetry property may appear a bit different from {\em e.g.}~\cite{muravyov} when equal indexes are present. Indeed, in~\cite{muravyov} one can read for example $g^p_{ps}=2g^s_{pp}$. This is the only consequence of the initial choice since in~\cite{muravyov} one has  $g^p_{sp}=0$ for $s>p$. In our case, the relationship reads $g^p_{ps}=g^s_{pp}$ and $g^p_{sp}=g^s_{pp}$.

By using such symmetry property, we can also simplify $\bar{g}^p_{ps}=g^p_{ps}=g^s_{pp}$ 
and substituting the value of the modal derivative in modal space $\theta^s_{pp}=-2g^s_{pp}/(\omega^2_{s}-\omega^2_{p})$ when $s\neq p$ and $\theta^p_{pp}=0$ in Eq.~\eqref{eq:dynaromMDSDOF}, one obtains:
\begin{equation}\label{eq:dynaromFMD}
\begin{split}
&
\ddot{\q}_p
\;+\;
\omega^2_{p} \q_p
\;+\;
g^p_{pp}\;\q_p^2
\;+\;
h^p_{ppp}\q_p^3
\;-\;
\underset{s\neq p}{\sum_{s=1}^N }
(g^s_{pp})^2
\frac{2}{\omega^2_{s}-\omega^2_{p}}
\left(
\frac{\omega^2_{s}-2\omega^2_{p}}{\omega^2_{s}-\omega^2_{p}}\;\q_p^3
\;-\;
\frac{2}{\omega^2_{s}-\omega^2_{p}}\;
(\dot{\q}_p^2 \q_p+\ddot{\q}_p \q_p^2
)
\right)\;=\;0.
\end{split}
\end{equation}

If the value of the SMD is used instead of the MD, then the reduced-order dynamics writes:
\begin{equation}\label{eq:dynaromSMD}
\begin{split}
\ddot{\q}_p
\;+\;
\omega^2_{p} \q_p
\;-\;
g^p_{pp}\frac{2}{\omega^2_{p}}\;(\omega^2_{p}\q_p^2+\dot{\q}_p^2 +2\q_p\ddot{\q}_p)
\;+\;
h^p_{ppp}\q_p^3
\;-\;
\sum_{s=1}^N 
(g^s_{pp})^2\frac{2}{\omega^2_{s}}
\left(
\q_p^3
\;-\;
\frac{2}{\omega^2_{s}}(\dot{\q}_p^2 \q_p+\ddot{\q}_p \q_p^2
)
\right)\;=\;0.
\end{split}
\end{equation}
For the explicit comparison, we rewrite the reduced-order dynamics derived with the normal form approach, Eq.~\eqref{dynsingledof}, where the $A^p_{ppp}$ and $B^p_{ppp}$ terms have been expanded:
\begin{equation}\label{eq:singleNF}
\begin{split}
\ddot{\q}_p + \omega_p^2 \q_p +
h_{ppp}^p \q_p^3 -
\sum_{s=1}^N 
({g}^s_{pp})^2
\frac{2}{\omega^2_{s}}
\left(
\frac{\omega^2_s-2 \omega^2_p}{\omega^2_s-4 \omega^2_p}\;\q_p^3
\;-\;
\frac{2}{\omega^2_s-4 \omega^2_p}\;\dot{\q}_p^2 \q_p
\right)\;=\;0.
\end{split}
\end{equation}

\red{Note that the remark on the order of the truncations given at the end of section \ref{subsec:NNMNF} may be better understood from these single-mode reduced dynamics. Eq.~\eqref{eq:singleNF} is the ROM given by normal form, be the calculation of the nonlinear change of coordinate truncated at order two or at order three. Consequently this equation gives the third-order reduced dynamics produced by truncating the normal form at second order. In the same line, Eqs.~\eqref{eq:dynaromFMD} and~\eqref{eq:dynaromSMD} are the third-order reduced-dynamics provided by the quadratic manifold approach. Hence comparing the predictions given by these reduced dynamics is correct since the same order of asymptotic developments is at hand. The only difference one can estimate in the analysis thus relies in the nonlinear mapping, which can be pushed at third-order easily in the normal form approach since the calculation has already been proposed in the past. This means that in the comparisons, the only difference will be on the geometry of the manifold in phase space and the reconstruction formula, but not on the reduced-order dynamics.}

In order to have a better view on the reduced-order dynamics for each of the methods, the general nonlinear oscillator equation describing the dynamics on the reduced subspace can be written under the general form as:
\begin{equation}
\ddot{\q}_p + \omega_p^2 \q_p + C_1 \q_p^2 + C_2 \dfrac{\dot{\q}_p^2}{\omega^2_{p}}+ C_3 \dfrac{\ddot{\q}_p \q_p}{\omega^2_{p}}+ C_4 \q_p^3+ C_5 \dfrac{\dot{\q}_p^2 \q_p}{\omega^2_{p}}+ C_6 \dfrac{\ddot{\q}_p \q_p^2}{\omega^2_{p}}
\;=\;0,
\label{eq:red_sys_general}
\end{equation}
with $C_1$ to $C_6$ different coefficients, which values are summarized in Tables~\ref{tab:coeff_quad},~\ref{tab:coeff_cubic} for the three different methods.
\begin{table}[h!]
\begin{tabular}{llll}\hline
\rule[-13pt]{0pt}{36pt}  & 
$C_1$&$C_2$&$C_3$
\\ \hline
MD\rule[-13pt]{0pt}{36pt} & $g_{pp}^p$ & 0 & 0 
\\ \hline
SMD\rule[-13pt]{0pt}{36pt} & -$2g_{pp}^p$ & 
-$2g^p_{pp} $& 
-$4g^p_{pp}$ 
\\ \hline
NF\rule[-13pt]{0pt}{36pt} & 0 & 
 0 & 0 
\\ \hline
\end{tabular}
\caption{Table of coefficients of the reduced system given by the three methods}
\label{tab:coeff_quad}
\end{table}

\begin{table*}[h!]
\begin{tabular}{llll}\hline
\rule[-13pt]{0pt}{36pt}  & 
$C_4$&$C_5$&$C_6$
\\ \hline
MD\rule[-13pt]{0pt}{36pt}  &
$h_{ppp}^p-\sum_{\underset{s\neq p}{s=1}}^N
(g^s_{pp})^2
\dfrac{2(\omega^2_{s}-2\omega^2_{p})}{(\omega^2_{s}-\omega^2_{p})^2}
$ & $\sum_{\underset{s\neq p}{s=1}}^N
(g^s_{pp})^2
\dfrac{4\;\omega^2_{p}}{(\omega^2_{s}-\omega^2_{p})^2}$ & 
$\sum_{\underset{s\neq p}{s=1}}^N
(g^s_{pp})^2
\dfrac{4\;\omega^2_{p}}{(\omega^2_{s}-\omega^2_{p})^2}$
\\ \hline
SMD\rule[-13pt]{0pt}{36pt} &
$h_{ppp}^p-\sum_{s=1}^N 
(g^s_{pp})^2
\dfrac{2}{\omega^2_{s}}
$ & $\sum_{s=1}^N 
(g^s_{pp})^2
\dfrac{4\;\omega^2_{p}}{\omega^4_{s}}$ & 
$\sum_{s=1}^N 
(g^s_{pp})^2
\dfrac{4\;\omega^2_{p}}{\omega^4_{s}}$
\\ \hline
NF\rule[-13pt]{0pt}{36pt} &
$h_{ppp}^p-\sum_{s=1}^N 
(g^s_{pp})^2
\dfrac{2(\omega^2_{s}-2\omega^2_{p})}{\omega^2_{s}(\omega^2_{s}-4\omega^2_{p})}
$ & $\sum_{s=1}^N
(g^s_{pp})^2
\dfrac{4\;\omega^2_{p}}{\omega^2_{s}(\omega^2_{s}-4\omega^2_{p})}$ & 
0
\\ \hline
\end{tabular}
\caption{Table of coefficients of the reduced system given by the three methods}
\label{tab:coeff_cubic}
\end{table*}
As already remarked, only the normal form approach is able to cancel the quadratic nonlinearity and produce a parsimonious, cubic-order reduced dynamics, depending on two separate coefficients only. Using SMDs creates the larger number of coefficients while only 4 are needed for MDs. Most importantly, the closeness of the results given by the three methods can be underlined in the case where a slow/fast decomposition can be assumed between the master mode $p$ and the slave modes $s$. This case is often encountered in mechanical vibrations since one has often to deal with a large number of modes with very high eigenfrequencies. Let us assume that all the slave modes $s$ are well separated from the master mode, so that for all $s$ one has $\omega_s \gg \omega_p$. It is then very easy to verify on the coefficients given in  Tables~\ref{tab:coeff_quad},~\ref{tab:coeff_cubic} that those provided by the normal form and the MD method tends to the values given by the SMD approach. More specifically, $C_4$ and $C_5$ from normal form exactly match those from the SMD, so that the only difference between the two reduced-order dynamics lies in the additional terms $C_1$, $C_2$, $C_3$ and $C_6$ for the SMD method. On the other hand, using the slow/fast approximation for the coefficients provided by the MD shows that $C_4$, $C_5$ and $C_6$ tends exactly to the values obtained with SMDs, the only difference being in the summation, where the $p$ term is excluded in the MD approach whereas it is not in the SMD, as a direct consequence from Eq.~\eqref{eq:mdmodalexpr}. Indeed, Eq.~\eqref{eq:mdmodalexpr-b} shows that for MD, the $g^p_{pp}$ term is taken into account in the amplitude-frequency relationship, and not in the reconstruction of the vector as given by Eq.~\eqref{eq:MD}. On the other hand for SMD, the $g^p_{pp}$ term is taken into account in the vector defining the SMD, Eq.~\eqref{eq:SMD}, but not in the frequency dependence on amplitude. This important difference between the two methods will have consequences that are commented further in the next sections, and the $g^p_{pp}$ will be denoted further as the self-quadratic term.

In order to better understand the observed differences on the reduced-order dynamics, a fair comparison has to be given not onto a term-by-term comparison, since the meaning of the reduced variables is not the same, but on the general predictions given by each reduction method on the most important nonlinear features. The next sections are thus devoted to comparing the prediction of the type of nonlinearity provided by each method ({\em i.e.} the first term in the amplitude-frequency relationship that dictates the hardening or softening behaviour), as well as the mode shape dependence on amplitude.

\subsubsection{Hardening/softening behaviour}\label{sec:HS}

The generic reduced-order dynamics, Eq.~\eqref{eq:red_sys_general}, can be solved with a perturbation method in order to derive the type of nonlinearity predicted by each method. Keeping the general notation with the $C_i$ coefficients for the ease of reading, the general solution up to second order in amplitude reads:
\begin{equation}
\q_p=a_0 \cos[\,\omega_p\,t\,( 1 +  \Gamma a_0^2)] +a_0^2 \left(  \frac{C_1 - C_2 - C_3}{6\,\omega_p^2}\cos[\,2\,\omega_p\,t\,( 1 + a_0^2 \Gamma)]
-\frac{C_1 + C_2 - C_3}{2\,\omega_p^2}\right)
+
\mathcal{O}(a_0^3),
\label{eq:multiple_scales_sol}
\end{equation}
with $a_0$ the amplitude, and $\Gamma$ the general coefficient that dictates the hardening/softening behaviour. Indeed, one can introduce the nonlinear frequency $\omega_{\mathrm{NL}} = \omega_p( 1 +  \Gamma a_0^2)$. If $\Gamma > 0$ then the system is hardening. The general expression for $\Gamma$ with all the $C_i$ coefficients writes:
\begin{equation}
\Gamma=-\frac{1}{24\,\omega_p^4}\left(
  10\,C_1^2
+ 10\,C_1 C_2
+  4\,C_2^2
-  7\,C_2 C_3
+   \,C_3^2
- 11\,C_1 C_3
\right)
+\frac{1}{8\,\omega_p^2}\left(
3\,C_4
+\,C_5
-3\,C_6\right).
\end{equation}

One can note in particular that with the normal form approach, one has $C_1=C_2=C_3=0$ since the method fully cancels the quadratic terms, so that there is no second harmonic term in the reduced-order dynamics and Eq.~\eqref{eq:multiple_scales_sol} reduces to its first term at order two. However, since quadratic terms are present in the nonlinear change of coordinates, this simplification does not imply that the second harmonic is not present in the reconstructed displacements, as it will be shown in the next section. Once again, these two last equations show that normal form approach produces a parsimonious representation of the reduced dynamics which is generally easier to read and interpret.

Replacing the values of the $C_i$ coefficients obtained for each method (MD, SMD or NF for normal form), one arrives at the prediction of the type of nonlinearity provided by each reduced-order model as:
\begin{subequations}\label{eq:zifullgamma}
\begin{align}
\Gamma_\text{MD}=&
-
\dfrac{5}{12\;\omega^2_{p}}
\left(\dfrac{g^p_{pp}}{\omega_{p}}\right)^2
+
\dfrac{3}{8\;\omega^2_{p}}\left(
h_{ppp}^p
-
\sum_{\underset{s\neq p}{s=1}}^N
2\left(\dfrac{g^s_{pp}}{\omega_{s}}\right)^2
\left(
1+
\dfrac{\omega^2_p(4\,\omega^2_s-3\,\omega^2_p)}{3(\omega^2_s-\omega^2_p)^2}
\right)\right), \label{eq:zifullgammaA}
\\
\Gamma_\text{SMD}=&
-
\dfrac{5}{12\;\omega^2_{p}}
\left(\dfrac{g^p_{pp}}{\omega_{p}}\right)^2
+
\dfrac{3}{8\;\omega^2_{p}}\left(
h_{ppp}^p
-
\sum_{\underset{s\neq p}{s=1}}^N
2\left(\dfrac{g^s_{pp}}{\omega_{s}}\right)^2
\left(
1+\dfrac{4\,\omega^2_p}{3\,\omega^2_s}
\right)\right), \label{eq:zifullgammaB}
\\
\Gamma_\text{NF}=&
-
\dfrac{5}{12\;\omega^2_{p}}
\left(\dfrac{g^p_{pp}}{\omega_{p}}\right)^2
+
\dfrac{3}{8\;\omega^2_{p}}\left(
h_{ppp}^p
-
\sum_{\underset{s\neq p}{s=1}}^N
2\left(\dfrac{g^s_{pp}}{\omega_{s}}\right)^2
\left(
1+
\dfrac{4\,\omega^2_p}{3(\omega^2_s-4\,\omega^2_p)}
\right)\right). \label{eq:zifullgammaC}
\end{align}
\end{subequations}

One can note that the first terms of the prediction are the same, while the difference arise from the way the slave (or neglected) coordinates are taken into account in order to predict the type of nonlinearity. This feature is however key in order to give a correct prediction since there is a strong need to take properly into account the curvature of the manifolds in phase space, otherwise incorrect predictions are given~\cite{touze03-NNM}.

In order to give more insights into Eqs.~\eqref{eq:zifullgamma}, let us first notice that in the summed terms, the first one is always the same since the different expressions all start with $1+...$. Let us isolate this term and introduce the following notation~:
\begin{equation}
\mathcal{C}^s_\text{SC} = 2\left(\dfrac{g^s_{pp}}{\omega_{s}}\right)^2.
\end{equation}
One can notice that this correction term is the one obtained by using static condensation, as already shown for example in~\cite{Vizza3d,YichangICE}, thus the subscript $\text{SC}$. Denoting as $\mathcal{C}_\text{MD}$, $\mathcal{C}_\text{SMD}$ and $\mathcal{C}_\text{NF}$ the correction factors given by each method ({\em i.e.} the term in the summation), one can then simply compares all these terms to $\mathcal{C}^s_\text{SC}$ in order to have an expression depending only on the eigenfrequencies. Assuming that there is only one slave mode $s$ in the summation in order to highlight the contribution brought by each term,  the following ratios can be written:
\begin{subequations}\label{eq:Cratio}
\begin{align}
\dfrac{\mathcal{C}_\text{MD}}{\mathcal{C}_\text{SC}}\;\,=\;&
1+\dfrac{4}{3}
\dfrac{\rho^2-3/4}{(\rho^2-1)^2},
\\
\dfrac{\mathcal{C}_\text{SMD}}{\mathcal{C}_\text{SC}}=\;&
1+\dfrac{4}{3}
\dfrac{1}{\rho^2},
\\
\dfrac{\mathcal{C}_\text{NF}}{\mathcal{C}_\text{SC}}\;\,=\;&
1+\dfrac{4}{3}
\dfrac{1}{\rho^2-4},
\end{align}
\end{subequations}
where $\rho=\omega_s/\omega_p$ has been introduced in order to highlight their behaviour with respect to the fulfilment of the slow/fast partition. These expressions clearly underline the fact that each method refine the correction factor of static condensation by an additional term. One can also observe that the refinement of the $\mathcal{C}^s_\text{SMD}$ comes from the inertia and velocity terms $C_5$ and $C_6$, whereas the term $C_4$ is exactly the one from static condensation. Consequently, using SMD without quadratic manifold could lead to erroneous predictions since inertial and velocity corrections could be missed. This remark should be particularly relevant in a case of geometric nonlinearity involving inertia, as {\em e.g.} in the case of a cantilever beam.

To better assess the quality of the predictions given by the three methods, Eqs.~\eqref{eq:Cratio} can be Taylor-expanded by using the slow/fast assumption $\omega_s \gg \omega_p$ for the slave modes $s$. This assumption allows introducing a small parameter $\omega_p/\omega_s$, or, equivalently, considering the expansion under the assumption $\rho\rightarrow\infty$. One then obtains:
\begin{subequations}\label{eq:CratioTaylor}
\begin{align}
\dfrac{\mathcal{C}_\text{MD}}{\mathcal{C}_\text{SC}}\;\,=\;&
1+\dfrac{4}{3}\dfrac{1}{\rho^2}+
\sum_{i=2}^\infty \dfrac{3+i}{3\;\rho^{2i}},
\\
\dfrac{\mathcal{C}_\text{SMD}}{\mathcal{C}_\text{SC}}=\;&
1+\dfrac{4}{3}
\dfrac{1}{\rho^2},
\\
\dfrac{\mathcal{C}_\text{NF}}{\mathcal{C}_\text{SC}}\;\,=\;&
1+\dfrac{4}{3}\dfrac{1}{\rho^2}+
\sum_{i=2}^\infty \dfrac{4^i}{3\;\rho^{2i}}.
\end{align}
\end{subequations}
These formulas show in particular that all the methods predict the same first two terms in the expansion that assumes slow/fast partition, and the limit for  $\rho\rightarrow\infty$ is the same for all methods, including static condensation, since the ratios tends to 1 in this case. This means that a formal equivalence in the prediction of the type of nonlinearity is obtained only in the limit case of $\omega_s \gg \omega_p$ for all the studied methods. Fig.~\ref{fig:gamma_flat} illustrates this convergence and shows that it is obtained rapidly, indicating in particular that from the value $\omega_s / \omega_p \simeq 4$, all methods are almost converged in terms of type of nonlinearity, thus quantifying more properly the value from which the slow/fast partition is effective so that one can use the methods based on modal derivatives safely. \red{In order to be a bit more quantitative, one can remark that the relative difference between $\mathcal{C}_\text{MD}$ and $\mathcal{C}_\text{NF}$ is equal to $5\%$ for $\rho=3.25$ and $1\%$ for $\rho=4.6$, so that the proposed bound $\omega_s / \omega_p \simeq 4$ has not to be understood as a strict one. Moreover, the error on $\Gamma$ will be smaller than the error on the correction factor $\mathcal{C}$, being $\Gamma$ composed of other terms that are not affected by the reduction method. The conclusion is that  $\rho \in[3,4]$ can be understood as a transition region, and converged results thanks to slow/fast assumption can be faithfully obtained over 4, but below 3 caution has to be exercised. }


\begin{figure}[h!]
\centering
\begin{subfigure}{.44\textwidth}
\includegraphics{Gamma_taylor.tikz}
\end{subfigure}
\caption{Evolution of the ratios $\dfrac{\mathcal{C}_\text{MD}}{\mathcal{C}_\text{SC}}$, $\dfrac{\mathcal{C}_\text{SMD}}{\mathcal{C}_\text{SC}}$ and $\dfrac{\mathcal{C}_\text{NF}}{\mathcal{C}_\text{SC}}$, defined in Eqs.~\eqref{eq:Cratio}, as a function of the parameter $\rho=\omega_s/\omega_p$, from which the behaviour of the type of nonlinearity defined by the $\Gamma$ coefficients in Eqs.~\eqref{eq:zifullgamma}, can be directly deduced. Dashed grey line is the (constant) value predicted by static condensation (all curves are normalized with respect to this value). Yellow curve: $\dfrac{\mathcal{C}_\text{SMD}}{\mathcal{C}_\text{SC}}$ predicted by static modal derivatives; orange curve:  $\dfrac{\mathcal{C}_\text{MD}}{\mathcal{C}_\text{SC}}$ computed from modal derivatives, blue curve: $\dfrac{\mathcal{C}_\text{NF}}{\mathcal{C}_\text{SC}}$ given by normal form theory.}
\label{fig:gamma_flat}
\end{figure}

Fig.~\ref{fig:gamma_flat} shows also other interesting features on the behaviour of the type of nonlinearity. Besides the convergence of all curves in the limit $\rho\rightarrow\infty$, important differences occur in the regions where the methods have a singularity. The normal form approach displays a singular behaviour in the vicinity of the 1:2 internal resonance when $\omega_s \simeq 2\omega_p$. This fact is logical and has already been commented in numerous prior publications. Indeed, when such a resonance exists, then a strong coupling arises between the two modes, so that reducing the dynamics to a single master mode has no meaning anymore, and the minimal model should be composed at least by these two internally resonant modes. The divergence in the behaviour of $\mathcal{C}_\text{NF} / \mathcal{C}_\text{SC}$ reflects this fact, meaning that in this zone the definition of the type of nonlinearity is of no more use since another dynamical regime takes place. Previous publications also clearly underlines that the prediction given by $\Gamma_\text{NF}$ in Eq.~\eqref{eq:zifullgammaC} is correct~\cite{touze03-NNM}, which has been confirmed with comparisons to direct simulations of the full-order model, and this prediction of the type of nonlinearity has then been used for continuous structures such as cables and shells~\cite{touze-shelltypeNL,regacarbo00,Arafat03,Pellicano02}.

On the other hand, the prediction given by MD displays a divergence at the 1:1 resonance, when the slave and master modes have close eigenfrequencies, $\omega_s \simeq \omega_p$. This divergence does not rely on a firm theoretical result from dynamical systems. Indeed, even though in the case a 1:1 internal resonance exists so that the two modes need to be taken into account to study the coupled dynamics, uncoupled solutions still exist and the backbone curves of these uncoupled solutions can be computed, thus preserving the meaning of the $\Gamma$ coefficients defined in Eqs.~\eqref{eq:zifullgamma}, see {\em e.g.} \cite{touze01-JSV,manevitch2003free,Givois11}. Thus the divergence of $\mathcal{C}_\text{MD} / \mathcal{C}_\text{SC}$ is interpreted as a failure of the method.  Finally, for small values of $\rho$, one can observe that the SMD method shows a singular behaviour, and will predict unreasonably stiff behaviour. On the other hand, MD method gives a finite value, which is a bit different from the correct one given by normal form approach. All these results underline that MD and SMD can be used safely only when the assumption $\omega_s > 4\omega_p$ is fulfilled, otherwise unreliable predictions may be given by these two methods.

\subsubsection{Drift and mode shapes}\label{sec:drift}

A second comparison on the global outcomes of the three method can be provided by contrasting the mode shape dependence on amplitude. Indeed, assuming a single mode motion with $R_s = 0$ for all $s\neq p$ (only master mode $p$ participates to the vibration), allows recovering the amplitude dependence of the $p$-th mode shape. At small amplitude, the three methods recover the usual eigenmode, but they then differ in the way they are taking into account the cross-couplings with slave modes. Let us denote as $\vec{u}_\text{MD}$, $\vec{u}_\text{SMD}$ and $\vec{u}_\text{NF}$ the physical displacement following single-mode motion for each of the three methods. Using the previous formula allows one to reconstruct
\begin{align}
\vec{u}_\text{MD}(t)=&\vec{\phi}_p \q_p(t)
-
\underset{s\neq p}{\sum_{s=1}^N }
\frac{g^s_{pp}}{\omega_s^2-\omega_p^2}
\q_p^2(t) \vec{\phi}_s,\label{eq:reconsuMD}
\\
\vec{u}_\text{SMD}(t)=&\vec{\phi}_p \q_p(t)
-
\vec{\phi}_p
\frac{g^p_{pp}}{\omega_p^2}
\q_p^2(t)
-
\underset{s\neq p}{\sum_{s=1}^N}
\frac{g^s_{pp}}{\omega_s^2}
\q_p^2(t) \vec{\phi}_s, \label{eq:reconsuSMD}
\\
\vec{u}_\text{NF}(t)=&\vec{\phi}_p \q_p(t)
-
\vec{\phi}_p
\frac{g^p_{pp}}{\omega_p^2}
\frac{1}{3}
\left(
\q_p^2(t)
+
\frac{2}{\omega_p^2}\dot{\q}_p^2(t)
\right)
-
\underset{s\neq p}{\sum_{s=1}^N}
\frac{g^s_{pp}}{\omega_s^2}
\left(
\frac{\omega_s^2-2\omega_p^2}{\omega_s^2-4\omega_p^2}
\q_p^2(t)
-
\frac{2}{\omega_s^2-4\omega_p^2}
\dot{\q}_p^2(t)
\right) \vec{\phi}_s.\label{eq:reconsuNF}
\end{align}

Comparing the mode shapes given by MD and SMD, one can already underline that the summed term given by MD reduces to that given by SMD if one considers the slow/fast assumption with $\omega_s \gg \omega_p$. However a difference persists in the two methods since with SMD an added quadratic term, depending on mode $p$ only, is present (second term in \eqref{eq:reconsuSMD}). This comes again from the treatment of the self-quadratic $g^p_{pp}$ term in Eqs.~\eqref{eq:mdmodalexpr}, already underlined in Sect.~\ref{sec:MDmodal}. Indeed, the $g^p_{pp}$ term  for the MD method is not present in the reconstruction, but in the dependence of the nonlinear frequency with amplitude instead, while the SMD method distributes the influence of this $g^p_{pp}$ term on the spatial reconstruction, but not on the amplitude-frequency relationship. This explains why the prediction of the hardening/softening behaviour appears to be more general for the MD method than for the SMD. Comparing now with the normal form approach, one can see that NF reduction gives rise to velocity-dependent terms in these formula, a feature that is not present in the other method, which is a direct consequence of the fact that NF method takes into account both independent displacement and velocity variables as it should be from a dynamical system perspective.

Again, one can also observe that the summed term in \eqref{eq:reconsuNF} reduces (at first significant order) to that provided by SMD when the slow/fast assumption is at hand, showing that the SMD method provides the most simplified expressions.

From the general expressions given in \eqref{eq:reconsuMD}-\eqref{eq:reconsuNF}, one can isolate the constant term (zero-th harmonic) which is produced by the quadratic nonlinearity, in order to compare more closely one term of this general expansion. This constant term is known as a drift since it corresponds to the fact that due to quadratic nonlinearity, the oscillations are no more centred around zero, and it has already been compared for different reduction methods, see {\em e.g.}~\cite{nayfehcarbo97,touze03-NNM}. One can then simply replace  $\q_p(t)$ by the expression given by Eq.~\eqref{eq:multiple_scales_sol}; while the values of $\q_p^2(t)$ and $\dot{\q}_p^2(t)$ up to second order write:
\begin{align}
\q_p^2(t)=&
\frac{a_0^2}{2}\left(1+\cos[2\omega_\text{NL}t]\right)+\mathcal{O}(a_0^3)
\\
\dot{\q}_p^2(t)=&
\frac{a_0^2}{2}\omega^2_\text{NL}\left(1-\cos[2\omega_\text{NL}t]\right)+\mathcal{O}(a_0^3)
\end{align}
where the nonlinear frequency $\omega_\text{NL}=\omega_p( 1 + a_0^2 \Gamma)$ has been introduced. Isolating the constant term leads to the following expressions for the drift $d$ predicted by each reduction method:
\begin{align}
\vec{d}_\text{MD}=&
\frac{a_0^2}{2}
\left(
-
\frac{g^p_{pp}}{\omega_p^2}\vec{\phi}_p
-
\underset{s\neq p}{\sum_{s=1}^N }
\frac{g^s_{pp}}{\omega_s^2-\omega_p^2}
\vec{\phi}_s
\right),
\\
\vec{d}_\text{SMD}=&
\frac{a_0^2}{2}
\left(
-
\frac{g^p_{pp}}{\omega_p^2}\vec{\phi}_p
-
\underset{s\neq p}{\sum_{s=1}^N }
\frac{g^s_{pp}}{\omega_s^2}
\vec{\phi}_s
\right),
\\
\vec{d}_\text{NF}=&
\frac{a_0^2}{2}
\left(
-
\frac{g^p_{pp}}{\omega_p^2}
\left(
\frac{1}{3}+\frac{2}{3}\frac{\omega^2_\text{NL}}{\omega_p^2}
\right)
\vec{\phi}_p
-
\underset{s\neq p}{\sum_{s=1}^N }
\frac{g^s_{pp}}{\omega_s^2}
\left(
1-\frac{2(\omega^2_\text{NL}-\omega_p^2)}{\omega_s^2-4\omega_p^2}
\right)
\vec{\phi}_s
\right).
\end{align}
One can observe that assuming slow/fast dynamics, the drift predicted by MD reduces to that given by SMD. On the other hand, one can also see that in order to retrieve the drift predicted by SMD from $\vec{d}_\text{NF}$, one has to assume that the deviation of the nonlinear frequency $\omega_\text{NL}$ is small as compared to the linear frequency so that $\omega_\text{NL} \simeq \omega_p $. Hence the prediction of the mode shape dependence on amplitude given by SMD is reliable only in the case where the backbone curve does not depart severely from the linear resonance, which is a strong assumption. 

In order to point out a last difference on the theoretical expressions which will have important consequences in the next sections, let us also follow the first harmonic of the solution in the reconstruction procedure. Using Eq.~\eqref{eq:multiple_scales_sol} to define the harmonic content of the master variable, and going back to the harmonic content of the modal coordinates $X_i$ defined using either the QM method, Eq.~\eqref{eq:QMNLmap_modal_comp}, or the normal form approach, Eqs.~\eqref{eq:transformNF}-\eqref{eq:invarman}, one can easily follow the first harmonic and retrieve its expression in the modal coordinates. Since $p$ is the master mode and at lowest order $X_p=R_p$, then the most important contribution is present in $X_p$ as compared to other $X_k$'s. Let us denote as  $[X_p^{(H1)}]_\text{MD}$ the first harmonic for the MD approach (and SMD and NF for the other two methods), these expressions write:
\begin{subequations}\label{eq:X_1H1expr}
\begin{align}
[X_p^{(H1)}]_\text{MD}=&
a_0
\cos(\omega_\text{NL} t)
\left(1+\mathcal{O}(a_0^4)\right)
,\label{eq:X_1H1exprMD}
\\
[X^{(H1)}_p]_\text{SMD}=&
a_0
\cos(\omega_\text{NL} t)
\left(1-a_0^2\;\dfrac{2}{3}\left(\dfrac{g^p_{pp}}{\omega_p^2}\right)^2
+\mathcal{O}(a_0^4)
\right)
,\label{eq:X_1H1exprSMD}
\\
[X^{(H1)}_p]_\text{NF}=&
a_0
\cos(\omega_\text{NL} t)
\left(1+\mathcal{O}(a_0^4)\right).\label{eq:X_1H1exprNF}
\end{align}
\end{subequations}
They underline the importance of the treatment of the self-quadratic $g^p_{pp}$ term between MD and SMD method. Indeed, whereas the amplitude $a_0$ defined from \eqref{eq:multiple_scales_sol} corresponds, for the MD and NF cases, to the amplitude of the first harmonic in $X_p$, this is not the case for the QM derived from SMD. In that case, the amplitude has an extra term implying the self-quadratic coupling term. Importantly, this term appears as a difference so that the amplitude of the first harmonic can tend to small values with increasing $a_0$. Whereas all the comparisons led in this section  shows that the methods tend to be equivalent under a slow/fast assumption, this last expression highlights the fact that, for the SMD method, the amplitude of the master mode can be very different from the amplitude of the initial coordinate. The consequence of this finding will be more clearly illustrated  in the next sections on examples, and will be key to understand why the SMD method can fail even under the slow/fast assumption. 
\section{Comparison on two degrees of freedom systems}

In this section, the comparisons drawn out on the theoretical expressions are illustrated on two dofs systems, in order to highlight the main differences on simple cases. Two different models are selected. The first one is derived from the equations of motion of a beam, and is selected in order to mimic the nonlinearities present in a flat symmetric system, where these simplifying assumptions help in letting the methods based on SMD work properly. The second example has important quadratic couplings and better accounts from the problems arising with curved structures such as arches and shells.

\subsection{A two-dof model representing a flat symmetric structure}\label{sec:2dofsFLAT}

\subsubsection{Presentation of the model}
\label{sec:2dofsFLATmodel}

The particular nature of the nonlinear couplings in the case of flat symmetric structures such as beams and plates, relies on the simplifying facts that flexural and in-plane modes are linearly uncoupled, and their nonlinear couplings involve simple terms that can be easily traced from the \vonkar models. These simplifications have been used in numerous recent papers in order to explain why a number of methods for producing ROMs are able to predict very good results in this case, see {\em e.g.}~\cite{Jain2017,givois2019,Vizza3d,VERASZTO}. In order to propose a simple two-dofs system mimicking these particular relationships, the \vonkar model for slender beams is used and simplified to two vibration modes, one flexural and one longitudinal, in order to produce the simplified system, from which the coefficients can be related to meaningful quantities of the beam and in particular to its slenderness.

A non-prestressed beam of length $L$ is thus considered, with a uniform rectangular cross section of area $S=bh$ ($h$ being the thickness and $b$ the width) and second moment of area $I=bh^3/12$ , made in an homogeneous and isotropic material of Young's modulus $E$ and density $\delta$. Boundary conditions are clamped at $X=0$ and $X=L$.

The equations of motion for the transverse displacement $W(X,T)$, and the longitudinal displacement $U(X,T)$ ($X$ and $T$ being the dimensional space and time variables), assuming \vonkar theory, reads~\cite{Nayfeh79,givois2019}:
\begin{subequations}\label{eq:VKbeamdim}
\begin{align}
&  \ddot{W}+\dfrac{E I}{\delta S} W^{''''}-\dfrac{E}{\delta} \left(U^{'}W^{'}+\dfrac{1}{2}W^{'\,3}\right)^{'}=0,\label{eq:VKbeamdim-a}\\
&  \ddot{U} - \dfrac{E}{\delta} (U{''}+W^{'}W^{''})=0.\label{eq:VKbeamdim-b}
\end{align}
\end{subequations}
A particular feature of Eqs.~\eqref{eq:VKbeamdim} is that the longitudinal displacements are only quadratically coupled with the transverse, as shown in~\eqref{eq:VKbeamdim-b}. On the other hand, the only nonlinear terms appearing on the equations of motion for the flexural term $W$ are: (i) a quadratic coupling involving a product between one in-plane and one transverse component, and a cubic term with only transverse components, see Eq.~\eqref{eq:VKbeamdim-a}.

Following \cite{givois2019}, the equations of motion can be made nondimensional so that the resulting system depends only on two physically meaningful parameters:  the slenderness ratio $\sigma=h/L$, and the wavelength $\beta$ appearing naturally when solving the eigenvalue problem. Indeed, focusing on the linear problem for the transverse motion, the eigenvalue problem $ \phi^{''''} = \omega^2 \frac{\delta S}{EI} \phi$ is solved by using a combination of sine, cosine, hyperbolic sine and hyperbolic cosine functions of $kx$, with $k$ dimensional wavelength such that $k^4 =  \frac{\delta S}{EI} \omega^2$ and $\beta = kL$. Assuming clamped boundary conditions the characteristic equation for $\beta$, from which the eigenfrequencies are deduced, reads : $\cos(\beta)\cosh(\beta)=1$. The reader is referred to Appendix~\ref{app:VK2dofs} for the details of this classical derivation.

Introducing the thickness $h$ as characteristic length, so that the nondimensional displacements are as $w=W/h$ and $u=U/h$, normalizing time using  $t=T (\beta^2/L^2\sqrt{EI/\delta S})$ and the space variable with the beam length, $x=X/L$; Eqs.~\eqref{eq:VKbeamdim} are rewritten as follows:
\begin{subequations}\label{eq:VKbeamNONdim}
\begin{align}
&  w_{,tt}+\dfrac{1}{\beta^4} w_{,xxxx}
-\dfrac{12}{\beta^4 \sigma }
\left(u_{,x}w_{,x}\right)_{,x}
-\dfrac{6}{\beta^4}
\left({w_{,x}}^3\right)_{,x}=0,\\
&  u_{,tt} - \frac{12}{\beta ^4 \sigma ^2}u_{,xx}
-\frac{12}{\beta^4 \sigma }w_{,x}w_{,xx}=0.
\end{align}
\end{subequations}
In order to derive a minimal two dofs system from these equations, we select the first flexural eigenmode and the most important longitudinal mode coupled to the first flexural. From previous studies, see {\em{e.g.}}~\cite{SOMBROEK2018,givois2019,Vizza3d}, it is known that the fourth in-plane mode is strongly coupled to the first flexural. Let us denote as $q_1$ the modal amplitude of the first transverse mode and $p_4$ the modal amplitude of the fourth in-plane mode (see Appendix~\ref{app:VK2dofs} for the details). Using a standard Galerkin projection (see {\em e.g.}~\cite{givois2019}), Eqs.~\eqref{eq:VKbeamNONdim} can be rewritten as
\begin{subequations}\label{eq:vk2dofsV1}
\begin{align}
&\ddot{q}_1+ q_1+ \dfrac{2G}{\sigma}p_4 q_1 + D q_1^3=0,\label{eq:vk2dofsV1-a}\\
&\ddot{p}_4+\dfrac{(4\pi)^212}{\beta^4 \sigma^2} p_4 + \dfrac{G}{\sigma}q_1^2=0,\label{eq:vk2dofsV1-b}
\end{align}
\end{subequations}
where $D$ and $G$ are the nonlinear coupling coefficients arising from the Galerkin projection, and involves integral on the length of products of derivatives of the mode shape functions, see~\cite{givois2019} for the general calculation and  Appendix~\ref{app:VK2dofs} for the detailed expression of these two coefficients. One  can note in particular that, due to the choice of the nondimensional time to arrive at Eqs.~\eqref{eq:VKbeamNONdim}, the eigenfrequency of the first flexural mode is 1, while the natural frequency of the fourth in-plane mode reads $\omega_2^2 = \dfrac{(4\pi)^212}{\beta^4 \sigma^2}$. Due to the normalisation selected (involving $\omega_1 = 1$ for the fundamental mode), the term in factor of $p_4$ in Eq.~\eqref{eq:vk2dofsV1-a} can be easily interpreted as the square of the ratio $\rho=\omega_2/\omega_1$, recovering the term introduced in Sect.~\ref{sec:HS}. Thanks to its explicit expression, $\rho$ can now be directly related to the slenderness ratio:
\begin{equation}\label{eq:slenderness}
\rho=\dfrac{4\pi\;\sqrt{12}}{\beta^2}\;\dfrac{1}{\sigma}\approx 1.95 \dfrac{1}{\sigma}.
\end{equation}

So that the final two-dofs system that will be used for the investigations reads:
\begin{align}
&\ddot{X}_1+ X_1+ 2\,\bar{G}\,\rho \,X_1 X_2 + D X_1^3=0,\\
&\ddot{X}_2+\rho^2 X_2 + \bar{G}\,\rho X_1^2=0,
\end{align}
where $\bar{G}=G\,\beta^2/(4\pi\sqrt{12})$ has been introduced for the ease of reading. Also the notation for the variables has been changed with $X_1=q_1$ and $X_2=p_4$ for the sake of simplicity. \red{A particular feature of this system is that the coupling between master and slave mode is purely quadratic. Consequently the potential third-order tensors from the normal form approach are all vanishing. In this case, the two nonlinear mappings are thus exactly at the same order due to the very simplified shape of the starting equations.}


This system is now investigated in order to see how the methods under study behaves when reducing the system to its first (flexural) mode  using different nonlinear mappings. The advantage of this formulation is that all coefficients are related to a physical problem so that some insights can be given to the results obtained with this simplistic model with regard to continuous problems. In particular, Sect.~\ref{sec:compaSF} underlined that all methods show a convergence on some properties when a slow/fast assumption is assumed, 
which has been quantified precisely on the type of nonlinearity as occurring for $\rho > 4$. Also, divergent behaviours has been underlined and explained for $\rho \simeq 1$ (case of MD) and $\rho \simeq 2$ (case of normal form). Consequently the system will be studied for three different values close to these points, namely $\rho=1.25$,  $\rho=2.5$ and $\rho=10$. Note that a beam is generally considered as slender if $\sigma \leq 1/20$. Thanks to Eq.~\eqref{eq:slenderness}, this means that $\rho \geq 40$. The consequence of this remark is that in all slender beams the slow/fast assumption is very well fulfilled, and our study concerns specific cases occurring for very thick beams. Regarding the nonlinear coefficients $G$ and $D$, they only depend on the modes selected in the expansion. In our study, we will always consider the first flexural and fourth axial, so that $G$ and $D$ are constants since they only depend on the nondimensional shape functions of the selected modes. In the remainder of the study, we have selected  $D=2.67$, $\bar{G}=0.63$.

\subsubsection{Results}

The comparisons between the different methods are drawn out on the geometry of the manifolds, as well as on the dynamics onto these manifolds, described by the frequency-amplitude relationship (backbone curve). All the solutions are computed thanks to a numerical continuation method using the asymptotic-numerical method, implemented in the software Manlab, where the unknowns are represented thanks to the harmonic balance method~\cite{COCHELIN2009,LAZARUS20105,Guillot2019}.  After a convergence study, the number of harmonics retained in the computations is 7. In each case, the master mode is the fundamental one, $X_1$, and the slave mode $X_2$. The dynamics onto the reduced subspaces is given by Eq.~\eqref{eq:dynaromFMD} when using the MD approach, Eq.~\eqref{eq:dynaromSMD} if one considers SMD instead, and Eq.~\eqref{eq:singleNF} with the normal form method, with $R_1$ the master coordinates. For the reduced models, continuation is performed on the master coordinate in order to compute the frequency-amplitude relationships. From these values, the nonlinear mappings, given either by Eqs.~\eqref{eq:zecvNL} for the normal form approach, or by Eqs.~\eqref{eq:QMNLmap_modal_comp} for the QM method, allows to retrieve the initial modal amplitude $X_1$ and $X_2$. From all these data, one can plot either the geometry of the manifolds in phase space $(X_1,Y_1,X_2,Y_2)$, or the backbone curves.

\begin{figure*}[h!]\centering
\begin{subfigure}[b]{.32\textwidth}
\includegraphics[width=\textwidth]{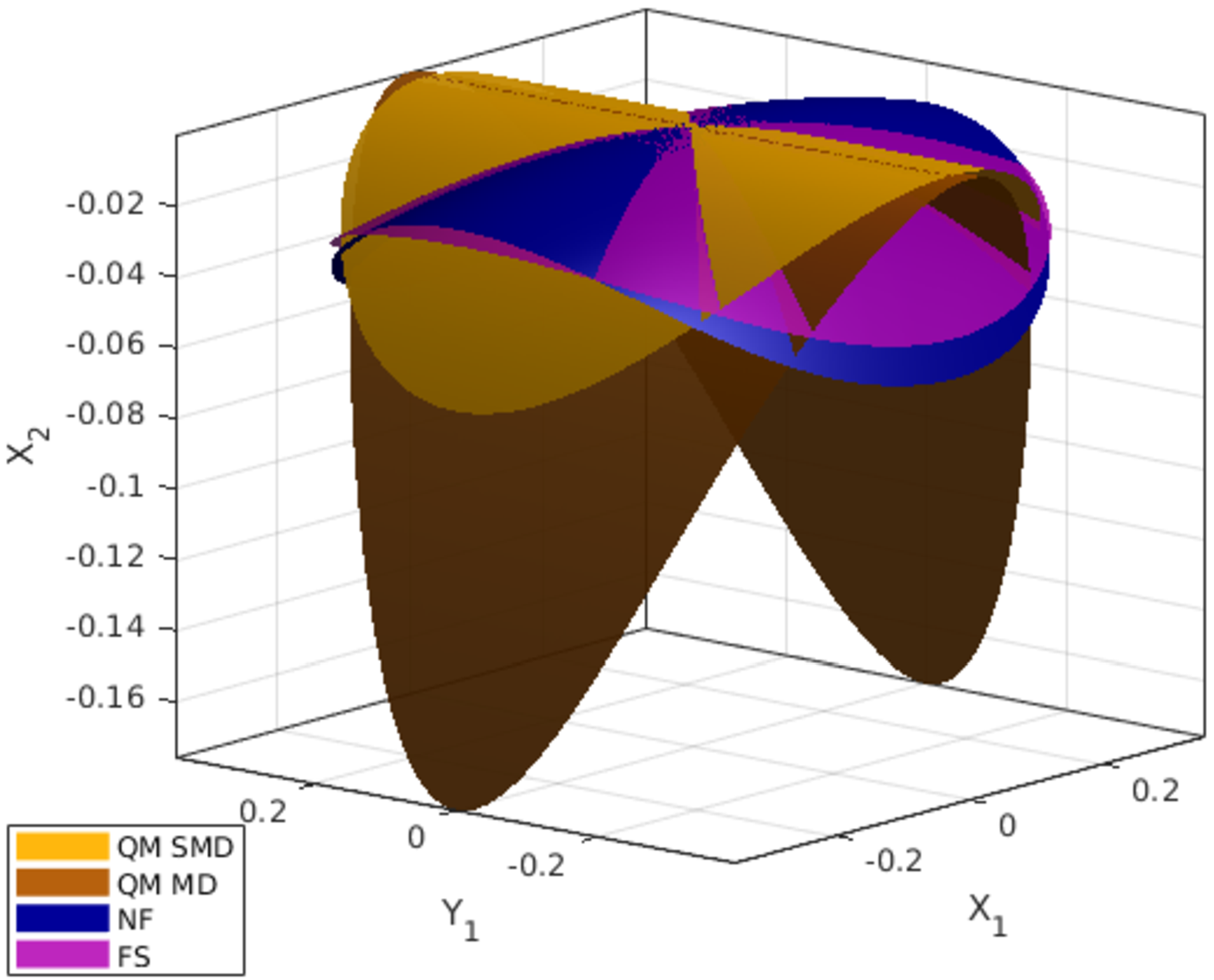}
\caption{$\rho=1.25$.}
\label{fig:ex1manifold-a}
\end{subfigure}
\begin{subfigure}[b]{.32\textwidth}
\includegraphics[width=\textwidth]{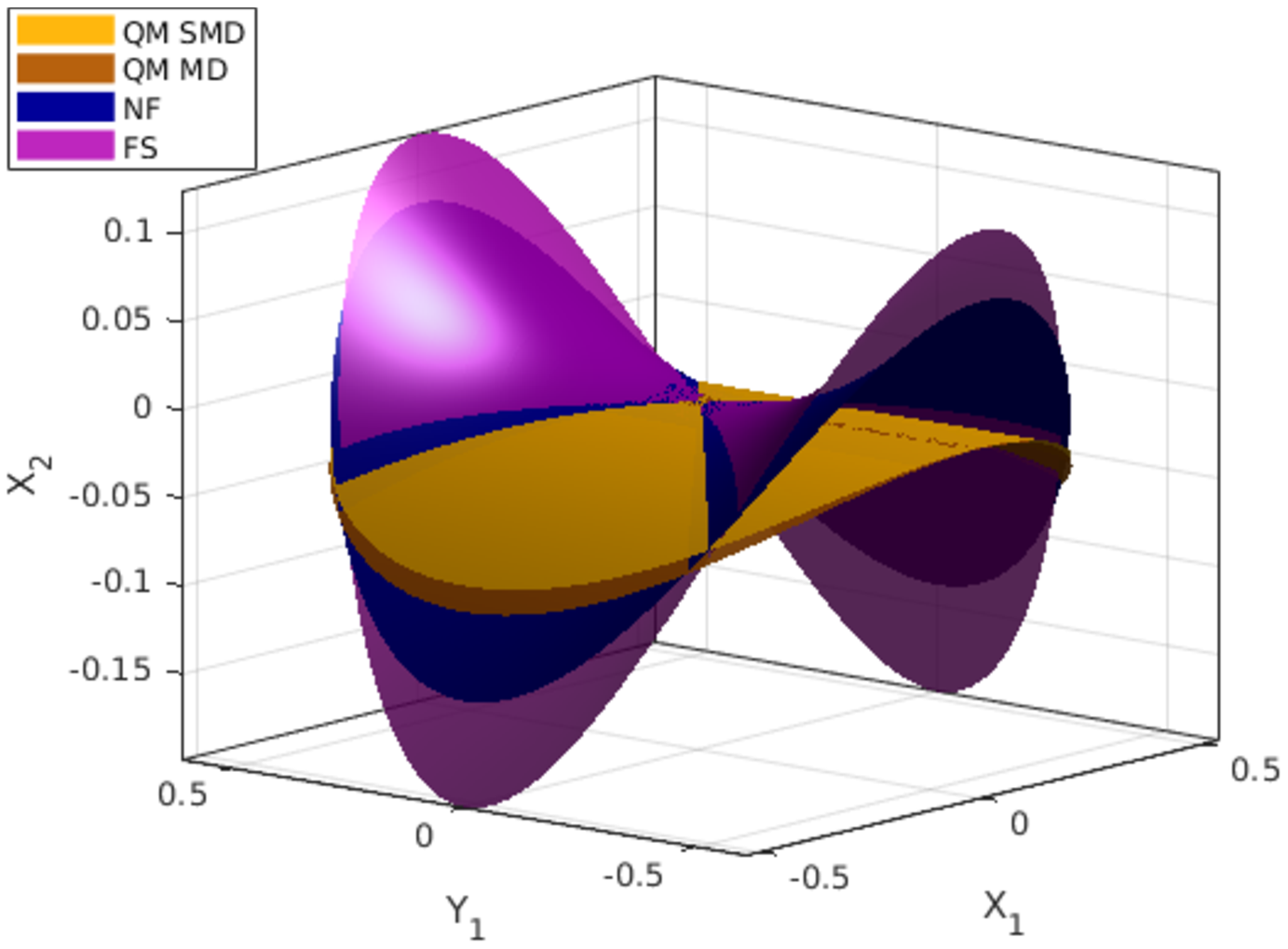}
\caption{$\rho=2.5$.}
\label{fig:ex1manifold-b}
\end{subfigure}
\begin{subfigure}[b]{.32\textwidth}
\includegraphics[width=\textwidth]{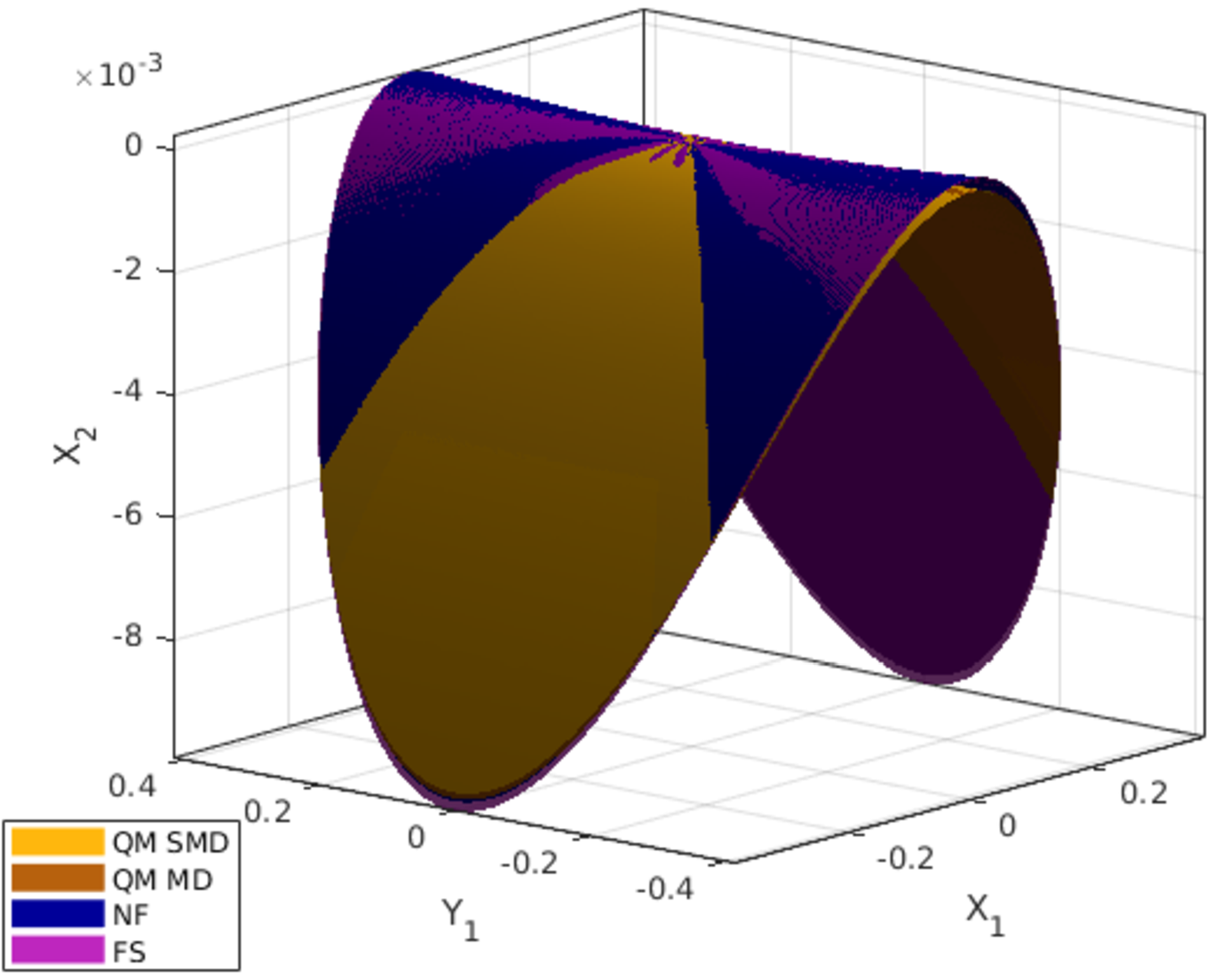}
\caption{$\rho=10$.}
\label{fig:ex1manifold-c}
\end{subfigure}
\caption{Comparison of manifolds in phase space for the first example, and for three different values of $\rho=\omega_2/\omega_1$. The exact NNM, represented in violet (full system solution: FS), is compared to the reduction manifolds obtained by QM from MDs (dark orange), from SMDs (yellow), and normal form  (blue). (a) $\rho=1.25$, (b) $\rho=2.5$, (c) $\rho=10$ with slow/fast assumption fulfilled.}
\label{fig:ex1manifold}
\end{figure*}

Fig.~\ref{fig:ex1manifold} shows the geometry of the manifolds obtained for this first system, when one increases the values of $\rho$ so as to meet the slow/fast assumption. One can remark that the reduced subspaces produced by the quadratic manifold method don't show a dependence on the velocity. Increasing the values of $\rho$ it is observed that the real manifold obtained from the full system loses this velocity dependence so that this approximation is less and less wrong. On the other hand, the manifold produced by normal form has two important advantages: it is an invariant manifold of the full system by construction, and it has this velocity dependence, hence allowing for a correct prediction of the reduction subspace, whatever the value of $\rho$. As a matter of fact the only limitation of the normal form approach is that it relies on a Taylor expansion, so that for large amplitudes, the solution departs from the exact manifold. But in any case the correct invariant subspace is approximated. \red{As already remarked, due to the fact that only quadratic couplings are present between master and slave coordinates, the manifolds shown in Fig.~\ref{fig:ex1manifold} for the normal form are  obtained thanks to the second-order expansion, the third-order terms being all equal to zero.}

Fig.~\ref{fig:ex1manifold-a} shows also that the quadratic manifold produced by MD encounters a problem near the 1:1 resonance, which is here underlined since $\rho$ has been selected close to 1. Comparison with a full order solution clearly shows that this is a failure of the method. On the other hand, Fig.~\ref{fig:ex1manifold-c} shows that when the slow/fast assumption is verified, then all methods converge to the same reduced subspace, in line with the theoretical results.

We now turn to the prediction given on the backbone curves. First of all, one can compare the values of the $\Gamma$ coefficients dictating the type of nonlinearity. Eqs.~\eqref{eq:zifullgamma} have thus been rewritten for the present two-dofs system and now read, as a function of the ratio $\rho=\omega_2/\omega_1$:
\begin{subequations}\label{eq:gammasysteme1}
\begin{align}
&\Gamma_\text{MD}=
\frac{3D}{8}-
\frac{\bar{G}^2(3\rho^2-2)\rho^2}
{4(\rho^2-1)^2},
\\
&\Gamma_\text{SMD}=
\frac{3D}{8}-
\frac{\bar{G}^2(3\rho^2+4)}
{4\rho^2},
\\
&\Gamma_\text{NF}=
\frac{3D}{8}-
\frac{\bar{G}^2(3\rho^2-8)}
{4(\rho^2-4)}.
\end{align}
\end{subequations}

\begin{figure}[h!]
\centering
\begin{subfigure}{.46\textwidth}
\includegraphics{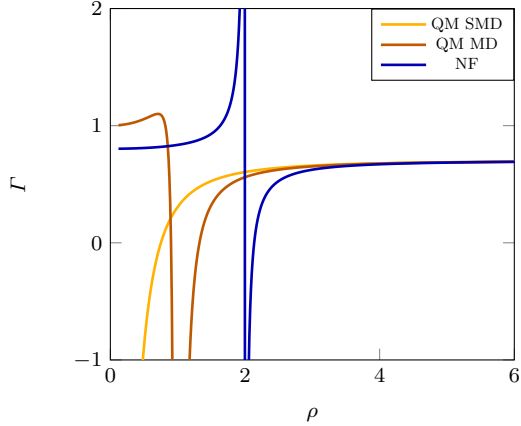}
\end{subfigure}
\caption{Values of the coefficient $\Gamma$ dictating the hardening/softening behaviour for the first two-dofs system. Comparison of $\Gamma_\text{MD}$, $\Gamma_\text{SMD}$ and $\Gamma_\text{NF}$, given respectively by QM with MDs, with SMDs, and normal form, Eqs.~\eqref{eq:gammasysteme1}, and for varying $\rho=\omega_2/\omega_1$ ratio.}
\label{fig:gamma_flat2dofs}
\end{figure}
These values are represented in Fig.~\ref{fig:gamma_flat2dofs}, which shows important similarities with Fig.~\ref{fig:gamma_flat}. Again the same divergent behaviours are observed, and the convergence of all methods for $\rho > 4$ is clearly observed. \red{To be more quantitative, the relative difference between $\Gamma_\text{MD}$ and $\Gamma_\text{NF}$ is  $5\%$ for $\rho=2.95$ and $1\%$ for $\rho=3.93$. On the other hand, the difference between $\Gamma_\text{SMD}$ and $\Gamma_\text{NF}$ is $5\%$ for $\rho=3.06$, and $1\%$ for $\rho=4.18$, underlining clearly that $\rho\in[3,4]$ has to be understood as a transition zone.} For very small values of $\rho$, the quadratic manifold based on SMD will predict incorrect result with a softening behaviour. Also, after its failure at $\rho=1$, the MD method will also produce an incorrect prediction with a softening behaviour. 


%
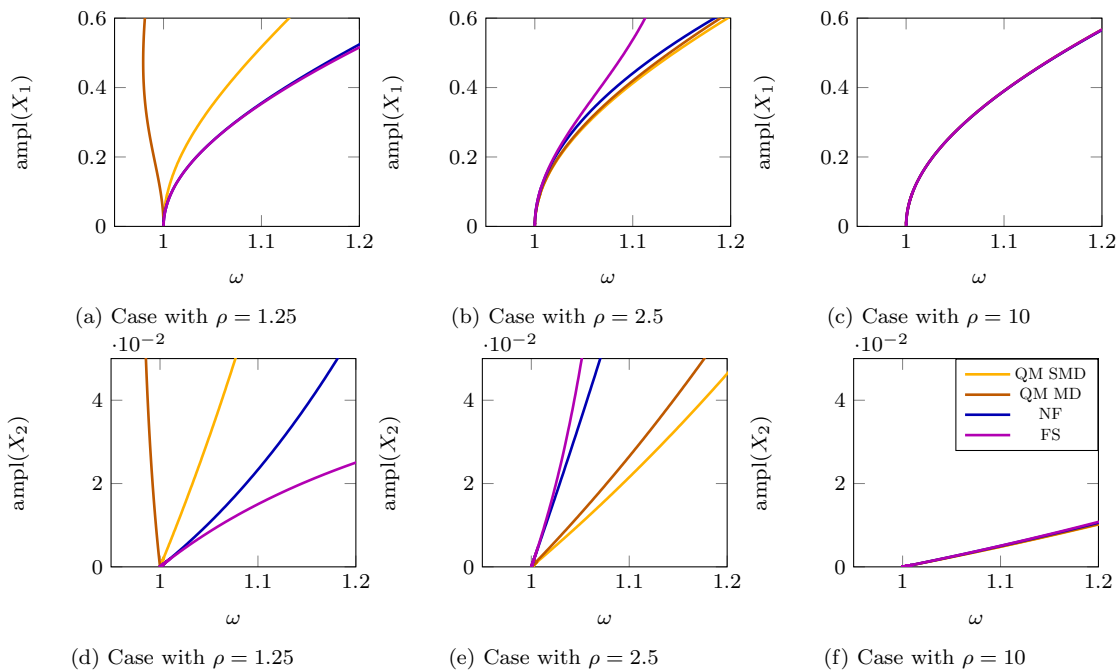
\begin{figure*}[h!]
\centering
\begin{subfigure}{.32\textwidth}
\begin{tikzpicture}
\begin{axis}[
unbounded coords=jump,
xlabel=$\omega$,
ylabel=ampl($X_1$),
width=\textwidth,
height=.9\textwidth,
domain=0:2,
xmin=.95, 
xmax=1.2,
ymin=0,ymax=.6,
legend style = {at={(1.,0.)},anchor=south east,font=\normalsize,nodes={scale=0.65, transform shape},},
x tick label style={rotate=0,anchor=north}
]
\addplot[Color_SMD,line width =1pt] 
table [
x=omega, 
y=x1_st, 
col sep=comma,
mark=none,]
{data_QM_case1_manlab.csv};
\addplot[Color_MD,line width =1pt] 
table [
x=omega, 
y=x1_st, 
col sep=comma,
mark=none,]
{data_QM_MD_case1_manlab.csv};
\addplot[Color_NF,line width =1pt] 
table 
[x=omega, 
y=x1_st, 
col sep=comma,
mark=none,] 
{data_NF_case1_manlab.csv};
\addplot[Color_FS,line width =1pt] 
table [
x=omega, 
y=x1_st, 
col sep=comma,
mark=none,]
{data_FS_case1_manlab.csv};
\end{axis}
\end{tikzpicture}
\caption{Case with $\rho=1.25$}
\label{fig:Backbones_flatdofx1_a}
\end{subfigure}
\begin{subfigure}{.32\textwidth}
\begin{tikzpicture}
\begin{axis}[
unbounded coords=jump,
xlabel=$\omega$,
ylabel=ampl($X_1$),
width=\textwidth,
height=.9\textwidth,
domain=0:2,
xmin=.95,  
xmax=1.2,
ymin=0,ymax=.6,
legend style = {at={(1.,0.)},anchor=south east,font=\normalsize,nodes={scale=0.65, transform shape},},
x tick label style={rotate=0,anchor=north}
]
\addplot[Color_SMD,line width =1pt] 
table [
x=omega, 
y=x1_st, 
col sep=comma,
mark=none,]
{data_QM_case2_manlab.csv};
\addplot[Color_MD,line width =1pt] 
table [
x=omega, 
y=x1_st, 
col sep=comma,
mark=none,]
{data_QM_MD_case2_manlab.csv};
\addplot[Color_NF,line width =1pt] 
table 
[x=omega, 
y=x1_st, 
col sep=comma,
mark=none,] 
{data_NF_case2_manlab.csv};
\addplot[Color_FS,line width =1pt] 
table [
x=omega, 
y=x1_st, 
col sep=comma,
mark=none,]
{data_FS_case2_manlab.csv};
\end{axis}
\end{tikzpicture}
\caption{Case with $\rho=2.5$}
\label{fig:Backbones_flatdofx1_b}
\end{subfigure}
\begin{subfigure}{.32\textwidth}
\begin{tikzpicture}
\begin{axis}[
unbounded coords=jump,
xlabel=$\omega$,
ylabel=ampl($X_1$),
width=\textwidth,
height=.9\textwidth,
domain=0:2,
xmin=.95,  
xmax=1.2,
ymin=0,ymax=.6,
legend style = {at={(1.,0.)},anchor=south east,font=\normalsize,nodes={scale=0.65, transform shape},},
x tick label style={rotate=0,anchor=north}
]
\addplot[Color_SMD,line width =1pt] 
table [
x=omega, 
y=x1_st, 
col sep=comma,
mark=none,]
{data_QM_case4_manlab.csv};
\addplot[Color_MD,line width =1pt] 
table [
x=omega, 
y=x1_st, 
col sep=comma,
mark=none,]
{data_QM_MD_case4_manlab.csv};
\addplot[Color_NF,line width =1pt] 
table 
[x=omega, 
y=x1_st, 
col sep=comma,
mark=none,] 
{data_NF_case4_manlab.csv};
\addplot[Color_FS,line width =1pt] 
table [
x=omega, 
y=x1_st, 
col sep=comma,
mark=none,]
{data_FS_case4_manlab.csv};
\end{axis}
\end{tikzpicture}
\caption{Case with $\rho=10$}
\label{fig:Backbones_flatdofx1_c}
\end{subfigure}
\begin{subfigure}{.32\textwidth}
\begin{tikzpicture}
\begin{axis}[
unbounded coords=jump,
xlabel=$\omega$,
ylabel=ampl($X_2$),
width=\textwidth,
height=.9\textwidth,
domain=0:2,
xmin=.95,
xmax=1.2,
ymin=0, 
ymax=.05,
legend style = {at={(1.,0.)},anchor=south east,font=\normalsize,nodes={scale=0.65, transform shape}},
x tick label style={rotate=0,anchor=north}
]
\addplot[Color_SMD,line width =1pt] 
table [
x=omega, 
y=x2_st, 
col sep=comma,
mark=none,]
{data_QM_case1_manlab.csv};
\addplot[Color_MD,line width =1pt] 
table [
x=omega, 
y=x2_st, 
col sep=comma,
mark=none,]
{data_QM_MD_case1_manlab.csv};
\addplot[Color_NF,line width =1pt] 
table 
[x=omega, 
y=x2_st, 
col sep=comma,
mark=none,] 
{data_NF_case1_manlab.csv};
\addplot[Color_FS,line width =1pt] 
table [
x=omega, 
y=x2_st, 
col sep=comma,
mark=none,]
{data_FS_case1_manlab.csv};
\end{axis}
\end{tikzpicture}
\caption{Case with $\rho=1.25$}
\label{fig:Backbones_flatdofx1_a}
\end{subfigure}
\begin{subfigure}{.32\textwidth}
\begin{tikzpicture}
\begin{axis}[
unbounded coords=jump,
xlabel=$\omega$,
ylabel=ampl($X_2$),
width=\textwidth,
height=.9\textwidth,
domain=0:2,
xmin=.95, 
xmax=1.2,
ymin=0, 
ymax=.05,
legend style = {at={(1.,0.)},anchor=south east,font=\normalsize,nodes={scale=0.65, transform shape}},
x tick label style={rotate=0,anchor=north}
]
\addplot[Color_SMD,line width =1pt] 
table [
x=omega, 
y=x2_st, 
col sep=comma,
mark=none,]
{data_QM_case2_manlab.csv};
\addplot[Color_MD,line width =1pt] 
table [
x=omega, 
y=x2_st, 
col sep=comma,
mark=none,]
{data_QM_MD_case2_manlab.csv};
\addplot[Color_NF,line width =1pt] 
table 
[x=omega, 
y=x2_st, 
col sep=comma,
mark=none,] 
{data_NF_case2_manlab.csv};
\addplot[Color_FS,line width =1pt] 
table [
x=omega, 
y=x2_st, 
col sep=comma,
mark=none,]
{data_FS_case2_manlab.csv};
\end{axis}
\end{tikzpicture}
\caption{Case with $\rho=2.5$}
\label{fig:Backbones_flatdofx1_b}
\end{subfigure}
\begin{subfigure}{.32\textwidth}
\begin{tikzpicture}
\begin{axis}[
unbounded coords=jump,
xlabel=$\omega$,
ylabel=ampl($X_2$),
width=\textwidth,
height=.9\textwidth,
domain=0:2,
xmin=.95, 
xmax=1.2,
ymin=0, 
ymax=.05,
legend style = {at={(1.,1.)},anchor=north east,font=\normalsize,nodes={scale=0.65, transform shape}},
x tick label style={rotate=0,anchor=north}
]
\addplot[Color_SMD,line width =1pt] 
table [
x=omega, 
y=x2_st, 
col sep=comma,
mark=none,]
{data_QM_case4_manlab.csv};
\addplot[Color_MD,line width =1pt] 
table [
x=omega, 
y=x2_st, 
col sep=comma,
mark=none,]
{data_QM_MD_case4_manlab.csv};
\addplot[Color_NF,line width =1pt] 
table 
[x=omega, 
y=x2_st, 
col sep=comma,
mark=none,] 
{data_NF_case4_manlab.csv};
\addplot[Color_FS,line width =1pt] 
table [
x=omega, 
y=x2_st, 
col sep=comma,
mark=none,]
{data_FS_case4_manlab.csv};
\legend{QM SMD,QM MD,NF,FS}
\end{axis}
\end{tikzpicture}
\caption{Case with $\rho=10$}
\label{fig:Backbones_flatdofx1_c}
\end{subfigure}
\caption{First mode backbone curves as a function of modal amplitude $X_1$ (first row), and $X_2$ (second row), and for different values of $\rho = \omega_2/\omega_1$. Comparisons between the exact solution (FS: full system, violet), that predicted by QM with MDs (dark orange), SMDs (yellow) and normal form (NF, blue).}
\label{fig:bbflat}
\end{figure*}

Fig.~\ref{fig:bbflat} shows the backbone curves obtained from the reduced-order dynamics and compared to that obtained from the full system. The comparison is drawn on the main modal amplitude $X_1$, which shows the largest values (first row), but also on the slave coordinate $X_2$ (second row). The first case selected, just after the 1:1 resonance with $\rho=1.25$, shows, as envisioned in Fig.~\ref{fig:gamma_flat2dofs}, that the QM produced from MD can be very wrong in this case and predict at first order a softening behaviour. When $\rho=2.5$, the three methods predicts a very similar behaviour and are almost undistinguishable. One can note that for large amplitude, the full system solution is less and less hardening. This is probably a consequence of the vicinity of the 2:1 internal resonance. Since $\omega_2 = 2.5\omega_1$ and the behaviour is hardening, the nonlinear frequency tends to approach the 2:1 ratio at higher amplitudes, which could explain this particular behaviour of the full system solution. Finally, for $\rho=10$, the three methods give the same predictions which are fully aligned with the full system.

The conclusion on this first example with simple nonlinearities are in the line of the theoretical results, since all methods tends to perform well in the limit of the slow/fast assumption, again estimated as a ratio of 4 between the eigenfrequencies of the master and slave mode. On the other hand, when this assumption is not fulfilled, the quadratic manifold is not reliable and can produce incorrect predictions, in contrary to the normal form approach, that gives a correct ROM up to the third-order, whatever the link between slave and master coordinates. These results explain also why the application of modal derivatives on slender structures that are flat and symmetric produce accurate results. Indeed, slenderness is fulfilled when $\rho$ is larger than 40, and our numerical experiments show that the slow/fast assumption can be considered as valid as soon as $\rho > 4$.

\subsection{A two-dof model representative of a shell structure}

\subsubsection{Equations of motion}

\red{In this section, a system composed of a mass connected to two springs representing geometric nonlinearity, is selected. This system has been used in a number of studies so that numerous results are already present in the literature, the interested reader is referred to~\cite{touze03-NNM} for the derivation of the equation of motions specifying the behaviour of the springs, and to~\cite{touze03-NNM,TOUZE:JSV:2006,RENSON2016,BreunungHaller18} for different results already published on this example system.} The equations of motion read:
\begin{equation}\label{eq:2dofscoche}
\begin{aligned}
\ddot{X}_1+\omega_1^2X_1+\frac{\omega_1^2}{2}(3X_1^2+X_2^2)+\omega_2^2X_1X_2+\frac{\omega_1^2+\omega_2^2}{2}X_1(X_1^2+X_2^2)=0,\\
\ddot{X}_2+\omega_2^2X_2+\frac{\omega_2^2}{2}(3X_2^2+X_1^2)+\omega_1^2X_1X_2+\frac{\omega_1^2+\omega_2^2}{2}X_2(X_1^2+X_2^2)=0.\\
\end{aligned}
\end{equation}
As compared to the previous example, this system has all quadratic nonlinear terms present in the equations of motion, and all the nonlinear coefficients are expressed directly from the two eigenvalues $\omega_1$ and $\omega_2$, so that the problem has only two parameters. Note that this model is not derived from a continuous shell structure like the previous example was derived from the \vonkar beam equations, however it is known that curved structures display strong   quadratic couplings that are found in this system. Moreover, the results will show that this model is sufficient to show important departures between the three tested methods, which are due to the way the quadratic terms are processed. 

\subsubsection{Results}

As for the preceding example, comparisons are drawn out on the geometry of the manifolds and the backbone curves. Numerical continuation is used to solve out the different systems and compare their outcomes. The eigenfrequency ratio $\rho=\omega_2/\omega_1$ is also used and the same values, namely 1.25, 2.5 and 10 are selected to observe the differences between the methods when tending to fulfil the slow/fast assumption. In the computation, $\omega_1=1$ in all cases so that one simply have $\omega_2=\rho$.

\begin{figure*}[h!]\centering
\begin{subfigure}[b]{.32\textwidth}
\includegraphics[width=\textwidth]{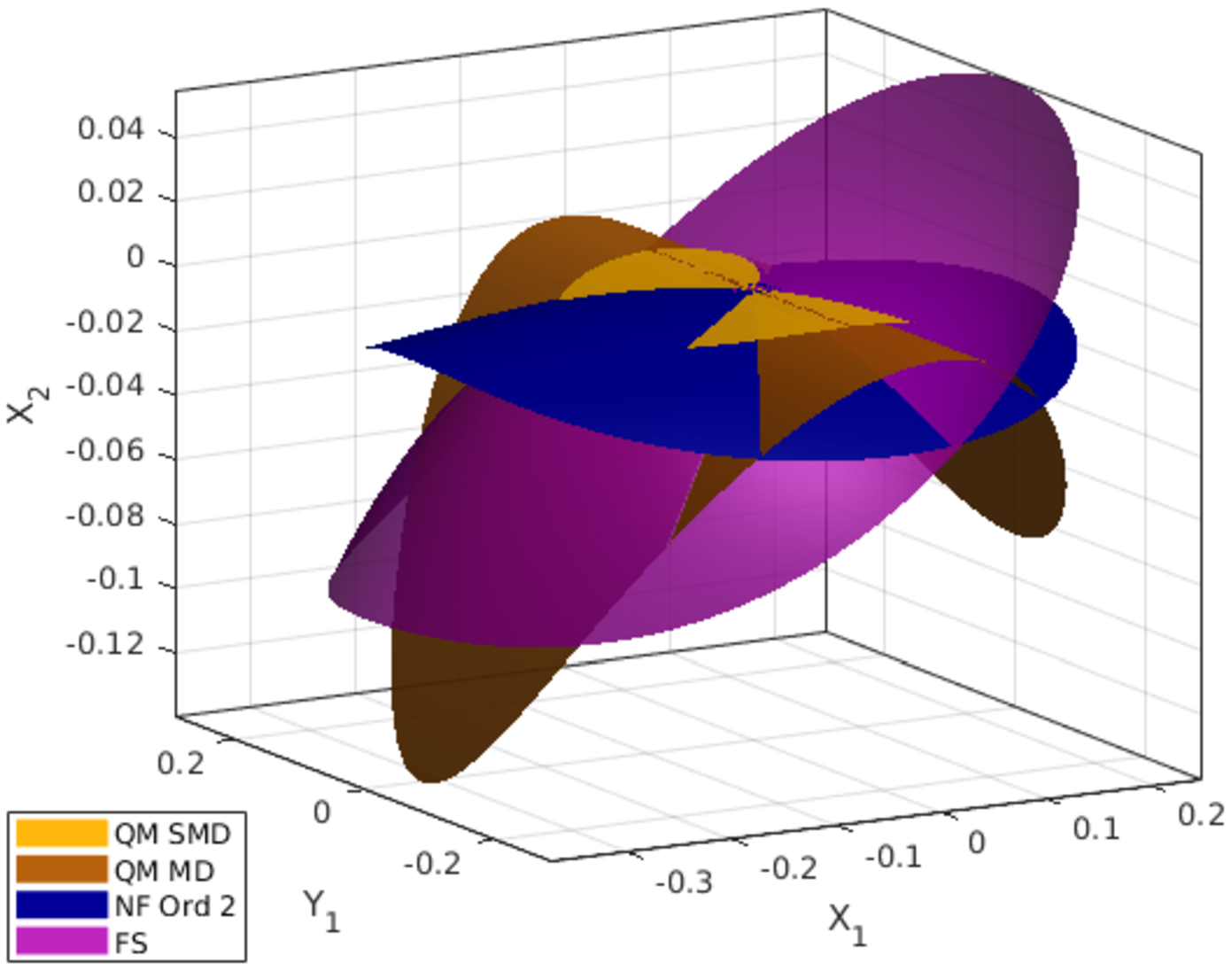}
\caption{Case with $\rho=1.25$.}
\label{fig:manif2dofshell-a}
\end{subfigure}
\begin{subfigure}[b]{.32\textwidth}
\includegraphics[width=\textwidth]{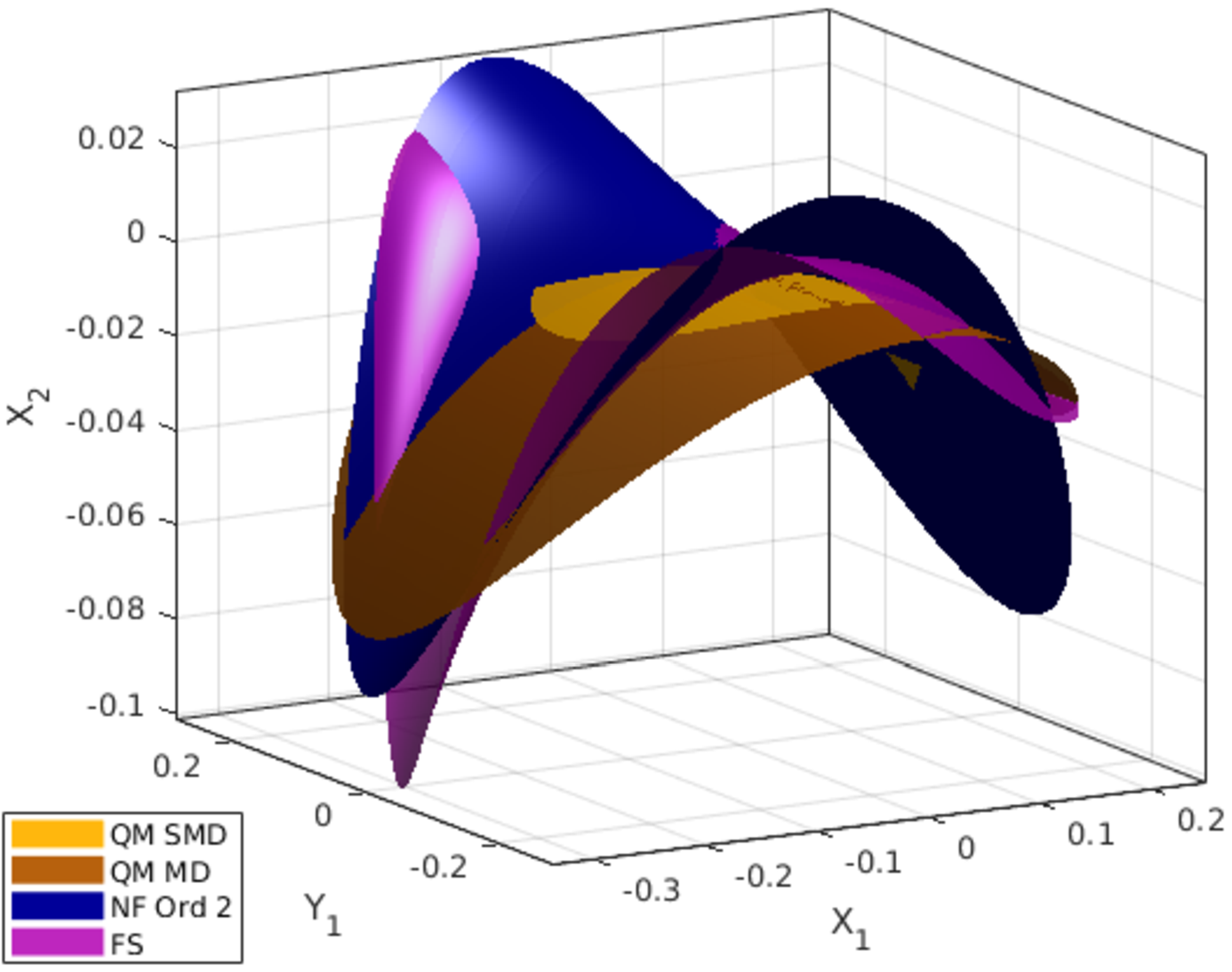}
\caption{Case with $\rho=2.5$.}
\label{fig:manif2dofshell-b}
\end{subfigure}
\begin{subfigure}[b]{.32\textwidth}
\includegraphics[width=\textwidth]{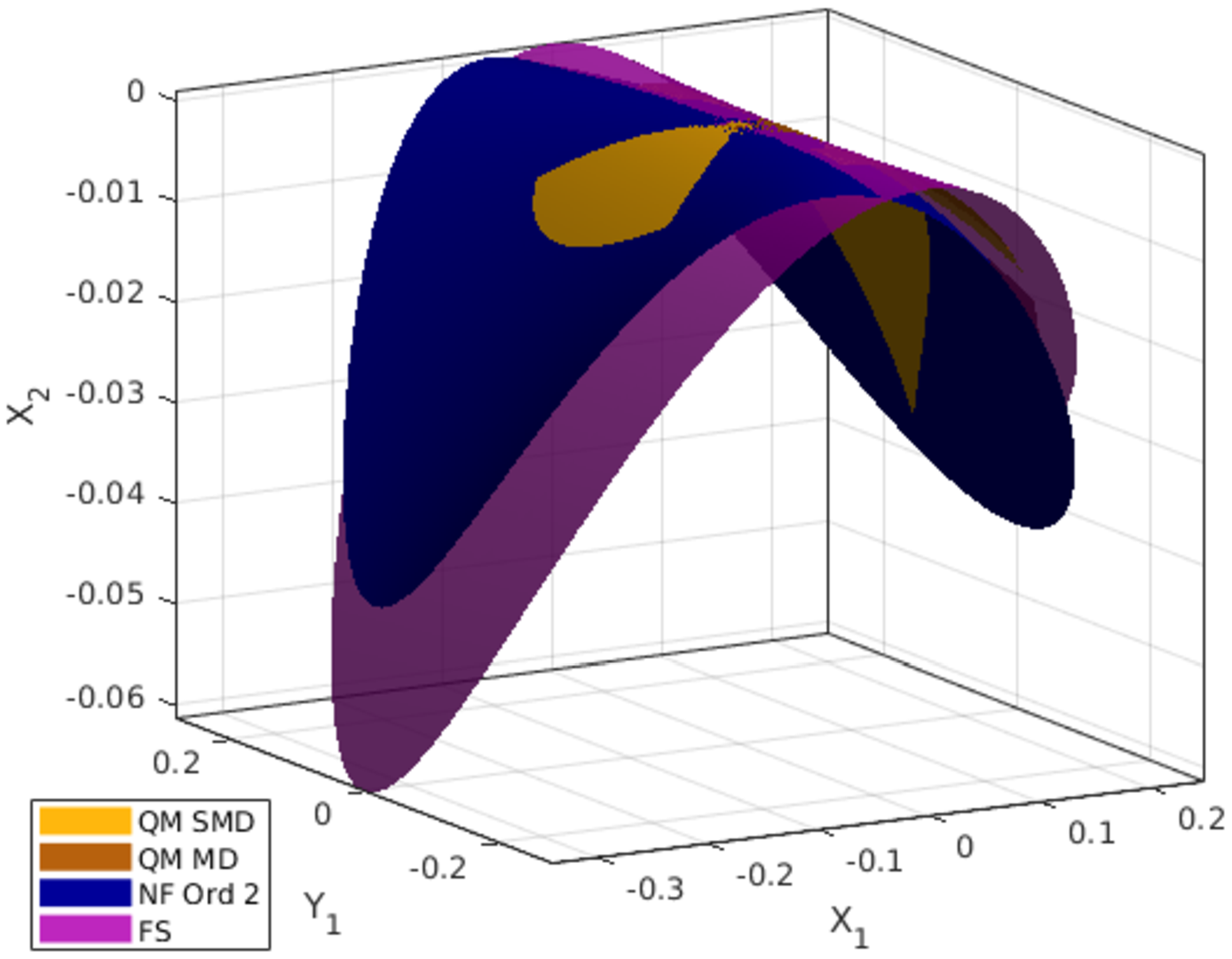}
\caption{Case with $\rho=10$.}
\label{fig:manif2dofshell-c}
\end{subfigure}
\begin{subfigure}[b]{.32\textwidth}
\includegraphics[width=\textwidth]{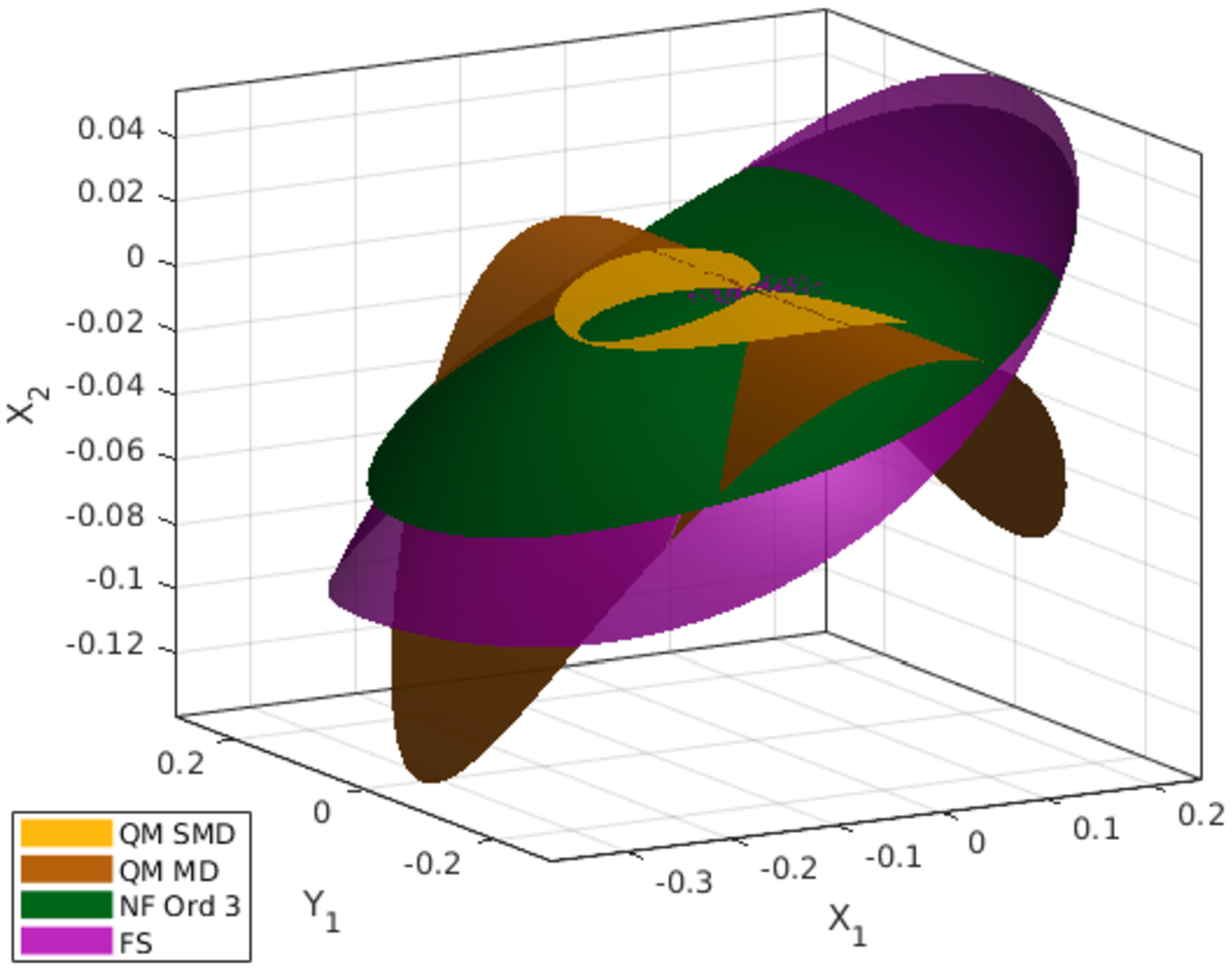}
\caption{Case with $\rho=1.25$.}
\label{fig:manif2dofshell-d}
\end{subfigure}
\begin{subfigure}[b]{.32\textwidth}
\includegraphics[width=\textwidth]{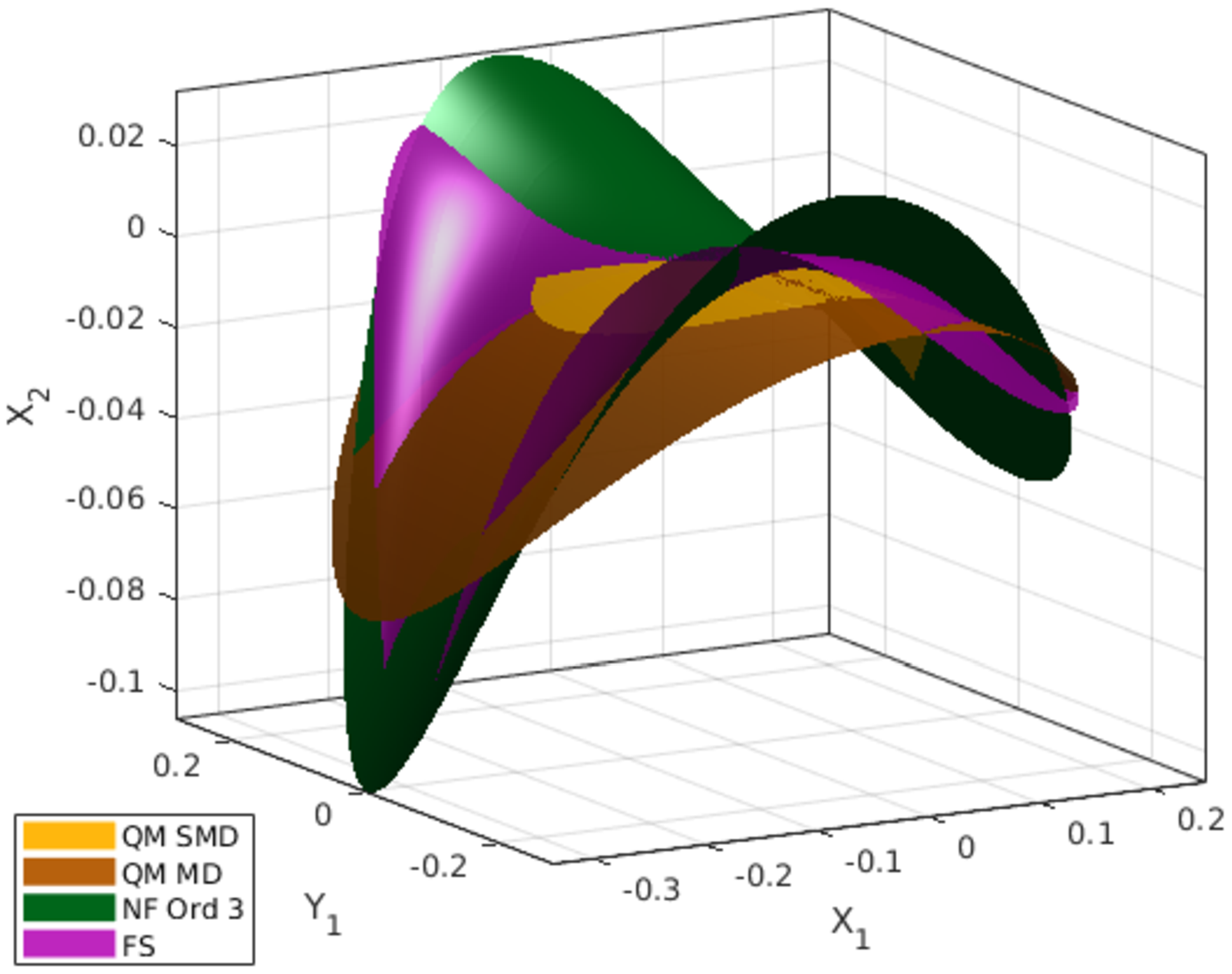}
\caption{Case with $\rho=2.5$.}
\label{fig:manif2dofshell-e}
\end{subfigure}
\begin{subfigure}[b]{.32\textwidth}
\includegraphics[width=\textwidth]{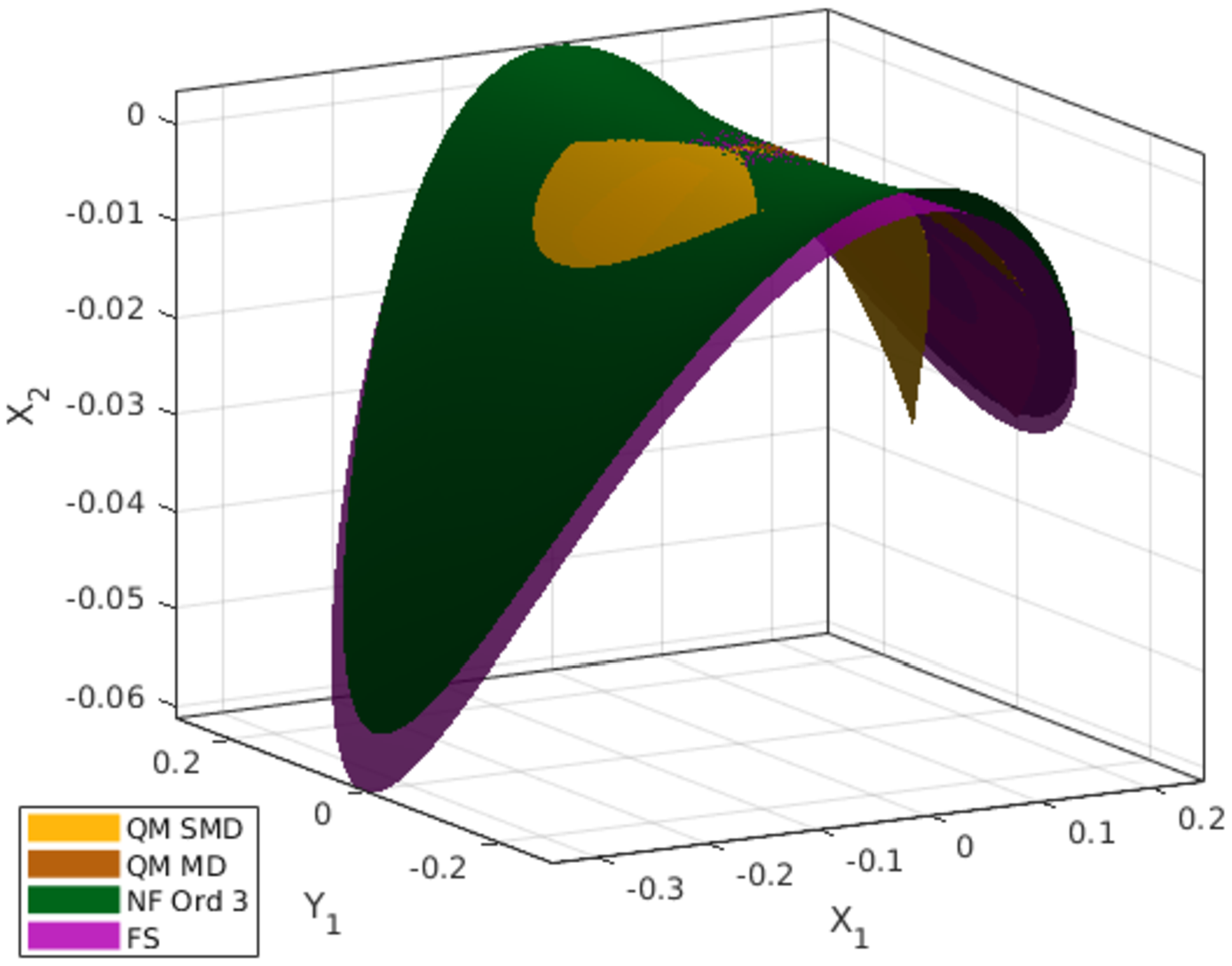}
\caption{Case with $\rho=10$.}
\label{fig:manif2dofshell-f}
\end{subfigure}
\caption{Comparison of manifolds in phase space for the second two-dofs example, and for three different values of $\rho=\omega_2/\omega_1$. The exact NNM, represented in violet (full system solution: FS), is compared to the reduction manifolds obtained by QM from MDs (dark orange), from SMDs (yellow), \red{and normal form up to the second order (blue manifold in the first line, plots a-b-c) and third order (green manifolds, second line in plots d-e-f) are given}. (a-d) $\rho=1.25$, (b-e) $\rho=2.5$, (c-f) $\rho=10$ with slow/fast assumption fulfilled. In all cases $\omega_1=1$.}
\label{fig:manif2dofshell}
\end{figure*}
%
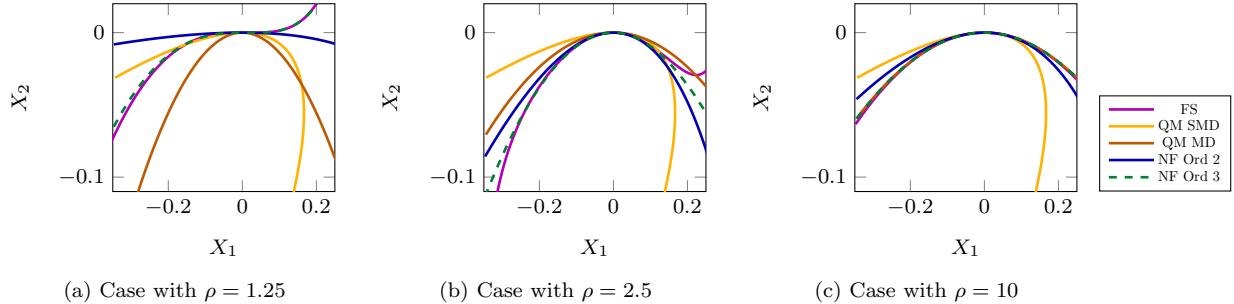
\begin{figure*}[h!]
\centering
\begin{subfigure}{.3\textwidth}
\begin{tikzpicture}
\begin{axis}[
unbounded coords=jump,
xlabel=$X_1$,
ylabel=$X_2$,
width=\textwidth,
height=.9\textwidth,
ytick = {-0.1,0},
ymin=-0.11,ymax=.02,
xmin=-0.35,xmax=.25,
legend style = {at={(1.,0.)},anchor=south east,font=\large,nodes={scale=0.5, transform shape},},
x tick label style={rotate=0,anchor=north}
]
\addplot[Color_FS,line width =1pt] 
table [
x=X1_cut,
y=X2_cut, 
col sep=comma,
mark=none,]
{data_FS_case1_manlab_manifolds.csv};
\addplot[Color_SMD,line width =1pt] 
table [
x=X1_cut,
y=X2_cut, 
col sep=comma,
mark=none,]
{data_QM_case1_manlab_manifolds.csv};
\addplot[Color_MD,line width =1pt] 
table [
x=X1_cut,
y=X2_cut, 
col sep=comma,
mark=none,]
{data_QM_MD_case1_manlab_manifolds.csv};
\addplot[Color_NF,line width =1pt] 
table [
x=X1_cut,
y=X2_cut, 
col sep=comma,
mark=none,] 
{data_NF_case1_manlab_manifolds.csv};
\addplot[Color_NFTO,dashed,line width =1pt] 
table [
x=X1_cut,
y=X2_cut, 
col sep=comma,
mark=none,] 
{data_NF_case1_manlab_manifolds_O3.csv};
%
\end{axis}
\end{tikzpicture}
\caption{Case with $\rho=1.25$}
\label{fig:Manifolds_shelldofx1_a}
\end{subfigure}
\quad
\begin{subfigure}{.3\textwidth}
\begin{tikzpicture}
\begin{axis}[
unbounded coords=jump,
xlabel=$X_1$,
ylabel=$X_2$,
width=\textwidth,
height=.9\textwidth,
ytick = {-0.1,0},
ymin=-0.11,ymax=.02,
xmin=-0.35,xmax=.25,
legend style = {at={(1.,0.)},anchor=south east,font=\large,nodes={scale=0.5, transform shape},},
x tick label style={rotate=0,anchor=north}
]
\addplot[Color_FS,line width =1pt] 
table [
x=X1_cut,
y=X2_cut, 
col sep=comma,
mark=none,]
{data_FS_case2_manlab_manifolds.csv};
\addplot[Color_SMD,line width =1pt] 
table [
x=X1_cut,
y=X2_cut, 
col sep=comma,
mark=none,]
{data_QM_case2_manlab_manifolds.csv};
\addplot[Color_MD,line width =1pt] 
table [
x=X1_cut,
y=X2_cut, 
col sep=comma,
mark=none,]
{data_QM_MD_case2_manlab_manifolds.csv};
\addplot[Color_NF,line width =1pt] 
table [
x=X1_cut,
y=X2_cut, 
col sep=comma,
mark=none,] 
{data_NF_case2_manlab_manifolds.csv};
\addplot[Color_NFTO,dashed,line width =1pt] 
table [
x=X1_cut,
y=X2_cut, 
col sep=comma,
mark=none,] 
{data_NF_case2_manlab_manifolds_O3.csv};
%
\end{axis}
\end{tikzpicture}
\caption{Case with $\rho=2.5$}
\label{fig:Manifolds_shelldofx1_b}
\end{subfigure}
\quad
\begin{subfigure}{.3\textwidth}
\begin{tikzpicture}
\begin{axis}[
unbounded coords=jump,
xlabel=$X_1$,
ylabel=$X_2$,
width=\textwidth,
height=.9\textwidth,
ytick = {-0.1,0},
ymin=-0.11,ymax=.02,
xmin=-0.35,xmax=.25,
legend style = {at={(1.1,0.)},anchor=south west,font=\large,nodes={scale=0.5, transform shape},},
x tick label style={rotate=0,anchor=north}
]
\addplot[Color_FS,line width =1pt] 
table [
x=X1_cut,
y=X2_cut, 
col sep=comma,
mark=none,]
{data_FS_case4_manlab_manifolds.csv};
\addplot[Color_SMD,line width =1pt] 
table [
x=X1_cut,
y=X2_cut, 
col sep=comma,
mark=none,]
{data_QM_case4_manlab_manifolds.csv};
\addplot[Color_MD,line width =1pt] 
table [
x=X1_cut,
y=X2_cut, 
col sep=comma,
mark=none,]
{data_QM_MD_case4_manlab_manifolds.csv};
\addplot[Color_NF,line width =1pt] 
table [
x=X1_cut,
y=X2_cut, 
col sep=comma,
mark=none,] 
{data_NF_case4_manlab_manifolds.csv};
\addplot[Color_NFTO,dashed,line width =1pt] 
table [
x=X1_cut,
y=X2_cut, 
col sep=comma,
mark=none,] 
{data_NF_case4_manlab_manifolds_O3.csv};
\legend{FS,QM SMD,QM MD,NF Ord 2,NF Ord 3}
\end{axis}
\end{tikzpicture}
\caption{Case with $\rho=10$}
\label{fig:Manifolds_shelldofx1_c}
\end{subfigure}
\caption{First mode invariant manifolds cut on the $Y_1=0$ plane, evaluated with the quadratic manifold method (QM)  (either with MD in dark orange, and SMD in yellow) and normal form (NF) approach, \red{where the distinction between NF up to second order (blue line) and third-order (dashed green line) is reported},  and compared to the numerical solution obtained with the full system (FS). In all cases $\omega_1=1$.}\label{fig:manifoldcut2d}
\end{figure*}
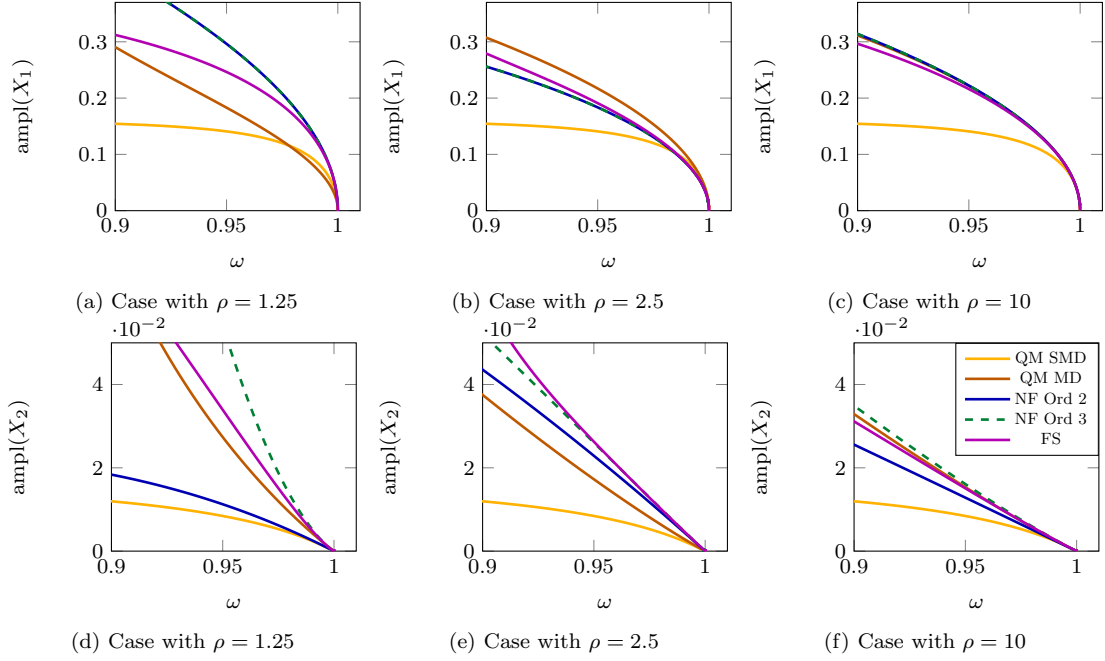
\begin{figure*}[h!]
\centering
\begin{subfigure}{.32\textwidth}
\begin{tikzpicture}
\begin{axis}[
unbounded coords=jump,
xlabel=$\omega$,
ylabel=ampl($X_1$),
width=\textwidth,
height=.9\textwidth,
domain=0:2,
xmin=.9, 
ymin=0,ymax=.37,
legend style = {at={(1.,1.)},anchor=north east,font=\normalsize,nodes={scale=0.65, transform shape},},
x tick label style={rotate=0,anchor=north}
]
\addplot[Color_SMD,line width =1pt] 
table [
x=omega, 
y=x1_st, 
col sep=comma,
mark=none,]
{data_QM_case1_manlab_shell.csv};
\addplot[Color_MD,line width =1pt] 
table [
x=omega, 
y=x1_st, 
col sep=comma,
mark=none,]
{data_QM_MD_case1_manlab_shell.csv};
\addplot[Color_NF,line width =1pt] 
table 
[x=omega, 
y=x1_st, 
col sep=comma,
mark=none,] 
{data_NF_case1_manlab_shell.csv};
\addplot[Color_NFTO,dashed,line width =1pt] 
table 
[x=omega, 
y=x1_st, 
col sep=comma,
mark=none,] 
{data_NF_case1_manlab_shell_O3.csv};
\addplot[Color_FS,line width =1pt] 
table [
x=omega, 
y=x1_st, 
col sep=comma,
mark=none,]
{data_FS_case1_manlab_shell.csv};
\end{axis}
\end{tikzpicture}
\caption{Case with $\rho=1.25$}
\label{fig:Backbones_shelldofx1_a}
\end{subfigure}
\begin{subfigure}{.32\textwidth}
\begin{tikzpicture}
\begin{axis}[
unbounded coords=jump,
xlabel=$\omega$,
ylabel=ampl($X_1$),
width=\textwidth,
height=.9\textwidth,
domain=0:2,
xmin=.9,  
ymin=0,ymax=.37,
legend style = {at={(1.,1.)},anchor=north east,font=\normalsize,nodes={scale=0.65, transform shape},},
x tick label style={rotate=0,anchor=north}
]
\addplot[Color_SMD,line width =1pt] 
table [
x=omega, 
y=x1_st, 
col sep=comma,
mark=none,]
{data_QM_case2_manlab_shell.csv};
\addplot[Color_MD,line width =1pt] 
table [
x=omega, 
y=x1_st, 
col sep=comma,
mark=none,]
{data_QM_MD_case2_manlab_shell.csv};
\addplot[Color_NF,line width =1pt] 
table 
[x=omega, 
y=x1_st, 
col sep=comma,
mark=none,] 
{data_NF_case2_manlab_shell.csv};
\addplot[Color_NFTO,dashed,line width =1pt] 
table 
[x=omega, 
y=x1_st, 
col sep=comma,
mark=none,] 
{data_NF_case2_manlab_shell_O3.csv};
\addplot[Color_FS,line width =1pt] 
table [
x=omega, 
y=x1_st, 
col sep=comma,
mark=none,]
{data_FS_case2_manlab_shell.csv};
\end{axis}
\end{tikzpicture}
\caption{Case with $\rho=2.5$}
\label{fig:Backbones_shelldofx1_b}
\end{subfigure}
\begin{subfigure}{.32\textwidth}
\begin{tikzpicture}
\begin{axis}[
unbounded coords=jump,
xlabel=$\omega$,
ylabel=ampl($X_1$),
width=\textwidth,
height=.9\textwidth,
domain=0:2,
xmin=.9,  
ymin=0,ymax=.37,
legend style = {at={(1.,1.)},anchor=north east,font=\normalsize,nodes={scale=0.65, transform shape},},
x tick label style={rotate=0,anchor=north}
]
\addplot[Color_SMD,line width =1pt] 
table [
x=omega, 
y=x1_st, 
col sep=comma,
mark=none,]
{data_QM_case4_manlab_shell.csv};
\addplot[Color_MD,line width =1pt] 
table [
x=omega, 
y=x1_st, 
col sep=comma,
mark=none,]
{data_QM_MD_case4_manlab_shell.csv};
\addplot[Color_NF,line width =1pt] 
table 
[x=omega, 
y=x1_st, 
col sep=comma,
mark=none,] 
{data_NF_case4_manlab_shell.csv};
\addplot[Color_NFTO,dashed,line width =1pt] 
table 
[x=omega, 
y=x1_st, 
col sep=comma,
mark=none,] 
{data_NF_case4_manlab_shell_O3.csv};
\addplot[Color_FS,line width =1pt] 
table [
x=omega, 
y=x1_st, 
col sep=comma,
mark=none,]
{data_FS_case4_manlab_shell.csv};
\end{axis}
\end{tikzpicture}
\caption{Case with $\rho=10$}
\label{fig:Backbones_shelldofx1_c}
\end{subfigure}
\begin{subfigure}{.32\textwidth}
\begin{tikzpicture}
\begin{axis}[
unbounded coords=jump,
xlabel=$\omega$,
ylabel=ampl($X_2$),
width=\textwidth,
height=.9\textwidth,
domain=0:2,
xmin=.9,
ymin=0, 
ymax=.05,
legend style = {at={(1.,1.)},anchor=north east,font=\normalsize,nodes={scale=0.65, transform shape}},
x tick label style={rotate=0,anchor=north}
]
\addplot[Color_SMD,line width =1pt] 
table [
x=omega, 
y=x2_st, 
col sep=comma,
mark=none,]
{data_QM_case1_manlab_shell.csv};
\addplot[Color_MD,line width =1pt] 
table [
x=omega, 
y=x2_st, 
col sep=comma,
mark=none,]
{data_QM_MD_case1_manlab_shell.csv};
\addplot[Color_NF,line width =1pt] 
table 
[x=omega, 
y=x2_st, 
col sep=comma,
mark=none,] 
{data_NF_case1_manlab_shell.csv};
\addplot[Color_NFTO,dashed,line width =1pt] 
table 
[x=omega, 
y=x2_st, 
col sep=comma,
mark=none,] 
{data_NF_case1_manlab_shell_O3.csv};
\addplot[Color_FS,line width =1pt] 
table [
x=omega, 
y=x2_st, 
col sep=comma,
mark=none,]
{data_FS_case1_manlab_shell.csv};
\end{axis}
\end{tikzpicture}
\caption{Case with $\rho=1.25$}
\label{fig:Backbones_shelldofx1_a}
\end{subfigure}
\begin{subfigure}{.32\textwidth}
\begin{tikzpicture}
\begin{axis}[
unbounded coords=jump,
xlabel=$\omega$,
ylabel=ampl($X_2$),
width=\textwidth,
height=.9\textwidth,
domain=0:2,
xmin=.9, 
ymin=0, 
ymax=.05,
legend style = {at={(1.,1.)},anchor=north east,font=\normalsize,nodes={scale=0.65, transform shape}},
x tick label style={rotate=0,anchor=north}
]
\addplot[Color_SMD,line width =1pt] 
table [
x=omega, 
y=x2_st, 
col sep=comma,
mark=none,]
{data_QM_case2_manlab_shell.csv};
\addplot[Color_MD,line width =1pt] 
table [
x=omega, 
y=x2_st, 
col sep=comma,
mark=none,]
{data_QM_MD_case2_manlab_shell.csv};
\addplot[Color_NF,line width =1pt] 
table 
[x=omega, 
y=x2_st, 
col sep=comma,
mark=none,] 
{data_NF_case2_manlab_shell.csv};
\addplot[Color_NFTO,dashed,line width =1pt] 
table 
[x=omega, 
y=x2_st, 
col sep=comma,
mark=none,] 
{data_NF_case2_manlab_shell_O3.csv};
\addplot[Color_FS,line width =1pt] 
table [
x=omega, 
y=x2_st, 
col sep=comma,
mark=none,]
{data_FS_case2_manlab_shell.csv};
\end{axis}
\end{tikzpicture}
\caption{Case with $\rho=2.5$}
\label{fig:Backbones_shelldofx1_b}
\end{subfigure}
\begin{subfigure}{.32\textwidth}
\begin{tikzpicture}
\begin{axis}[
unbounded coords=jump,
xlabel=$\omega$,
ylabel=ampl($X_2$),
width=\textwidth,
height=.9\textwidth,
domain=0:2,
xmin=.9, 
ymin=0, 
ymax=.05,
legend style = {at={(1.,1.)},anchor=north east,font=\normalsize,nodes={scale=0.65, transform shape}},
x tick label style={rotate=0,anchor=north}
]
\addplot[Color_SMD,line width =1pt] 
table [
x=omega, 
y=x2_st, 
col sep=comma,
mark=none,]
{data_QM_case4_manlab_shell.csv};
\addplot[Color_MD,line width =1pt] 
table [
x=omega, 
y=x2_st, 
col sep=comma,
mark=none,]
{data_QM_MD_case4_manlab_shell.csv};
\addplot[Color_NF,line width =1pt] 
table 
[x=omega, 
y=x2_st, 
col sep=comma,
mark=none,] 
{data_NF_case4_manlab_shell.csv};
\addplot[Color_NFTO,dashed,line width =1pt] 
table 
[x=omega, 
y=x2_st, 
col sep=comma,
mark=none,] 
{data_NF_case4_manlab_shell_O3.csv};
\addplot[Color_FS,line width =1pt] 
table [
x=omega, 
y=x2_st, 
col sep=comma,
mark=none,]
{data_FS_case4_manlab_shell.csv};
\legend{QM SMD,QM MD,NF Ord 2,NF Ord 3,FS}
\end{axis}
\end{tikzpicture}
\caption{Case with $\rho=10$}
\label{fig:Backbones_shelldofx1_c}
\end{subfigure}
\caption{First mode backbone curves for the second two-dofs example with quadratic nonlinearity, as a function of modal amplitude $X_1$ (first row), and $X_2$ (second row), and for different values of $\rho = \omega_2/\omega_1$. Comparisons between the exact solution (FS: full system, violet), that predicted by QM with MDs (dark orange), SMDs (yellow) and normal form (NF, blue, NF third order, dashed green).}\label{fig:BackbonesShell}
\end{figure*}

Fig.~\ref{fig:manif2dofshell} shows the geometry of the manifolds in phase space, as compared to the exact invariant manifold defining the first NNM of the system. The comment on the velocity dependence, already raised in the previous example, still holds: while for small values of $\rho$ the quadratic manifolds are not able to catch the correct curvature in this direction, for large values of $\rho$ the velocity dependence vanishes. Note that in all the three figures, the manifold produced by the SMD method has a smaller range in amplitude. This maximal range used for the representation has been fixed  from the frequency-amplitude relationships (see Fig.~\ref{fig:BackbonesShell}, when the nonlinear frequency has decreased of ten percent and reaches the value 0.9 --a softening behaviour is at hand in the considered cases--), so that all manifolds spans the same frequency range, but corresponds to different amplitudes.  This underlines in particular that even if the correct manifold is approximated, which is the case for $\rho=10$, the amplitude-frequency relationship may be not.

\red{Since the only difference between second- and third-order normal form can be appreciated from the nonlinear mapping and not the reduced dynamics, Fig.~\ref{fig:manif2dofshell} illustrates the case. In the first line, the manifold produced by the second-order normal form (in blue) is contrasted to the other methods, while the third-order is shown in the second line (in green). One can observe that the effect of retaining the cubic term is especially important for the smallest values of $\rho=1.25$, but then the differences between second- and third-order are barely visible.} Interestingly, this example also  shows that the quadratic manifolds produced by MD and SMD does not tend to the same geometries, even under the assumption of slow/fast dynamics. This can be appreciated in Fig.~\ref{fig:manif2dofshell}, but is more clearly evidenced in Fig.~\ref{fig:manifoldcut2d} where a section of the manifolds in space $(X_1,X_2)$ is shown, without the amplitude limit given by the frequency, used in the 3d plot.


Unlike Fig.~\ref{fig:manif2dofshell}, Fig.~\ref{fig:manifoldcut2d} has been directly obtained from the manifolds  expressions given by  Eq.~\eqref{eq:QMNLmap_modal_comp} for the QM approach, and \eqref{eq:transformNF}-\eqref{eq:invarman} for the normal form approach, by simply prescribing the values of $R_1$ and compute the resulting $(X_1,X_2)$ values. More specifically, let us underline the main difference between the MD and SMD method in this case. Using Eqs.~\eqref{eq:QMNLmap_modal_comp} with \eqref{eq:mdmodalexpr}, the reconstruction of $(X_1,X_2)$ from the QM method derived from MD reads:
\begin{subequations}\label{eq:QMMDcoche}
\begin{align}
&X_1=\q_1, \label{eq:QMMDcoche-a}\\
&X_2=-\frac{g^2_{11}}{\omega_2^2 - \omega_1^2}\q_1^2 = - \frac{\omega_2^2}{2(\omega_2^2 - \omega_1^2)}\q_1^2.\label{eq:QMMDcoche-b}
\end{align}
\end{subequations}
On the other hand, using SMD in the QM leads to:
\begin{subequations}\label{eq:QMSMDcoche}
\begin{align}
&X_1=\q_1  - \frac{g^1_{11}}{\omega_1^2}\q_1^2 = \q_1  - \frac{3}{2}\q_1^2, \label{eq:QMSMDcoche-a}\\
&X_2=-\frac{g^2_{11}}{\omega_2^2 }\q_1^2 = - \frac{1}{2}\q_1^2.\label{eq:QMSMDcoche-b}
\end{align}
\end{subequations}
One can first notice that for this specific example, the manifold produced with the SMD method does not depend on the parameters $(\omega_1,\omega_2)$. Consequently, the cut of this manifold in $(X_1,X_2)$ plane in Figs.~\ref{fig:manifoldcut2d}(a-c) for different values of $\rho$, is always the same.
The second comment is on the slow/fast approximation: even though the value given for $X_2$ tends to be the same under the slow/fast assumption $\omega_2 \gg \omega_1$, this is not the case for $X_1$. This is a major difference between the two methods, so that a persistent error on the manifold is done when using SMD, whereas MD tends to the solution provided by the NF and full system when $\rho$ increases. The last interesting comment is on the fact that the manifold produced by SMD shows a constant folding point. Indeed,  $X_1$ from Eq.~\eqref{eq:QMSMDcoche-a} cannot exceed the value of $1/6$ (achieved at $\q_1=1/3$) after which the quadratic term in Eq.~\eqref{eq:QMSMDcoche-a} is larger than the linear one.

This is a direct consequence of the different treatment of the self-quadratic coupling term $g^1_{11}$, already underlined at the end of Sect.~\ref{sec:drift}, leading to the fact that even under the slow/fast assumption, the QM built on SMD can lead to erroneous results. This point is further commented on the backbone curves comparison. First, Eqs.~\eqref{eq:zifullgamma} are written for this specific system, leading to the following predictions, as a function of the ratio $\rho=\omega_2/\omega_1$:
\begin{align}
&\Gamma_\text{MD}=
-\frac{16\rho^4
-27\rho^2+12}
{16(\rho^2-1)^2},
\\
&\Gamma_\text{SMD}=
-1,
\\
&\Gamma_\text{NF}=
-\frac{\rho^2-3}{\rho^2-4}.
\end{align}
In line with the constant manifold found with SMD, the method also predicts  a constant type of nonlinearity, independent of the variations of the eigenfrequencies  $(\omega_1,\omega_2)$. Assuming slow/fast partition, $\rho \rightarrow \infty$, then all three methods tends to predict the same $\Gamma$ coefficient dictating the hardening/softening behaviour. However, as underlined  at the end of Sect.~\ref{sec:drift}, the amplitude of the first harmonics for each method is different. Since in this case $g^1_{11}\neq 0$, a direct consequence of \eqref{eq:X_1H1exprSMD} is that the backbone curves for the SMD method will show a saturation effect, the amplitude $X_1$ being unable to overcome a maximum value. This phenomenon is clearly visible in Fig.~\ref{fig:BackbonesShell}, depicting  the backbone curves obtained for the three selected values of $\rho$. The constant value of $\Gamma_\text{SMD}$ has for direct consequence that the backbone predicted by the SMD quadratic manifold is almost unchanged with respect to variations of $\rho$. When the slow/fast assumption is fulfilled for $\rho=10$, Fig.~\ref{fig:BackbonesShell}(c),   the backbone predicted by SMD QM is in line with those predicted by the other methods at small amplitude level. However, at higher amplitude the SMD backbone moves away from the others and saturates to a limit value for all cases, since the amplitude is differently computed as shown in Eq.~\eqref{eq:X_1H1exprSMD}. On the other hand, the backbone predicted by the MD method tends to the correct values under the slow/fast approximation, while the normal form approach always produces a correct prediction. \red{More specifically, the prediction for the master $X_1$ component given by the normal form is the same if one considers a quadratic or cubic normal form expansion, see Eq.~\eqref{eq:transformNF}. On the other hand, the slave component $X_2$ is affected by the order and this is illustrated in Fig.~\ref{fig:BackbonesShell}(d-e-f), where one can observe that, as for the manifold approximation in phase space, the third-order terms bring about a better estimate.}


\section{Comparison on continuous structures}

\subsection{Presentation of the test cases}

\begin{figure*}[h!]
\centering
\begin{subfigure}[b]{.253\textwidth}
\includegraphics[scale=.33]{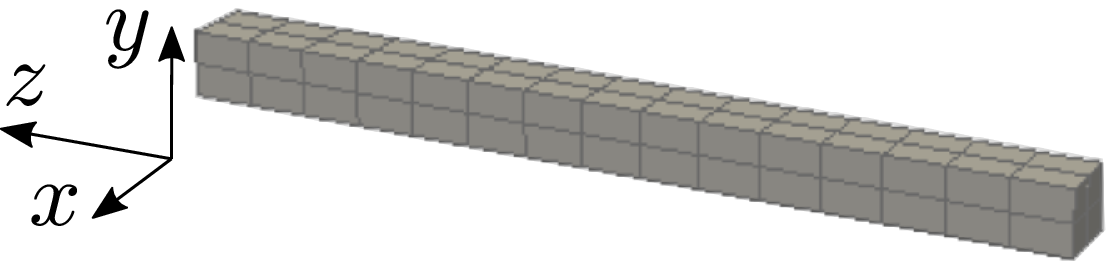}\\
\includegraphics[scale=.37]{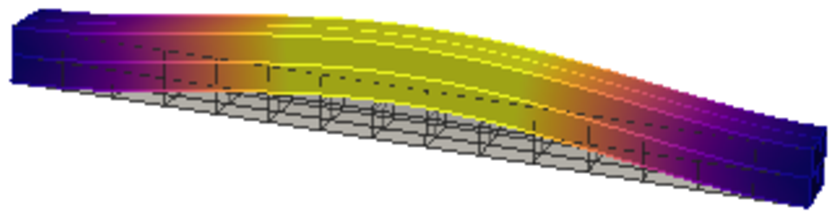}
\caption{Flat beam.}\label{fig:beams_flat}
\end{subfigure}
\begin{subfigure}[b]{.385\textwidth}
\includegraphics[scale=.33]{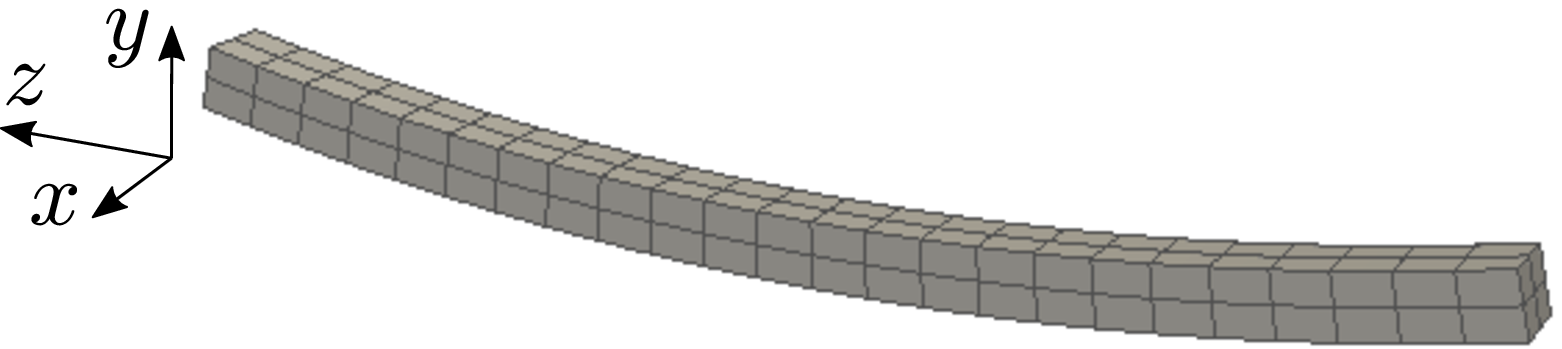}\\
\includegraphics[scale=.37]{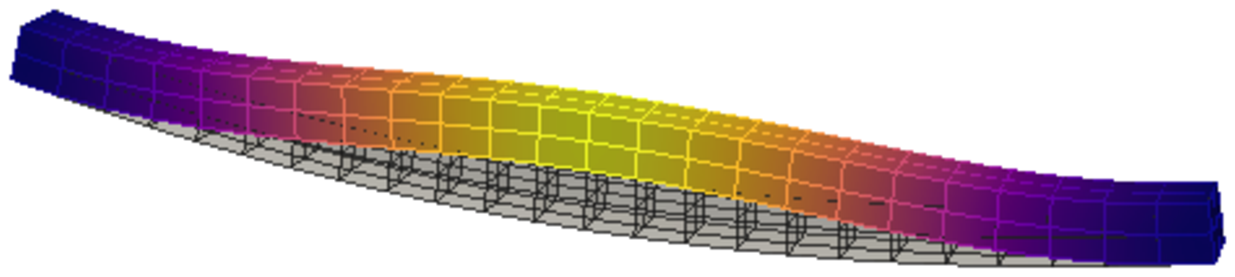}
\caption{Shallow arch.}\label{fig:beams_shallow}
\end{subfigure}
\begin{subfigure}[b]{.33\textwidth}
\includegraphics[scale=.33]{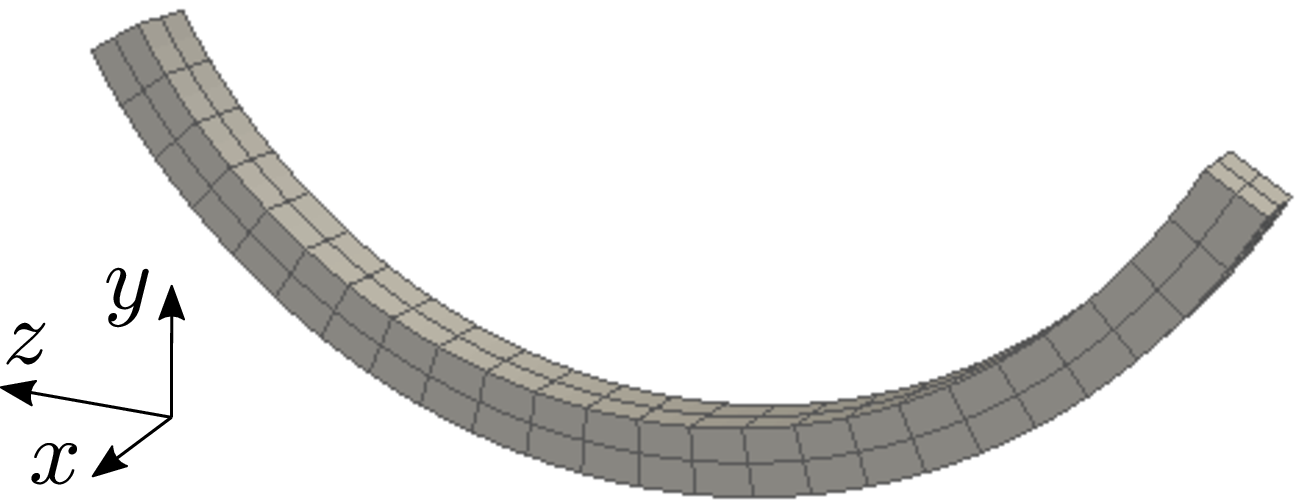}\\
\includegraphics[scale=.37]{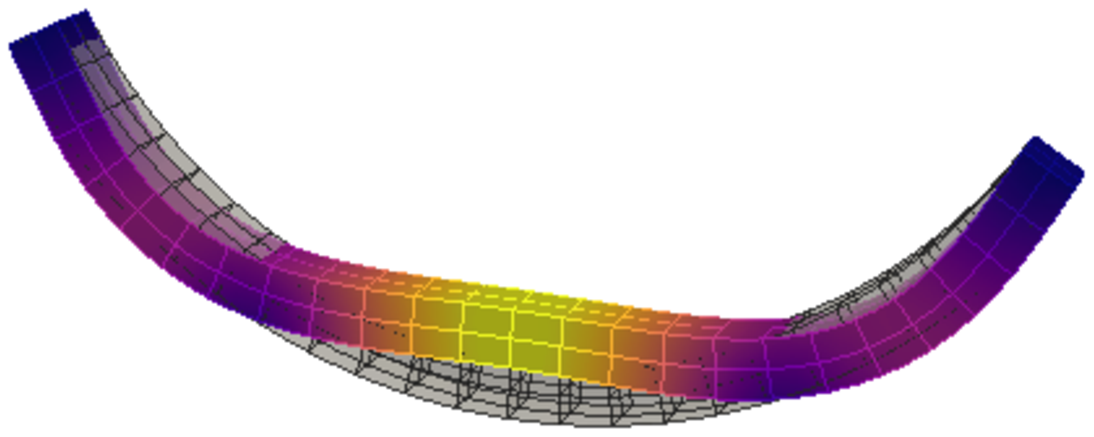}
\caption{Non-shallow arch.}\label{fig:beams_curved}
\end{subfigure}
\caption{Mesh and deformation along the mode under study for three different test cases: (a) a flat beam, (b) a shallow arch and (c) a non-shallow arch. For each test case the mode under study corresponds to the first flexural mode in the plane $y-z$ where the curvature is imposed if present.}
\label{fig:beams}
\end{figure*}

This section aims at drawing a comparison between the different methods when applied on a continuous structure discretised with three dimensional finite elements. In order to investigate how the results obtained in the previous section are confirmed in the general case, three beams are considered and shown  in Fig.~\ref{fig:beams}. They have been selected in order to fulfil different assumptions that have been highlighted on the two-dofs examples in order to achieve correct predictions from the ROMs. The first case, Fig.~\ref{fig:beams_flat}, is a slender flat symmetric beam. The two other examples, Fig.~\ref{fig:beams_shallow} and Fig.~\ref{fig:beams_curved} are arches, the first one being shallow while the third one is non-shallow. Adding curvature has two important effects. First, flexural and in-plane modes are no longer linearly uncoupled. Second, the curvature renders the restoring force asymmetric and an important quadratic nonlinearity appears between  the bending modes. This example  illustrates the fact that the slow/fast assumption is not enough to guarantee that the method based on static modal derivatives will converge. The curvature will be used in order to play on the slow/fast assumption as well as on the values of the quadratic coupling terms. 

In all three cases, the boundary conditions are clamped, the material parameters are selected as homogeneous linear elastic with Young modulus $E =124$ GPa, Poisson ratio $\nu = 0.3$ and density $\delta = 4400$ $\text{kg} \text{m}/\text{s}^2$. In each cases, an equal thickness $h$ and width $b$ is selected: $h=b=5\;cm$.  For the flat beam, the length is $L=0.7$ m. The arches have been built from a portion of a circle. For the shallow arch, the radius of curvature is set as $250$ cm, for an angular span of $2\pi/15$, resulting in a curvilinear length of $20\pi/3 \simeq 1.05$~m. The height of the static deflection at center is $5.5$ cm, {\em i.e.} almost equal to the thickness. For the non-shallow arch, the radius of curvature is set as $50$ cm, for an angular span of $2\pi/3$, resulting in the same curvilinear length of $20\pi/3 \simeq 1.05$~m, but with a static deflection of $25$ cm, {\em i.e.} 5 times the thickness. 
All beams are discretised with three-dimensional hexahedral 20 nodes finite elements. The flat beam uses 60 elements (4 in the section and 15 in the length), resulting in a total number of 1287 dofs. The two arches have 96 solid elements (4  in the section and 24 in the length) and 2097 dofs. A relative coarse mesh has been selected in order to have a limited number of degrees of freedom so that all the methods can be handled easily. Indeed, the key point here is not to look for converged and refined results on a large frequency range, but to compare the different reduction methods on the same test examples. Moreover, as already shown in~\cite{Vizza3d}, using 3D elements leads to couplings with very high-frequency thickness modes, so that truncations and convergence are difficult to observe in general.

In the three cases, the nonlinear behaviour of the first flexural mode in the curvature plane is investigated. The mode shape is shown in Fig.~\ref{fig:beams}. In the case of the flat beam, it corresponds to the first mode and its eigenfrequency is $545.60$ Hz. As already underlined in Sect.~\ref{sec:2dofsFLAT}, the most important coupling arises with the fourth in-plane mode, whose eigenfrequency is $15.19$ kHz, so that the ratio $\rho$ between the most important slave mode and the master mode is in this case equal to $27.83$. Consequently the slow/fast assumption and our criterion $\rho \geq 4$ is perfectly fulfilled. This example can be seen as an extension of the first two-dofs example, with the distinctive feature that now many more modes are coupled to the first bending, all of them being of higher frequencies than the fourth axial. Also, the nonlinear coupling terms have in this case a simplified form, following the general discussion given in Sect.~\ref{sec:2dofsFLAT}. In the case of the arches, for the shallow arch the first flexural mode corresponds to the second mode of the structure and its eigenfrequency is equal to $372.28$ Hz and, for the non-shallow arch, the first flexural mode corresponds to the fourth mode of the structure and its eigenfrequency is equal to $1003.99$ Hz. Contrary to the case of flat symmetric structures, the curvature renders the restoring force asymmetric and an important quadratic nonlinearity appears between the bending modes. 
Investigating the important couplings between the linear modes of the curved beams shows that the first bending mode is strongly coupled with the third one.  While the ratio between the first and third bending modes is 5.4 in the case of the flat beam, it decreases when the curvature increases. Consequently, for the case of the shallow arch, this ratio is equal to $3.44$ (eigenfrequency of third bending equal to $1283.33$ Hz), and $1.66$ for the non-shallow arch (eigenfrequency of third bending equal to $1665.11$ Hz). These two examples have thus been built as an extension of the second two-dofs example. For the shallow arch, the slow/fast assumption is almost fulfilled (3.44 is slightly smaller than the criterion we proposed with a limit value at 4), but now important quadratic couplings are present and in particular the self-quadratic term $g^p_{pp}$.   Finally, the case of the non-shallow arch allows testing a case where the slow/fast assumption does not hold, and important self-quadratic terms are present. 

\subsection{Amplitude-frequency relationships}

The methods are compared on their ability to predict the backbone curves. A reference solution is computed thanks to a numerical continuation on all the degrees of freedom of the structure, using a code with parallel implementation of harmonic balance method and pseudo arc-length continuation algorithm \cite{Blahos2020}.  In this computation, a small amount of mass proportional damping is added under the form $\zeta \omega_p \vec{M}$ so that a frequency-response function (FRF) is computed, in the vicinity of the eigenfrequency of the master mode (first flexural). The values of $\zeta$ are $0.18\%$, $0.27\%$, and $0.1\%$ for the flat beam, shallow, and non-shallow arches, respectively. The forcing is located in the central node of each mesh in the $y$-direction in order to excite the first flexural mode. The force amplitude is chosen in order to have a displacement amplitude at resonance comparable to the thickness so that its values are $5$ kN, $1.5$ kN, and $2.5$ kN for the flat beam, shallow, and non-shallow arches, respectively. It must be noticed that in the case of curved structures the value of amplitude of vibration equal to the thickness has not been achieved and the reported FRFs excite a maximum amplitude of approximately half of it. In fact, due to the long computational time that the full model FRF requires, approximately 1 day for each FRF, and due to its high chances to undergo internal resonance with higher modes, these values have been selected in order to stay in the limit of one-mode approximation without exciting more complex dynamics. However, with this level of amplitudes, the nonlinearity is sufficiently important so that its effect is clearly visible on the backbone curves.

The ROMS are built using QM or NNM approach, and their backbone curves are computed in the same manner than in the previous section, assuming a single master mode. \red{For the normal form approach, the third order coefficients have not been included in the computation. Indeed, the third-order tensors require the computation of huge number of coupling coefficients from the modal basis expressions, which would need for an important number of pre-computation steps. This choice has also been guided by the fact that comparing the two methods at the same order of accuracy is more meaningful.} The FRF of the ROMS have not been computed since taking into account the damping of the slave modes is important to achieve good results. If the normal form theory has been developed for that purpose, see {\em e.g.}~\cite{TOUZE:JSV:2006} where the effect of a small amount of damping of the slave modes on the FRF of the master mode is reported, the inclusion of the damping for the modal derivatives has not been derived theoretically yet. Hence it appears that a better comparison is given on the backbone curves only, and the FRF of the full model with a small amount of damping is used to underline if the nonlinearity is correctly addressed by the methods.
\begin{figure*}[h!]
\centering
\begin{subfigure}[b]{.323\textwidth}
\includegraphics[width=.9\textwidth]{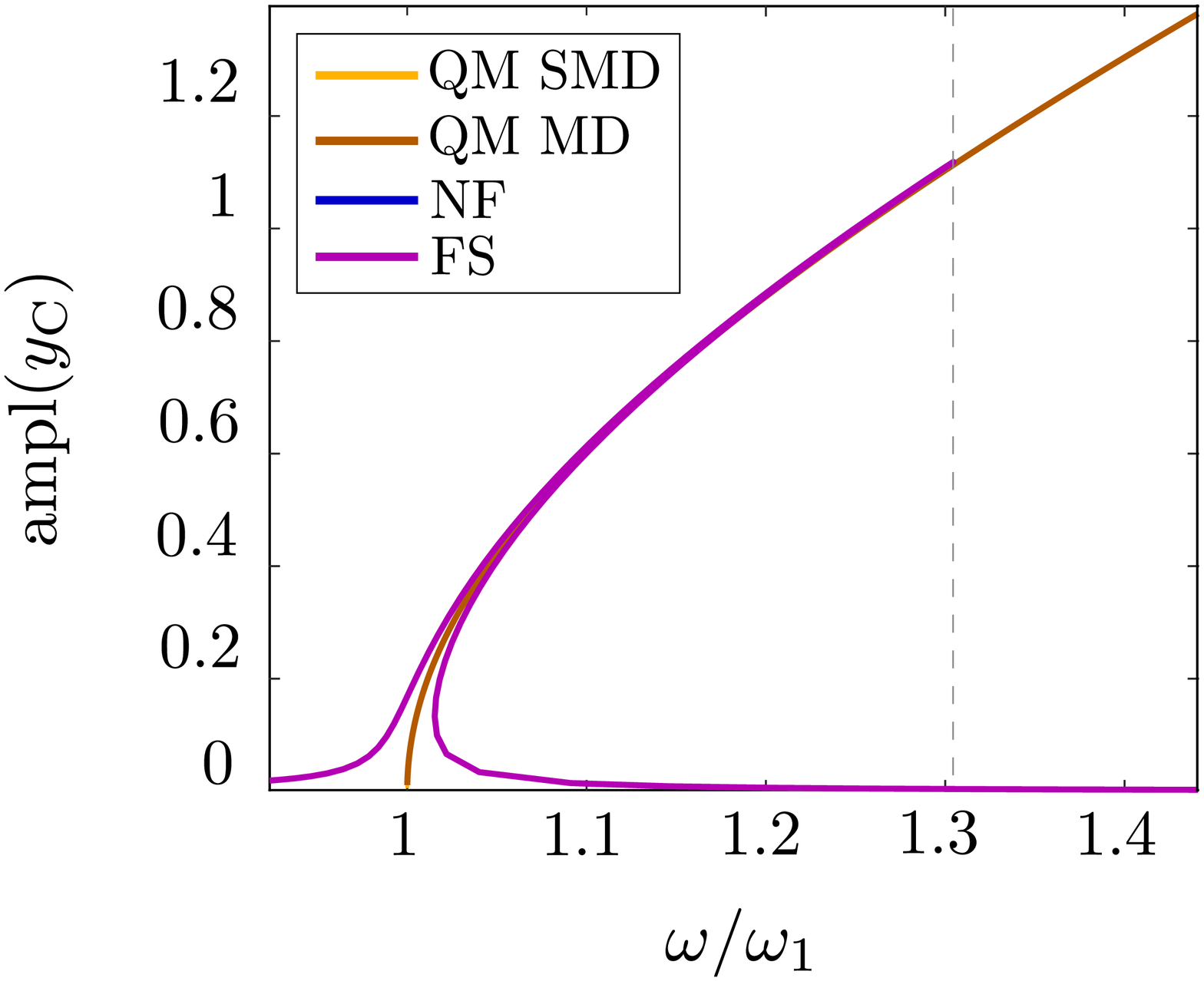}\caption{Flat beam.}\label{fig:frf_flat}
\end{subfigure}
\begin{subfigure}[b]{.323\textwidth}
\includegraphics[width=.9\textwidth]{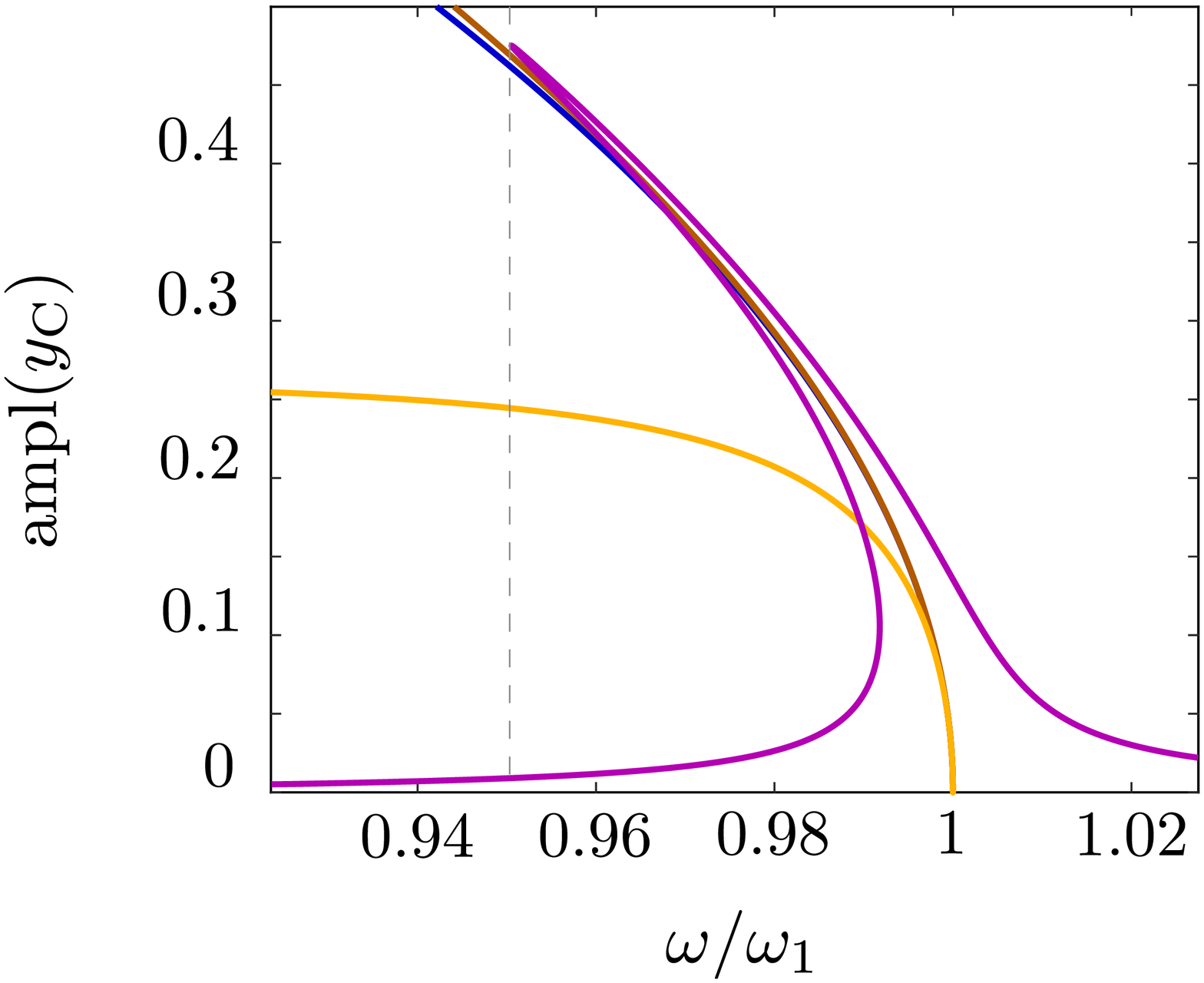}\caption{Shallow arch.}\label{fig:frf_shallow}
\end{subfigure}
\begin{subfigure}[b]{.323\textwidth}
\includegraphics[width=.9\textwidth]{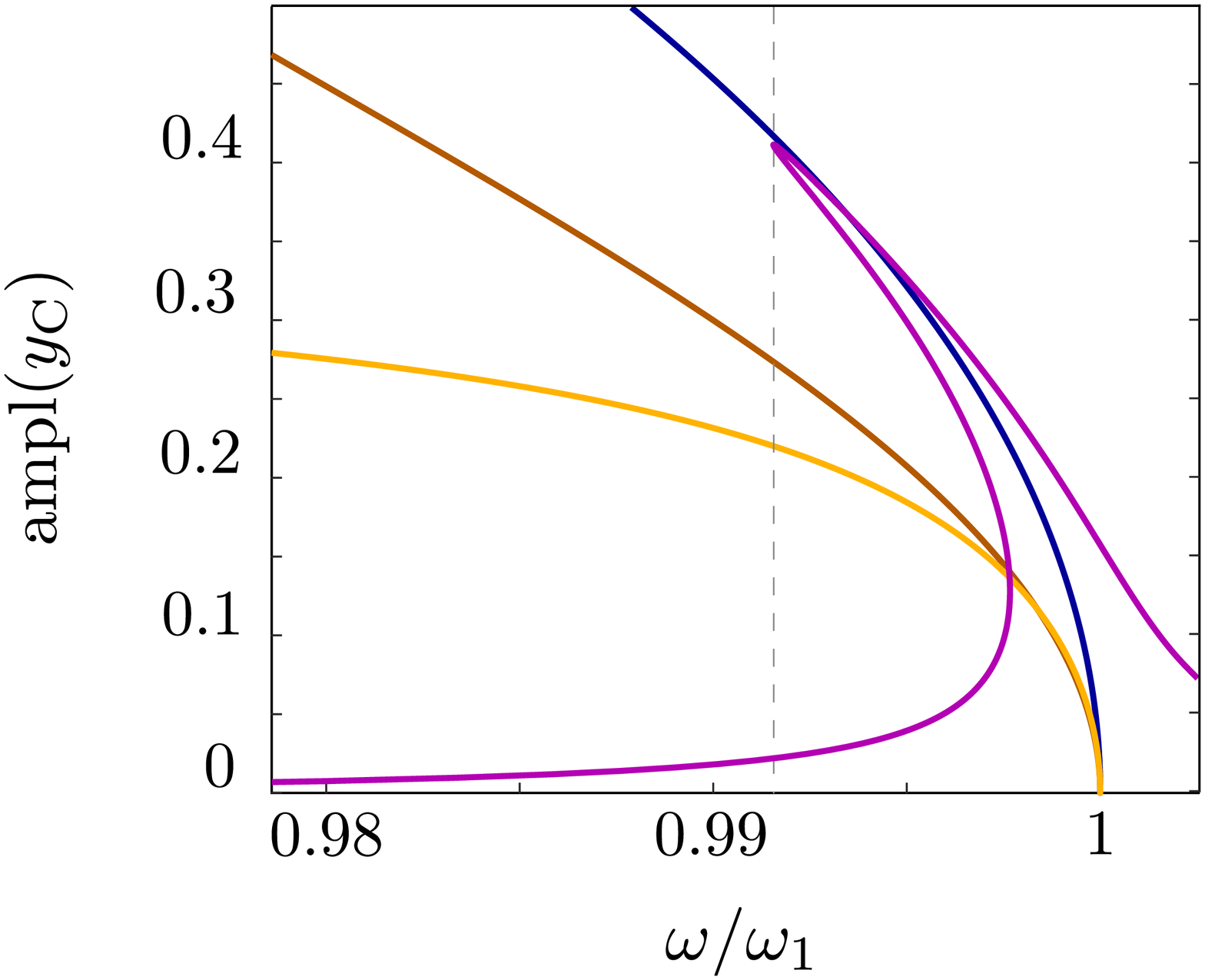}\caption{Non-shallow arch.}\label{fig:frf_curved}
\end{subfigure}
\caption{Comparison of backbone curves obtained from QM with MDs (dark orange), SMDs (yellow) and normal form approach (blue), for the three tested structures: (a) flat beam, (b) shallow arch, (c) non-shallow arch. Nondimensional amplitude of flexural displacement (along $y$, nondimensionalised with respect to the thickness) of the central node of each beam as a function of $\omega/\omega_1$ where $\omega_1$ refers to the  eigenfrequency of the first flexural mode studied shown in Fig.~\ref{fig:beams}. The backbone curves are contrasted to the FRF obtained on the full system (FS, violet) with numerical continuation and a small amount of damping, see text. The vertical gray dashed lines represent the frequency at which the saddle-node bifurcation of the forced response occurs.}
\label{fig:bbfrfbeams}
\end{figure*}
Fig.~\ref{fig:bbfrfbeams} shows the numerical results obtained for the three cases. The case of the flat beam is the one having the most assumptions fulfilled (slow/fast separation and no self-quadratic terms). Consequently, the three methods match very well and are all able to retrieve correctly the nonlinearity of the full model with a very good accuracy. In the case of the shallow arch, the slow/fast assumption is almost fulfilled (since being a little bit below the proposed criterion $\rho \geq 4$), and important self-quadratic coupling appears due to the curvature. The main consequence is that the QM built from SMD is not able anymore to predict the correct type of nonlinearity. As already found for the second two-dofs example, it overpredicts the softening behaviour and make appear again the saturation phenomenon in the amplitude of the backbone. On the other hand, both QM based on MD and normal form methods give a correct prediction. For the non-shallow arch, the slow/fast assumption does not hold anymore. The consequence is that the MD method does not predict the correct nonlinearity. This example again illustrates clearly that: (i) as soon as important self-quadratic terms appear (case of arches and shells), then the SMD method is not reliable anymore, whatever the slow/fast assumption is fulfilled or not, (ii) the MD can still give correct result but only if the criterion $\rho \geq 4$ for the slow/fast assumption is fulfilled. As soon as $\rho$ gets under this value, then the solution starts departing from the full-order model, and becomes unreliable when $\rho \leq 2$.

\subsection{Nonlinear modeshapes}

The different approximations made by the three methods are finally contrasted on the mode shape dependence on amplitude,  illustrating the Equations given in Sect.~\ref{sec:drift}. Recalling Eqs.~\eqref{eq:reconsuMD}~--~\eqref{eq:reconsuNF}, it is possible to see that, for each method, the contributions to the nonlinear modeshape can be divided into (i) a deformation along the master $p$ mode and (ii) a deformation that contains all the coupled modes but the $p$-th. In order to make the figures more illustrative, and since the amplitude of the deformation along $p$-th mode generally gives the dominant contribution, it is decided to compare the outcomes of the methods only on the (ii) part of the solution. Also,  since the normal form approach constructs the solution both with displacements and velocities, to draw a better comparison the focus will be on the time step where the reduced variable $\q_p(t)$ reaches its maximum and minimum values (i.e. a turning point such that $\dot{\q}_p(t)=0$).

Under these assumptions, let us define as  $\vec{u}^{\perp}$ the component of the nonlinear mode shape $\vec{u}$ that is orthogonal to $\vec{\phi}_p$. From Eqs.~\eqref{eq:reconsuMD}~--~\eqref{eq:reconsuNF}, it reads, for the three different methods:
\begin{subequations}\label{eq:uperprom}
\begin{align}
\vec{u}^{\perp}_\text{MD}(t^*)=&
-\q_p^2(t^*) 
\underset{s\neq p}{\sum_{s=1}^N }
\frac{g^s_{pp}}{\omega_s^2-\omega_p^2}
\vec{\phi}_s,
\\
\vec{u}^{\perp}_\text{SMD}(t^*)=&
-\q_p^2(t^*) 
\underset{s\neq p}{\sum_{s=1}^N}
\frac{g^s_{pp}}{\omega_s^2}
\vec{\phi}_s,
\\
\vec{u}^{\perp}_\text{NF}(t^*)=&
-\q_p^2(t^*) 
\underset{s\neq p}{\sum_{s=1}^N}
\frac{g^s_{pp}}{\omega_s^2}
\left(
\frac{\omega_s^2-2\omega_p^2}{\omega_s^2-4\omega_p^2}
\right) \vec{\phi}_s,
\end{align}
\end{subequations}
where $t^*$ is the time instant where $\q_p$ is either maximum or minimum. 

In order to compare to the full-order solution, the deformation must be first filtered out from its component along the $p$-th mode. One can thus define $\vec{u}^{\perp}_\text{FS}$ for the full system as:
\begin{equation}
\vec{u}^{\perp}_\text{FS}(t)=
\vec{u}_\text{FS}(t)-\dfrac{\vec{\phi}^\text{T}_p \vec{u}_\text{FS}(t)}
{\vec{\phi}^\text{T}_p \vec{\phi}_p}\vec{\phi}_p.
\end{equation}
Finally, given the quadratic nature of the deformation computed from the reduction methods based on second-order expansions (and clearly underlined by the dependence in $R_p^2$ in Eqs.~\eqref{eq:uperprom}), the third-order component should be also filtered out from $\vec{u}^{\perp}_\text{FS}$ for a closer comparison. In order to cancel the odd harmonics of the full-order solution, we thus define $\vec{u}^{\perp,sym}_\text{FS}$ as the symmetric part of $\vec{u}^{\perp}_\text{FS}$ with respect to amplitude:
\begin{equation}\label{eq:uperpfullsym}
\vec{u}^{\perp,sym}_\text{FS}=\dfrac{1}{2}(
\vec{u}^{\perp}_\text{FS}(t^{max})+
\vec{u}^{\perp}_\text{FS}(t^{min})).
\end{equation}
This value will be used as reference and compared to the prediction of the ROMS given by Eq.~\eqref{eq:uperprom}.

\begin{figure*}[h!]
\centering
\begin{subfigure}[b]{.24\textwidth}
\includegraphics[width=\textwidth]{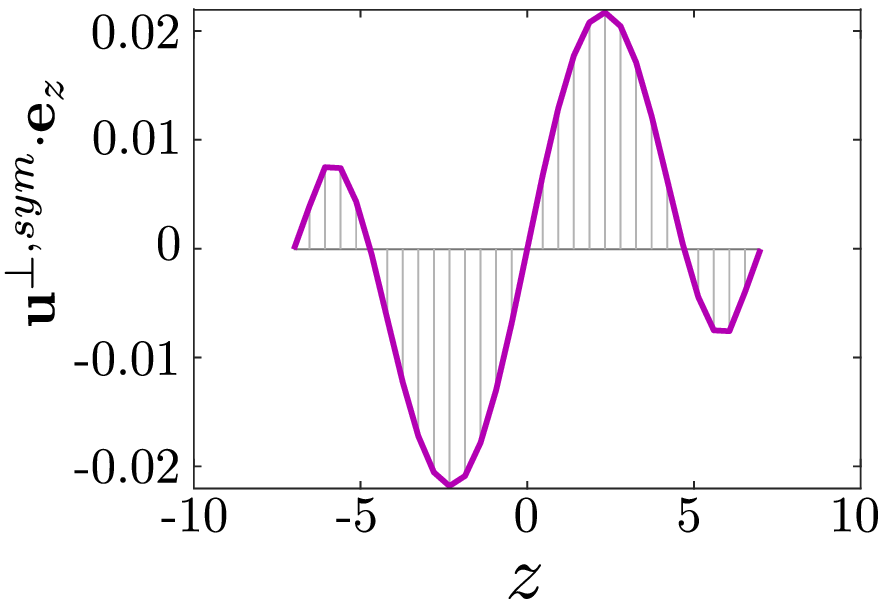}
\caption{Full Model.}
\end{subfigure}
\begin{subfigure}[b]{.24\textwidth}
\includegraphics[width=\textwidth]{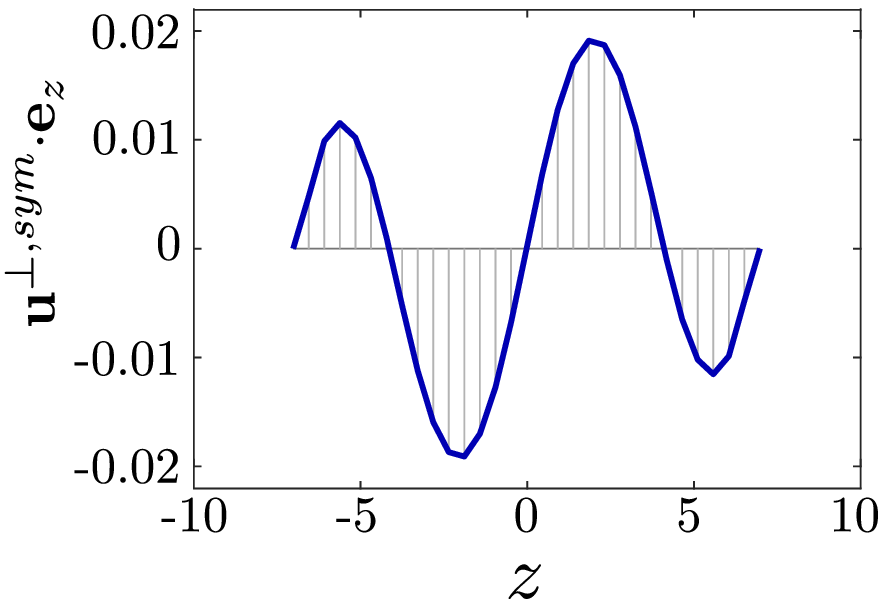}
\caption{Normal Form.}
\end{subfigure}
\begin{subfigure}[b]{.24\textwidth}
\includegraphics[width=\textwidth]{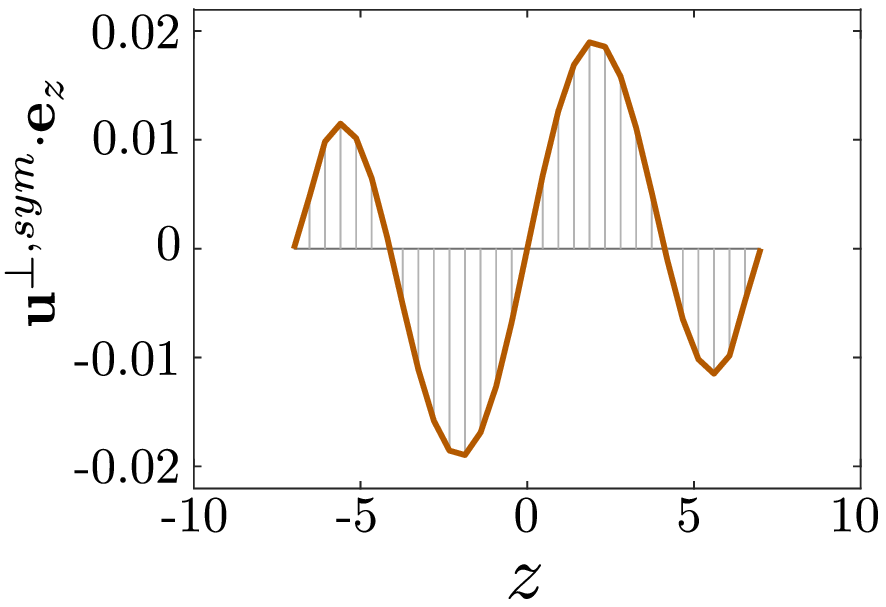}
\caption{Modal Derivatives.}
\end{subfigure}
\begin{subfigure}[b]{.24\textwidth}
\includegraphics[width=\textwidth]{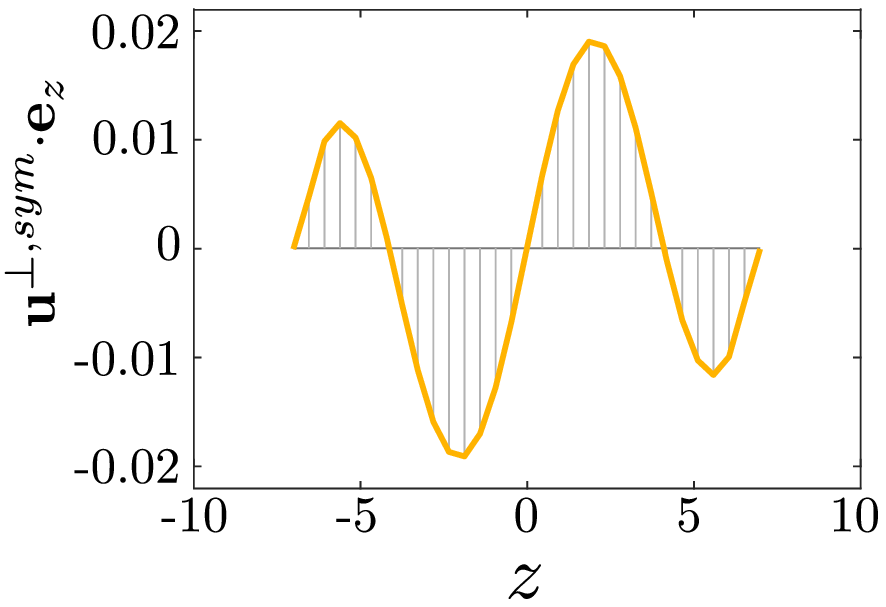}
\caption{Static MD.}
\end{subfigure}
\caption{Comparisons of the additional terms perturbing the linear mode shape solutions (deformation orthogonal to $\vec{\phi}_p$) for the case of the flat beam, computed at the saddle-node bifurcation point marked in Fig.~\ref{fig:frf_flat}, fixing the frequency at which they have been computed.  (a) full model solution, representation of the axial component of displacement $\vec{u}^{\perp,sym}_\text{FS}.\vec{e}_z$ of the centre line nodes, (b) normal form: $\vec{u}^{\perp}_\text{NF}(t^*).\vec{e}_z$, (c) Modal derivative : $\vec{u}^{\perp}_\text{MD}(t^*).\vec{e}_z$, (d) Static modal derivative $\vec{u}^{\perp}_\text{SMD}(t^*).\vec{e}_z$.}
\label{fig:3D_res_def_quad_flat}
\end{figure*}
Fig.~\ref{fig:3D_res_def_quad_flat} shows the comparison between the $\vec{u}^{\perp}$ defined by Eq.~\eqref{eq:uperpfullsym} for the full-order system, and those produced by the reduced-order models, Eqs.~\eqref{eq:uperprom}, for the case of the flat beam. Importantly enough, since the nonlinear couplings are with in-plane modes, the contributions of the $\vec{u}^{\perp}$ along the axial $z$ direction is shown in Fig.~\ref{fig:beams}, since  the most important contributions are along this direction. As it could be awaited from the previous analyses, Fig.~\ref{fig:3D_res_def_quad_flat} clearly shows that the three ROMs are all able to recover the correct spatial dependence of the contributions of coupled modes to the fundamental flexural NNM. Also, this contribution is mostly conveyed by the fourth in-plane mode, being the most importantly coupled to the fundamental flexural mode. Note that the amplitude used to construct this figure is the one corresponding to the upper saddle-node bifurcation point in the FRF of the full-order system, as shown by the gray vertical line in Fig.~\ref{fig:frf_flat}. At that point, the backbones and the FRF meet so that it can be used safely for a correct comparison. It also corresponds to an amplitude of one time the thickness for the mode shape.

\begin{figure*}[h!]
\centering
\begin{subfigure}[b]{.24\textwidth}
\includegraphics[width=\textwidth]{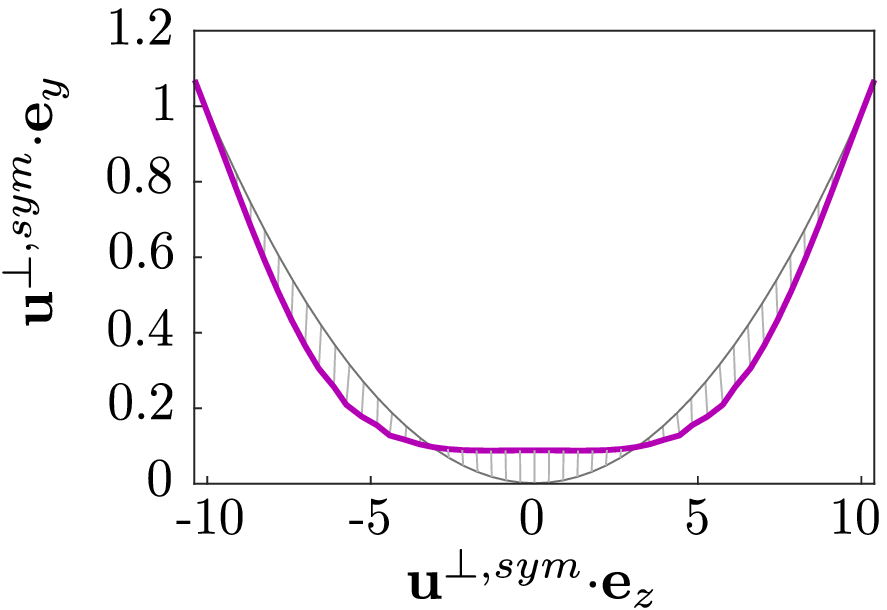}
\caption{Full Model.}
\end{subfigure}
\begin{subfigure}[b]{.24\textwidth}
\includegraphics[width=\textwidth]{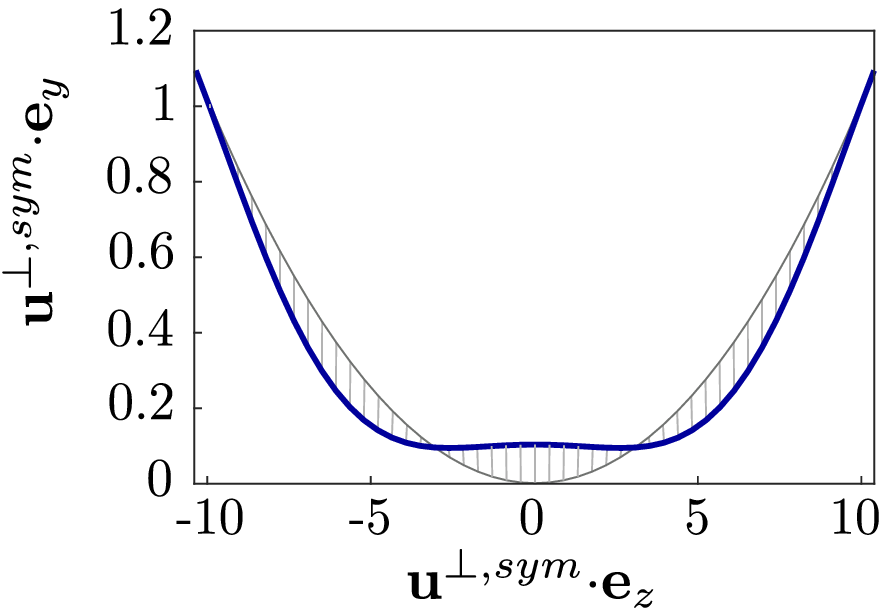}
\caption{Normal Form.}
\end{subfigure}
\begin{subfigure}[b]{.24\textwidth}
\includegraphics[width=\textwidth]{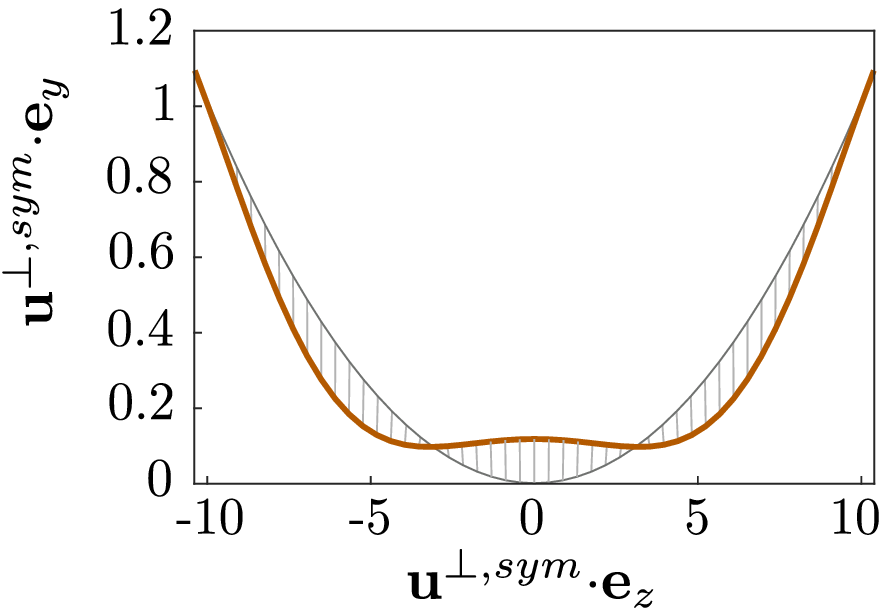}
\caption{Modal Derivatives.}
\end{subfigure}
\begin{subfigure}[b]{.24\textwidth}
\includegraphics[width=\textwidth]{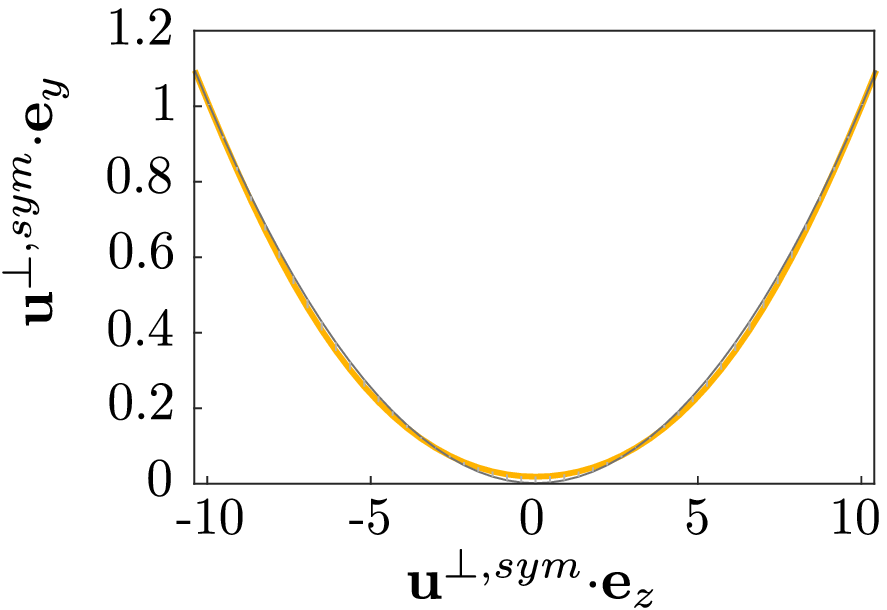}
\caption{Static MD.}
\end{subfigure}
\caption{Comparisons of the additional terms perturbing the linear mode shape solutions (deformation orthogonal to $\vec{\phi}_p$) for the case of the shallow arch, computed at the saddle-node bifurcation point marked in Fig.~\ref{fig:frf_shallow}, fixing the frequency at which they have been computed.  (a) full model solution, representation of the displacement $\vec{u}^{\perp,sym}_\text{FS}$ of the centre line in the $y-z$ plane; vertical axis $\vec{u}^{\perp,sym}_\text{FS}.\vec{e}_y$ (transverse direction) and horizontal axis $\vec{u}^{\perp,sym}_\text{FS}.\vec{e}_z$ (axial direction), (b) normal form: $\vec{u}^{\perp}_\text{NF}(t^*)$, (c) Modal derivative : $\vec{u}^{\perp}_\text{MD}(t^*)$, (d) Static modal derivative $\vec{u}^{\perp}_\text{SMD}(t^*)$. Gray lines:  position of the centre line of the beam at rest. Solution amplified of factor 15.}
\label{fig:3D_res_def_quad_shallow}
\end{figure*}
In the case of the shallow arch, some differences are appearing due to the self-quadratic coupling term, creating a deficiency in the prediction given by the SMD. This is underlined in the nonlinear mode shape dependence in Fig.~\ref{fig:3D_res_def_quad_shallow}, where in this case, since the most important coupling is between bending modes, the contributions of the different $\vec{u}^{\perp}$ are represented along the transverse $y$ direction. The amplitude used for the figure is illustrated in Fig.~\ref{fig:frf_shallow} with a gray line, and still corresponds to the upper saddle-node bifurcation point in the FRF of the full-order system. One can observe in Fig.~\ref{fig:3D_res_def_quad_shallow} that, in the line of the results found on the nonlinear amplitude-frequency relationships, normal form and MD methods are able to retrieve the correct spatial dependence for the contribution of the slave modes. On the other hand, the treatment of the self-quadratic term by the SMD approach prevents the correct prediction of this spatial dependence. 
\begin{figure*}[h!]
\centering
\begin{subfigure}[b]{.24\textwidth}
\includegraphics[width=\textwidth]{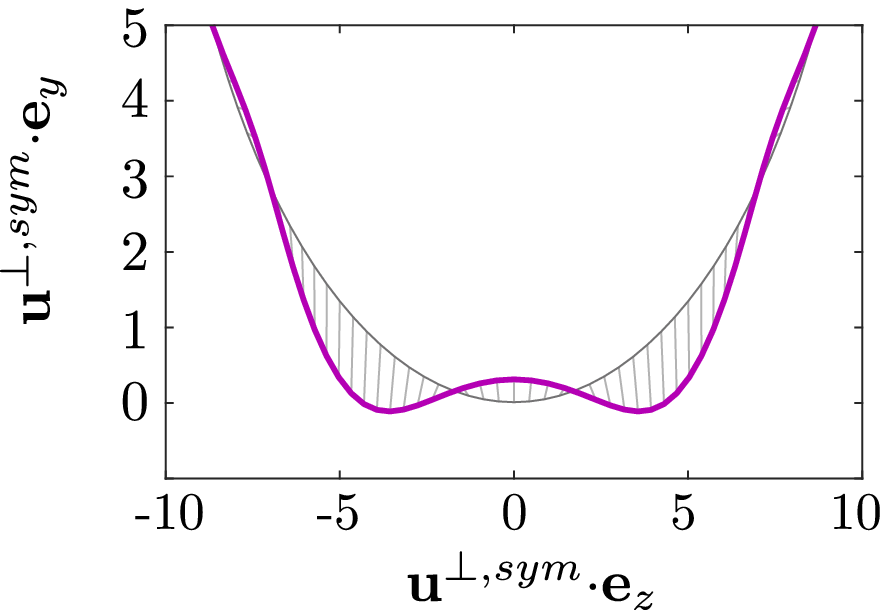}
\caption{Full Model.}
\end{subfigure}
\begin{subfigure}[b]{.24\textwidth}
\includegraphics[width=\textwidth]{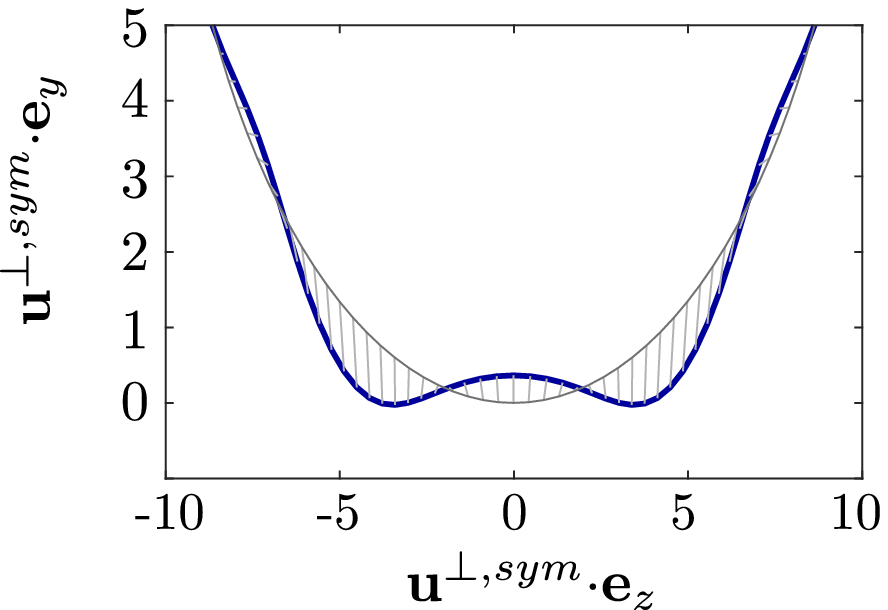}
\caption{Normal Form.}
\end{subfigure}
\begin{subfigure}[b]{.24\textwidth}
\includegraphics[width=\textwidth]{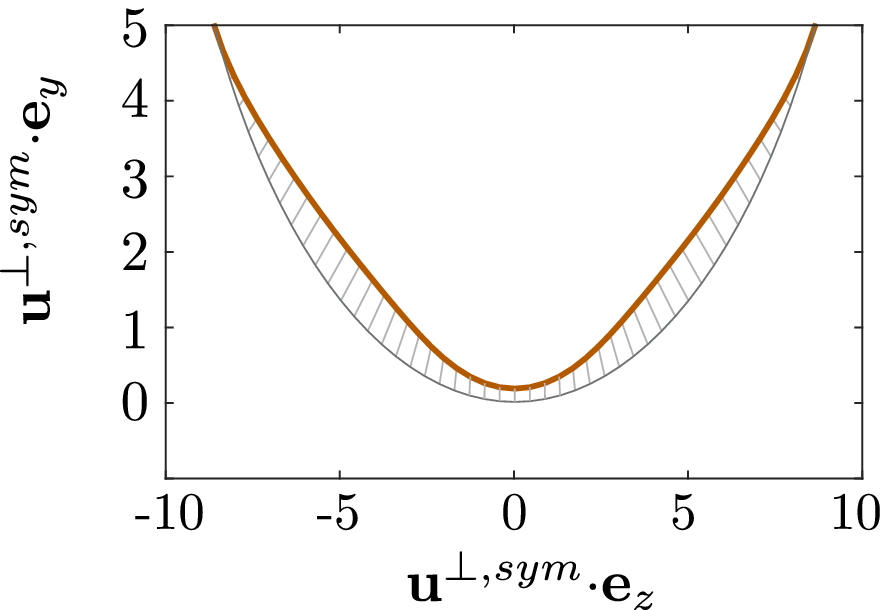}
\caption{Modal Derivatives.}
\end{subfigure}
\begin{subfigure}[b]{.24\textwidth}
\includegraphics[width=\textwidth]{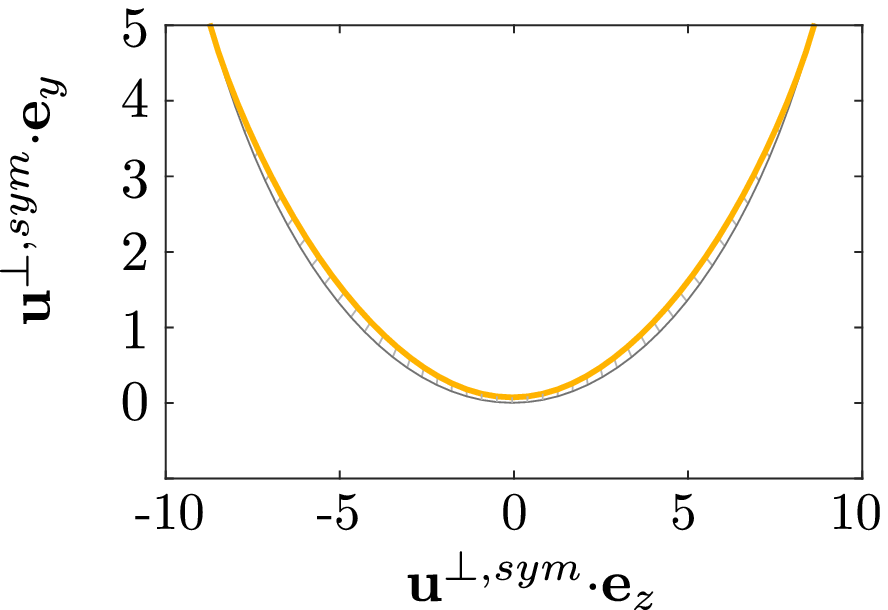}
\caption{Static MD.}
\end{subfigure}
\caption{Comparisons of the additional terms perturbing the linear mode shape solutions (deformation orthogonal to $\vec{\phi}_p$) for the case of the non-shallow arch, computed at the saddle-node bifurcation point marked in Fig.~\ref{fig:frf_curved}, fixing the frequency at which they have been computed.  (a) full model solution, representation of the displacement $\vec{u}^{\perp,sym}_\text{FS}$ of the centre line in the $y-z$ plane; vertical axis $\vec{u}^{\perp,sym}_\text{FS}.\vec{e}_y$ (transverse direction) and horizontal axis $\vec{u}^{\perp,sym}_\text{FS}.\vec{e}_z$ (axial direction), (b) normal form: $\vec{u}^{\perp}_\text{NF}(t^*)$, (c) Modal derivative : $\vec{u}^{\perp}_\text{MD}(t^*)$, (d) Static modal derivative $\vec{u}^{\perp}_\text{SMD}(t^*)$. Gray lines:  position of the centre line of the beam at rest. Solution amplified of factor 50.}
\label{fig:3D_res_def_quad_curved}
\end{figure*}

The case of the non-shallow arch is shown in Fig.~\ref{fig:3D_res_def_quad_curved}, for an amplitude of motion marked by the gray line in Fig.~\ref{fig:frf_curved}. Following the observation on the frequency, one can notice that only the normal form approach is able to retrieve the correct spatial dependence. On the other hand, SMD method fails because of the incorrect treatment of the self-quadratic term, while QM constructed from MD does not produce the correct result since the slow/fast assumption is not fulfilled anymore.

\section{Conclusion}
In this contribution, a detailed comparison of different methods proposed in the recent years in order to define nonlinear mappings with the aim of providing accurate reduced-order models for geometrically nonlinear structures, has been made. The quadratic manifold proposed from the definitions of modal derivatives has thus been contrasted to the normal form theory, related to the definition of nonlinear normal modes as invariant manifolds in phase space. While the quadratic manifold only contains the displacements as unknowns, the normal form approach takes into account displacements and velocities, thus giving a more complete link to the geometry in phase space. Secondly, the quadratic manifold is defined up to the second-order while current expressions of normal form are up to order three and can be continued to higher orders easily. Thirdly, normal form theory relies of firm mathematical theorems, ensuring a clean conceptual framework, while modal derivatives appear as an ad-hoc, yet efficient, method used in the vibration community.

The main outcomes of this article are the following. First, the theoretical derivations of the quadratic manifold using either MD or SMD, has been fully made explicit. These calculations have highlighted the fact that both methods do not handle the quadratic terms in the same manner, and especially the self-quadratic coupling terms arising between the master coordinates. This difference has been found to have important consequences on the global predictions of the methods. Secondly,  detailed comparisons between the three methods have been fully analysed on the mathematical expressions: nonlinear change of coordinates, reduced-order dynamics, and main predictive outcomes of the methods such as type of nonlinearity, drift and mode shape dependence on amplitude. To illustrate the results, two two-dofs systems have been used as starting example, and the results found from these have been extended to a continuous structure: a clamped-clamped beam with varying curvature.

A main result of our investigations is that the results predicted by the QM approach with MDs converge to those provided by the normal form approach, only in the case where a slow/fast assumption between master and slave coordinates, holds. \red{This result is fully in the line of general theorems provided in~\cite{HallerSF,VERASZTO}, and thus further illustrates the general findings given in these papers where a more general framework including damping is given, together with an exact result that do not rely on asymptotic expansion.} A first quantification of the limit value for the slow/fast assumption to hold has been provided, based on the predicted values for  the type of nonlinearity, showing that a small gap is needed: $\omega_s > 4 \omega_p$, thus justifying a posteriori the good results found by previously published papers using this method. However, the different treatment of the quadratic nonlinearity (and more specifically the self-quadratic coupling term) between MD and SMD, leads to the fact that even with a slow/fast assumption, the QM built from SMD can lead to erroneous predictions, as soon as an important self-quadratic coupling term is present. This result has important implications when one wants to build ROMs for slender curved structures such as arches and shells. This specific feature has been clearly highlighted on the two-dofs system, and found in the more general case of a non shallow arch. On the other hand, the robustness of the normal form approach has been underlined in each case.

These results argue for the use of the tools from dynamical system theory to derive safe and robust ROMS: invariant manifold, normal form theory and spectral submanifold. A limitation could be the use of these methods in the context of FE models where the need of computing, possibly in a non-intrusive manner, the nonlinear coefficients might be a difficult task, see {\em e.g.} all the literature related to the STEP method (Stiffness Evaluation Procedure, see {\em{e.g}}~\cite{muravyov,mignolet13,givois2019}). However recent developments show that the coefficients can be directly computed, for the case of spectral submanifold~\cite{VERASZTO}, or for the case of normal form in a non-intrusive manner~\cite{vizzaENOC}, so that this limitation does not hold anymore.

\appendix
\section{Definition of dot products for tensors}
\label{app:products}

The products between quadratic nonlinear tensors $\vec{G}$ and $\vec{g}$ and the displacement vectors $\vec{u}$ and $\vec{X}$, in physical and modal basis respectively, are defined as follows:
\begin{subequations}\begin{align}
\vec{G}\vec{u}\vec{u}
=&
\sum_{i=1}^N \sum_{j=1}^N
\vec{G}_{ij}{u}_{j}{u}_{i},\label{eq:product_Guu}\\
\vec{g}\vec{X}\vec{X}
=&
\sum_{i=1}^N \sum_{j=1}^N
\vec{g}_{ij}{X}_{j}{X}_{i},\label{eq:product_gXX}
\end{align}\end{subequations}
where $\vec{G}_{ij}$ and $\vec{g}_{ij}$ are the N-dimensional vector of coefficients $G^p_{ij}$ and $g^p_{ij}$, for $p=1..N$. Note that, in the course of this paper, the summation are kept complete so that all the terms in Eqs.~\eqref{eq:product_Guu}-\eqref{eq:product_gXX} are present. This is a choice of representation, the other choice (often realized in literature) consisting in symmetrizing the tensor, given the fact that products ${u}_{j}{u}_{i}$ can commute. In this case, a quadratic tensor is made symmetric such that for example $\vec{G}_{ij}=0$ when $j<i$, so that the summations can be written for $j\geq i$ only. Here we consider full tensors of coefficients without using their potential symmetry.

Similarly, the cubic nonlinear tensors $\vec{H}$ and $\vec{h}$, with current terms $H^p_{ijk}$ and $h^p_{ijk}$, contracts to an $N$-dimensional vector when multiplied with three displacements vectors. The cubic terms thus  explicitly writes, in indicial components, with $\vec{u}$ and $\vec{X}$ the two associated displacement vectors:
\begin{subequations}\begin{align}
\vec{H}\vec{u}\vec{u}\vec{u}
=&
\sum_{i=1}^N \sum_{j=1}^N \sum_{k=1}^N
\vec{H}_{ijk}{u}_{k}{u}_{j}{u}_{i},\label{eq:product_Huuu}\\
\vec{h}\vec{X}\vec{X}\vec{X}
=&
\sum_{i=1}^N \sum_{j=1}^N \sum_{k=1}^N
\vec{h}_{ijk}{X}_{k}{X}_{j}{X}_{i},\label{eq:product_hXXX}
\end{align}\end{subequations}
where again $\vec{H}_{ijk}$ is the $N$-dimensional vector with entries $H^p_{ijk}$, for $p=1...N$. The innermost product defined above coincides with a matrix product performed on the last index of the tensors; a more extensive definition of this notation is provided:
\begin{subequations}
\begin{align}
(\vec{G}
\vec{u})_i
=&
\sum_{j=1}^N
\vec{G}_{ij}
{u}_{j}\label{eq:product_Gu}\\
(\vec{H}
\vec{u})_{ij}
=&
\sum_{k=1}^N
\vec{H}_{ijk}
{u}_{k}\label{eq:product_Hu}\\
(\vec{H}
\vec{u}
\vec{u})_i
=&
\sum_{j=1}^N
\sum_{k=1}^N
\vec{H}_{ijk}
{u}_{j}
{u}_{k}\label{eq:product_Huu}
\end{align}
\end{subequations}
Same definition holds for the equivalent products in modal basis.

\red{
\section{Normal form coefficients}
\label{app:NF_coeff}
In this appendix, the nonlinear coefficients of the mapping given by the normal form theory are given in detail. For the sake of brevity, only the coefficients appearing in the simplified case where a single master mode $p$ is selected in the reduced model, are recalled. The interested reader can refer to~\cite{touze03-NNM,TOUZE:JSV:2006} for most complete expressions covering all the cases, including  also damping. We begin with the second-order coefficients, $a^s_{pp}$, $b^s_{pp}$ and $\gamma^s_{pp}$ coefficients, with $p$ the master mode and $s$ a slave mode~:
\begin{subequations}
\begin{align}
&a^s_{pp}=g^s_{pp}\dfrac{2\,\omega^2_p-\omega^2_s}{-\omega^2_s(4\,\omega^2_p-\omega^2_s)},\\
&b^s_{pp}=g^s_{pp}\dfrac{2}{-\omega^2_s(4\,\omega^2_p-\omega^2_s)},\\
&\gamma^s_{pp}=g^s_{pp}\dfrac{2}{4\,\omega^2_p-\omega^2_s}.
\end{align}
\end{subequations}
The third order nonlinear mapping coefficients are equal to zero in the specific case when  $s= p$. Their full expressions for $s\neq p$ read:
\begin{subequations}
\begin{align}
&r^s_{ppp}=
\dfrac{
(A^s_{ppp}+h^s_{ppp})(7 \omega_p^2-\omega_s^2)+
2 B^s_{ppp}(\omega_p^4)
}{(\omega_s^2-\omega_p^2)(\omega_s^2-9\omega_p^2)},\\
&u^s_{ppp}=
\dfrac{
6(A^s_{ppp}+h^s_{ppp})+
B^s_{ppp}(3\omega_p^2-\omega_s^2)
}{(\omega_s^2-\omega_p^2)(\omega_s^2-9\omega_p^2)},\\
&\mu^s_{ppp}=
\dfrac{
6(A^s_{ppp}+h^s_{ppp})+
B^s_{ppp}(3\omega_p^2-\omega_s^2)
}{(\omega_s^2-\omega_p^2)(\omega_s^2-9\omega_p^2)},\\
&\nu^s_{ppp}=
\dfrac{
3(A^s_{ppp}+h^s_{ppp})(3\omega_p^2-\omega_s^2)+
2 B^s_{ppp}(\omega_p^2\omega_s^2)
}{(\omega_s^2-\omega_p^2)(\omega_s^2-9\omega_p^2)},\\
\end{align}
\end{subequations}
where:
\begin{subequations}
\begin{align}
&A^s_{ppp} = \sum_{l=1}^N 2 g^s_{pl} a^l_{pp}\\
&B^s_{ppp} = \sum_{l=1}^N 2 g^s_{pl} b^l_{pp}
\end{align}
\end{subequations}
with $\vec{g}_{pl}$, the vector of quadratic coupling between the master mode $p$ and a generic mode $l$ of the structure. 
}
\section{Linear change of coordinates from physical to modal basis}\label{app:phys_to_mod}

The nonlinear force vector in physical basis reads:
\begin{equation}
\vec{F}(\vec{u})=
\vec{K}\vec{u}+
\vec{G}
\vec{u}
\vec{u}
+
\vec{H}
\vec{u}
\vec{u}
\vec{u}.
\label{eq:F_u}
\end{equation}
The nonlinear force vector in modal basis is in the form:
\begin{equation}
\vec{f}(\vec{X})
=
\vec{\Omega}^2\vec{X}+
\vec{g}
\vec{X}
\vec{X}
+
\vec{h}
\vec{X}
\vec{X}
\vec{X}
,
\label{eq:f_x}
\end{equation}
where the assumption of mass normalised eigenvectors is used to retrieve the squared eigenfrequencies on the diagonal of the matrix $\vec{\Omega}^2$.

The transformation from physical to modal basis uses the full linear eigenvector matrix $\vec{\Phi}$ and reads:
\begin{equation}
\vec{f}(\vec{X})
=
\vec{\Phi}^T\vec{F}(\vec{\Phi}\vec{X}).
\label{eq:f_x_F_PhiX}
\end{equation}
Expanding the right-hand side (RHS) term, it reads:
\begin{equation}
\vec{f}(\vec{X})
=
\vec{\Phi}^T\vec{K}
\;\vec{\Phi}\vec{X}
+
\vec{\Phi}^T\vec{G}
\;\vec{\Phi}\vec{X}\;
\vec{\Phi}\vec{X}
+
\vec{\Phi}^T
\vec{H}\;
\vec{\Phi}\vec{X}\;
\vec{\Phi}\vec{X}\;
\vec{\Phi}\vec{X}
.
\label{eq:f_x_F_PhiX_expanded}
\end{equation}
The relation between linear stiffness matrix in physical and modal coordinates can easily be found by comparing the linear terms in Eq.~\eqref{eq:f_x} and Eq.~\eqref{eq:f_x_F_PhiX_expanded}, allowing one to retrieve the classical formula:
\begin{equation}
\vec{\Omega}^2=\vec{\Phi}^T\vec{K}\vec{\Phi}.
\end{equation}
To relate the quadratic and cubic tensors in physical basis to those in modal basis, it is necessary to expand the term $\vec{\Phi}\vec{X}$ into the sum of all eigenvectors multiplied by their modal amplitudes:
\begin{equation}
\vec{\Phi}\vec{X}=\sum^N_{i=1} \vec{\phi}_i X_i,
\end{equation}
and substitute this sum into Eq.~\eqref{eq:f_x_F_PhiX_expanded}.
By doing so one obtains:
\begin{subequations}
\begin{align}
\vec{\Phi}^T
\vec{G}\;\vec{\Phi}\vec{X}\;
\vec{\Phi}\vec{X}
=&\;
\vec{\Phi}^T
\vec{G}
\;\sum_{i=1}^N \vec{\phi}_i X_i \;
\sum_{j=1}^N \vec{\phi}_j X_j
=
\sum_{i=1}^N \sum_{j=1}^N 
\vec{\Phi}^T
\vec{G}\vec{\phi}_i
\vec{\phi}_j 
\;X_i X_j,
\\
\vec{\Phi}^T\vec{H}\;
\vec{\Phi}\vec{X}\;
\vec{\Phi}\vec{X}\;
\vec{\Phi}\vec{X}=
=&\;
\vec{\Phi}^T
\vec{H}
\;\sum_{i=1}^N \vec{\phi}_i X_i \;
\sum_{j=1}^N \vec{\phi}_j X_j \;
\sum_{k=1}^N \vec{\phi}_k X_k
=
\sum_{i=1}^N \sum_{j=1}^N \sum_{k=1}^N 
\vec{\Phi}^T
\vec{H}\vec{\phi}_i
\vec{\phi}_j\vec{\phi}_k \; X_i X_j X_k
.
\end{align}
\label{eq:a_long_one}
\end{subequations}
where the last simplification comes from the rearrangement of summations made in order to isolate the modal amplitudes. Finally, comparing the RHS of Eqs.~\eqref{eq:a_long_one} with the definition of products given in Eq.~\eqref{eq:product_gXX} and Eq.~\eqref{eq:product_hXXX} of Appendix~\ref{app:products} leads to:
\begin{subequations}
\begin{align}
&\vec{g}_{ij}
=\vec{\Phi}^T
\vec{G}\vec{\phi}_i
\vec{\phi}_j,
\label{eq:g_ij}\\
&\vec{h}_{ijk}
=\vec{\Phi}^T
\vec{H}\vec{\phi}_i \vec{\phi}_j
\vec{\phi}_k.
\label{eq:h_ijk}
\end{align}
\end{subequations}
These last two equations allows expressing the quadratic and cubic coefficients of the modal basis from those computed in the physical basis. Please note that the obtained formula directly depend on the choice of the representation used for the coefficients. Since we have selected to keep full-order tensors of coefficients without exploiting the symmetries arising from the fact that the usual product is commutative, the obtained formula are as in~\eqref{eq:g_ij}-\eqref{eq:h_ijk}. If one chooses to use symmetric tensors for the coefficients, then, for the quadratic term, the relationship would have read  $\vec{g}_{ij}=2\vec{\Phi}^T
\vec{G}\vec{\phi}_i
\vec{\phi}_j$ and $\vec{g}_{ji}=\vec{0}$.

\section{First and second order derivatives of the nonlinear force vector}
\label{app:derivatives}
Given the definition of the nonlinear force tensor in physical basis, one can show that:		
\begin{alignat*}{5}
\left(\dfrac{\partial \vec{F}}{\partial \vec{u}}\right)^r_s=
\dfrac{\partial F^r}{\partial u_s}&=
\text{K}_{rs}&&\;+\;&&
\sum_{j=1}^N \;
({G}^r_{sj}+{G}^r_{js})
{u}_{j}&&\;+\;&&
\sum_{j=1}^N
\sum_{k=1}^N\;
({H}^r_{sjk}+{H}^r_{jsk}+{H}^r_{jks})
{u}_{j}
{u}_{k}\\
&=\text{K}^r_{s}&&\;+\;&&
\sum_{j=1}^N \;
\;2\;{G}^r_{sj}
{u}_{j}&&\;+\;&&
\sum_{j=1}^N
\sum_{k=1}^N\;
\;3\;{H}^r_{sjk}
{u}_{j}
{u}_{k}
\end{alignat*}
where the last simplification is derived from the symmetry of the quadratic and cubic tensors \cite{muravyov}. In compact form we can write:
\begin{equation}
\dfrac{\partial \vec{F}(\vec{u})}{\partial \vec{u}}=\vec{K}+2\vec{Gu}+3\vec{Hu}\vec{u}.
\end{equation}
ans similarly for the second order derivatives:
\begin{alignat*}{4}
\left(\dfrac{\partial^2 \vec{F}}{\partial \vec{u}\partial \vec{u}}\right)^r_{sp}=
\dfrac{\partial F^r}{\partial u_s \partial u_p}
&=&
&({G}^r_{sp}+{G}^r_{ps})&
&\;+\;&
&
\sum_{k=1}^N\;
({H}^r_{spk}+{H}^r_{psk}+{H}^r_{pks}+{H}^r_{skp}+{H}^r_{ksp}+{H}^r_{kps})
{u}_{k}\\
&=&
&\;\;2\;{G}^r_{sp}&
&\;+\;&
&
\sum_{j=1}^N
\sum_{k=1}^N\;
\;\;6\;{H}^r_{spk}
{u}_{k}
\end{alignat*}
and in compact form:
\begin{equation}
\dfrac{\partial^2 \vec{F}(\vec{u})}{\partial \vec{u}\partial \vec{u}}=
2\vec{G}+6\vec{Hu}
\label{d2F_du2}
\end{equation}

\section{Derivation of modal derivatives}
\label{app:taylor_MD}

In this appendix, we derive the Taylor expansion of the nonlinear eigenproblem defined by Eq.~\eqref{nonlinear_eigprob}, and recalled here for the sake of completeness:
\begin{equation}
\left(\dfrac{\partial\vec{F}(\vec{u})}{\partial \vec{u}} - \tilde{\omega}_i^2(\vec{u}) \vec{M}\right)\tilde{\vec{\phi}}_i(\vec{u}) = \vec{0},
\label{nonlinear_eigprob_annex}
\end{equation}
Assuming small motions in the vicinity of the position at rest, this nonlinear eigenproblem can be expanded up to second order as:
\begin{equation}
\begin{split}
&\left(\dfrac{\partial\vec{F}(\vec{u})}{\partial \vec{u}} - \tilde{\omega}_i^2(\vec{u}) \vec{M}\right)\tilde{\vec{\phi}}_i(\vec{u})=\\
&\qquad\left(
\left(\dfrac{\partial\vec{F}(\vec{u})}{\partial \vec{u}} - \tilde{\omega}_i^2(\vec{u}) \vec{M}\right)\tilde{\vec{\phi}}_i(\vec{u})\right)\bigg\rvert_{\vec{u}=\vec{0}}+
\sum^n_{j=1}\left(
\dfrac{\partial}{\partial \q_j}
\left(
\left(\dfrac{\partial\vec{F}(\vec{u})}{\partial \vec{u}} - \tilde{\omega}_i^2(\vec{u}) \vec{M}\right)\tilde{\vec{\phi}}_i(\vec{u})
\right)\right)
\bigg\rvert_{\vec{u}=\vec{0}}\q_j+\mathcal{O}(|\vec{\q}|^2),
\end{split}
\label{Taylor_expansion}
\end{equation}
where the second term has been expanded along the coordinates $R_j$ used for the reduced basis, and the expansion has been written up to $\mathcal{O}(|\vec{\q}|^2)$ terms, or equivalently  $\mathcal{O}(|\vec{u}|^2)$ terms. As the constant term of the expansion, one retrieves the linear eigensystem of Eq.~\eqref{linear_eigsys}, since the Jacobian of the nonlinear force vector at $\vec{u}(\vec{\q}=\vec{0})$ coincides with the linear stiffness matrix $\vec{K}$:
\begin{equation}
\dfrac{\partial\vec{F}(\vec{u})}{\partial \vec{u}}\bigg\rvert_{\vec{u}=\vec{0}}=\vec{K}.
\end{equation}
Consequently the first term of the expansion allows recovering the $i$-th eigenvalue $\omega_i$ as well as the  $i$-th eigenvector $\vec{\phi}_i$.

To verify Eq.~\eqref{nonlinear_eigprob} up to first order, not only the constant term, but also all the linear terms in $\q_j, \forall j=1,\dots,N$ must be zero. By expanding the $j$-th term, one obtains the condition:
\begin{equation}
\left(\left(
\dfrac{\partial}{\partial \q_j}
\left(
\dfrac{\partial \vec{F}(\vec{u})}{\partial \vec{u}}\right)\right)\bigg\rvert_{\vec{0}}
-
\dfrac{\partial \tilde{\omega}_i^2(\vec{\q})}{\partial \q_j}
\bigg\rvert_{\vec{0}}
\vec{M}
\right)
\vec{\phi}_i
+
\left(
\vec{K}
-
\omega_i^2 \vec{M}
\right)
\dfrac{\partial \tilde{\vec{\phi}}_i(\vec{\q})}{\partial \q_j}
\bigg\rvert_{\vec{0}}
=\vec{0}.
\label{linear_terms_in_j}
\end{equation}
In this equation, the sought modal derivatives is the vector $\frac{\partial \tilde{\vec{\phi}}_i(\vec{\q})}{\partial \q_j}\bigg\rvert_{\vec{0}}$ and the other unknown of the system is the value $\frac{\partial \tilde{\omega}_i^2(\vec{\q})}{\partial \q_j}
\bigg\rvert_{\vec{0}}$. 
Moreover, by noticing that:
\begin{equation}
\left(\dfrac{\partial}{\partial \q_j}
\left(
\dfrac{\partial \vec{F}(\vec{u})}{\partial \vec{u}}
\right)\right)\bigg\rvert_{\vec{0}}=
\left(\dfrac{\partial}{\partial \vec{u}}\left(
\dfrac{\partial \vec{F}(\vec{u})}{\partial \vec{u}}\right)
\dfrac{\partial \vec{u}}{\partial \q_j}\right)\bigg\rvert_{\vec{0}}=
\left(
\dfrac{\partial^2 \vec{F}(\vec{u})}{\partial \vec{u}\;\partial \vec{u}}
\right)\bigg\rvert_{\vec{0}}
\vec{\phi}_j,
\end{equation}
one can write Eq.~\eqref{linear_terms_in_j} as:
\begin{equation}
\left(
\dfrac{\partial^2 \vec{F}(\vec{u})}{\partial \vec{u}\;\partial \vec{u}}
\right)\bigg\rvert_{\vec{0}}
\vec{\phi}_j \vec{\phi}_i
-
\dfrac{\partial \tilde{\omega}_i^2(\vec{\q})}{\partial \q_j}
\bigg\rvert_{\vec{0}}
\vec{M}
\vec{\phi}_i
+
\left(\vec{K}-\omega_i^2 \vec{M}\right)
\dfrac{\partial \tilde{\vec{\phi}}_i(\vec{\q})}{\partial \q_j}
\bigg\rvert_{\vec{0}}
=\vec{0}.
\end{equation}
The first term can be further simplified by recalling the definition of the nonlinear force vector and the value of the second derivatives of it given in Eq.~\eqref{d2F_du2}, leading to:
\begin{equation}
2\vec{G}
\vec{\phi}_j \vec{\phi}_i
-
\dfrac{\partial \tilde{\omega}_i^2(\vec{\q})}{\partial \q_j}
\bigg\rvert_{\vec{0}}
\vec{M}
\vec{\phi}_i
+
\left(\vec{K}-\omega_i^2 \vec{M}\right)
\dfrac{\partial \tilde{\vec{\phi}}_i(\vec{\q})}{\partial \q_j}
\bigg\rvert_{\vec{0}}
=\vec{0}.
\label{main_system}
\end{equation}
This is now an undetermined system of equation in the unknowns $\frac{\partial \tilde{\vec{\phi}}_i(\vec{\q})}{\partial \q_j}
\bigg\rvert_{\vec{0}}$ and $\frac{\partial \tilde{\omega}_i^2(\vec{\q})}{\partial \q_j}
\bigg\rvert_{\vec{0}}$. To solve this system, the additional equation of mass normalisation must be introduced.

Following a similar approach, i.e. expanding in Taylor series the nonlinear mass normalisation equation one obtains:
\begin{equation}
\tilde{\vec{\phi}}_i(\vec{\q})^T 
\vec{M}
\tilde{\vec{\phi}}_i(\vec{\q})-1=
\left(\tilde{\vec{\phi}}_i(\vec{\q})^T 
\vec{M}
\tilde{\vec{\phi}}_i(\vec{\q})\right)\bigg\rvert_{\vec{0}}
-1
+
\sum^n_{i=1} 
\left(
\dfrac{\partial}{\partial \q_j}
\left(
\tilde{\vec{\phi}}_i(\vec{\q})^T 
\vec{M}
\tilde{\vec{\phi}}_i(\vec{\q})
\right)\right)
\bigg\rvert_{\vec{0}}
\q_j+
\mathcal{O}(|\vec{\q}|^2)
.
\end{equation}
The constant term is verified by the linear eigenvectors $\vec{\phi}_i$ whereas the linear terms must be equal to zero. The linear term in $\q_j$ becomes the required complement to Eq.~\eqref{main_system}. By expanding the derivatives in $\q_j$, it reads:
\begin{equation}
\dfrac{\partial \tilde{\vec{\phi}}_i(\vec{\q})}{\partial \q_j}
\bigg\rvert_{\vec{0}}^T 
\vec{M}
\vec{\phi}_i+\vec{\phi}_i^T 
\vec{M}
\dfrac{\partial \tilde{\vec{\phi}}_i(\vec{\q})}{\partial \q_j}
\bigg\rvert_{\vec{0}}
=0
\label{mass_normalisation}
\end{equation}
In the usual case of symmetric mass matrix, the LHS only reduces to one of its terms as they are equal. In light of this, the system that permits to evaluate the modal derivatives reads:
\begin{equation}
\begin{bmatrix}
\vec{K}-\omega_i^2 \vec{M}
&\;
-\vec{M}
\vec{\phi}_i
\\[10pt]
-\vec{\phi}_i^T 
\vec{M}
&
\; 0
\end{bmatrix}
\left\lbrace
\begin{matrix}
\frac{\partial \tilde{\vec{\phi}}_i(\vec{\q})}{\partial \q_j}\rvert_{\vec{0}}
\\[10pt]
\frac{\partial \tilde{\omega}_i^2(\vec{\q})}{\partial \q_j}\rvert_{\vec{0}}
\end{matrix}\right\rbrace
=
\left\lbrace
\begin{matrix}
2\vec{G}
\vec{\phi}_j \vec{\phi}_i
\\[15pt]
0
\end{matrix}\right\rbrace
\end{equation}

\section{Derivation of the two-dofs system from von K{\'a}rm{\'a}n beam equations}
\label{app:VK2dofs}

In this appendix, we give the detailed calculation for obtaining all the coefficients of the two-dof model used in section \ref{sec:2dofsFLAT}, from the \vonkar beam model.
As described in Sect.~\ref{sec:2dofsFLATmodel}, we only refer to the coupling between the first flexural mode and the fourth axial mode of a clamped clamped beam. The eigenfunctions and eigenvalues of these two modes can be found solving the linear eigenproblem for respectively flexural and longitudinal vibrations.

Recalling Eqs.~\eqref{eq:VKbeamdim} and focusing only on their linear part, reads:
\begin{subequations}\label{eq:VKbeamdim_lin}
\begin{align}
&  \ddot{W}(X,T)+\dfrac{E I}{\delta S} W^{''''}(X,T)=0,\label{eq:VKbeamdim_lin-a}\\
&  \ddot{U}(X,T) - \dfrac{E}{\delta} U{''}(X,T)=0.\label{eq:VKbeamdim_lin-b}
\end{align}
\end{subequations}
Eq.~\eqref{eq:VKbeamdim_lin-a}, can be solved assuming:
\begin{equation}
{W}(X,T)=Q(T)\Phi(X),
\end{equation}
where $\Phi(X)$ has to respect the clamped clamped boundary conditions, namely $\Phi(0)=\Phi(L)=\Phi^{'}(0)=\Phi^{'}(L)=0$. The first three conditions are respected by the eigenfunction $\Phi(X)$ of arbitrary amplitude $A$:
\begin{equation}
\Phi(X)=
A
(\cos{k X}-\cosh{k X})
(\sin{k L}-\sinh{k L})
-
A
(\sin{k X}-\sinh{k X})
(\cos{k L}-\cosh{k L}),
\end{equation}
whereas the last condition gives rise to the wavelength equation:
\begin{equation}\label{eq:wavelength_flex}
\cos{k L}\cosh{k L}=1.
\end{equation}
The first value of $k>0$ that verifies the transcendental Eq.~\eqref{eq:wavelength_flex} is the dimensional wavelength $k_1$ of the first flexural mode. Its eigenfrequency is then obtained by solving Eq.~\eqref{eq:VKbeamdim_lin-a} and reads:
\begin{equation}
\omega_{1f}^2=\dfrac{E I}{\delta S}k_1^4.
\end{equation}

As for the fourth longitudinal mode, a similar approach is followed. Imposing the separation of variables on $U$:
\begin{equation}
U(X,T) = P(T)\Psi(X),
\end{equation}
and imposing the clamped clamped boundary conditions $\Psi(0)=\Psi(L)=0$ one can find the eigenfunction that has now a simpler form being Eq.~\eqref{eq:VKbeamdim_lin-b} a second order differential equation in $X$. The first boundary condition is verified by:
\begin{equation}
\Psi(X) = B \sin{\kappa X},
\end{equation} 
and the second one by the wavelength equation:
\begin{equation}\label{eq:wavelength_long}
\sin{\kappa L}=0.
\end{equation}
The dimensional wavelength of the fourth axial mode is fourth value of $\kappa>0$ that respects Eq.~\eqref{eq:wavelength_long} equal to $\kappa_4=4 \pi/L$. The eigenfrequency of the fourth axial mode is then obtained from Eq.~\eqref{eq:VKbeamdim_lin-b} and reads:
\begin{equation}
\omega_{4a}^2=\dfrac{E}{\delta}\dfrac{(4 \pi)^2}{L^2}.
\end{equation}

Before operating the reduction that will produce the 2 dofs system of ODEs, it is convenient to make Eq.~\eqref{eq:VKbeamdim_lin} non-dimensional with the following identities:
\begin{align*}
& X = x L,\\
& T = t T_0,\\
& W = w h,\\
& U = u h,\\
& Q = q h,\\
& P = p h,
\end{align*}
and by introducing the additional quantities:
\begin{align*}
& c = \sqrt{E/\delta},\\
& \beta = k_1 L,\\
& \sigma = h/L,\\
& T_0 = 1/\omega_{1f} = \sqrt{12} L^2/(c h \beta^2),
\end{align*}
where the rectangular section assumption has been used in the last equation to simplify $I/S=h^2/12$.
It is now possible to rewrite Eqs.~\eqref{eq:VKbeamdim} with respect to the new variables $x,t,w,u$ as:
\begin{subequations}\label{eq:VKbeamNONdim_ugly}
\begin{align}
& 
\label{eq:VKbeamNONdim_ugly-a}
\dfrac{h}{T_0^2} 
\left(
w_{,tt}
+
\dfrac{1}{\beta^4} w_{,xxxx}
-
\dfrac{12}{\beta^4 \sigma}
\left(u_{,x}w_{,x}\right)_{,x}
-
\dfrac{6}{\beta^4}
\left({w_{,x}}^3\right)_{,x}
\right)
=0,\\
&\label{eq:VKbeamNONdim_ugly-b}
\dfrac{h}{T_0^2} 
\left(
u_{,tt} 
-
\frac{12}{\beta^4 \sigma^2}u_{,xx}
-
\frac{12}{\beta^4 \sigma}
 w_{,x}w_{,xx}
\right)=0.
\end{align}
\end{subequations}
which coincides with Eqs.~\eqref{eq:VKbeamNONdim} multiplied by the nonzero factor $h/T_0^2$.

This system of equations can be reduced to a system of ODEs using the equations for $w$ and $u$ now in their non-dimensional form:
\begin{align}
& w(x,t) =  q_1(t)\Phi_{1f}(x),\\
& u(x,t) =  p_4(t)\Psi_{4a}(x),
\end{align}
with:
\begin{align*}
& \Phi_{1f}(x) = \alpha_1
(\cos{\beta x}-\cosh{\beta x})
(\sin{\beta}-\sinh{\beta})
-
(\sin{\beta x}-\sinh{\beta x})
(\cos{\beta}-\cosh{\beta})\\
& \Psi_{4a}(x) = \alpha_4\sin{4 \pi x},
\end{align*}
by projecting Eq.~\eqref{eq:VKbeamNONdim_ugly-a} on the shapefunction $\Phi_{1f}$ and Eq.~\eqref{eq:VKbeamNONdim_ugly-b} on the shapefunction $\Psi_{4a}$.


The Galerkin projection of Eqs.~\eqref{eq:VKbeamNONdim_ugly} leads to:
\begin{subequations}\label{eq:VKbeam2dofProj}
\begin{align}
& 
\label{eq:VKbeam2dofProj-a}
q_{1,tt}
+
q_1
-
\dfrac{2}{\sigma}G
q_1\,p_4
+D
q_1^3
=0,\\
&\label{eq:VKbeam2dofProj-b}
p_{4,tt} 
+
\frac{12 (4\pi)^2}{\beta^4 \sigma^2}
p_4
-
\dfrac{1}{\sigma}C
q_1^2
=0.
\end{align}
\end{subequations}
With the coefficients $G$, $D$, $C$ being equal to:
\begin{align}
&G=-
\dfrac{6}{\beta^4}
\left(
\dfrac{
\int_0^1\Phi\left(\Psi_{,x}\Phi_{,x}\right)_{,x} dx}
{\int_0^1\Phi^2 dx}
\right)
\\
&C=-\frac{12}{\beta^4}
\left(
\dfrac{
\int_0^1\Psi \left(\Phi_{,x}\Phi_{,xx}\right) dx}
{\int_0^1\Psi^2 dx}
\right)
\\
&D=-
\dfrac{6}{\beta^4}
\left(
\dfrac{
\int_0^1\Phi\left({\Phi_{,x}}^3\right)_{,x} dx}
{\int_0^1\Phi^2 dx}
\right)
\end{align}

If the arbitrary amplitudes $\alpha_1$ and $\alpha_4$ are chosen to have mass normalised eigenfunctions:
\begin{align*}
&\alpha_1: \int_0^1\Phi^2 dx=1,\\
&\alpha_4: \int_0^1\Psi^2 dx=1,
\end{align*}
the quadratic coupling coefficients are symmetric $G=C=1.23$, the cubic coefficient is $D=2.67$ and  Eqs.~\ref{eq:vk2dofsV1} are finally retrieved.
\section*{Conflict of interest}
The authors declare that they have no conflict of interest.
\section*{Data availability statement}
The codes written to run most of the simulations presented in this paper can be available upon simple request to the authors.
\bibliographystyle{spmpsci}
\bibliography{biblio}

\end{document}